\documentclass[11pt]{article}

\usepackage[margin=1.3in]{geometry}

\usepackage{graphicx}
\usepackage{amsmath}
\usepackage{amssymb}
\usepackage{amsthm}
\usepackage{textcomp}
\usepackage{esint}
\usepackage{xcolor}
\usepackage{mathrsfs}
\usepackage{tikz}
\usetikzlibrary{arrows.meta,
                positioning,
                quotes
                }
\usepackage{subfigure}
\usepackage[colorlinks=true]{hyperref}
\hypersetup{
    colorlinks=blue,%
    citecolor=red,%
    filecolor=black,%
    linkcolor=blue,%
    urlcolor=blue
}
\usepackage[colorinlistoftodos]{todonotes}

\newtheorem{theorem}{Theorem}
\newtheorem{remark}{Remark}
\newtheorem{lemma}{Lemma}
\newtheorem{proposition}{Proposition}

\newtheorem{corollary}{Corollary}

\def \no{\nonumber}
\def\e{\epsilon}
\def\ve{\varepsilon}
 
\newcommand{\R}{\mathbb{R}}

\newcommand{\Rn}{\mathbb{R}^n}
\newcommand{\ud}{\mathrm{d}}
\newcommand{\Bn}{\mathbb{B}^{n}}
\newcommand{\Sn}{\mathbb{S}^n}
\newcommand{\Sp}{\mathbb{S}}
\newcommand{\N}{\mathbb{N}}

\newcommand{\pa}{\partial}

\makeatletter

\newdimen\bibspace
\renewenvironment{thebibliography}[1]{%
 \section*{\refname 
       \@mkboth{\MakeUppercase\refname}{\MakeUppercase\refname}}%
     \list{\@biblabel{\@arabic\c@enumiv}}%
          {\settowidth\labelwidth{\@biblabel{#1}}%
           \leftmargin\labelwidth
           \advance\leftmargin\labelsep
           \itemsep\bibspace
           \parsep\z@skip     %
           \@openbib@code
           \usecounter{enumiv}%
           \let\p@enumiv\@empty
           \renewcommand\theenumiv{\@arabic\c@enumiv}}%
     \sloppy\clubpenalty4000\widowpenalty4000%
     \sfcode`\.\@m}
    {\def\@noitemerr
      {\@latex@warning{Empty `thebibliography' environment}}%
     \endlist}

\makeatother

\makeatletter
 \allowdisplaybreaks

\begin{document}

\title{\Large \bf
Almost sharp Sobolev trace inequalities in the unit ball  under constraints}
\author{Xuezhang  Chen\thanks{X. Chen is partially supported by NSFC (No.11771204). Email: xuezhangchen@nju.edu.cn.},~ Wei Wei\thanks{W. Wei: wei\_wei@nju.edu.cn.}~ and Nan Wu\thanks{N. Wu: nwu2@nd.edu.}\\
 \small
$^\ast$$^\dag$Department of Mathematics \& IMS, Nanjing University, Nanjing
210093, P. R. China
}

\date{}

\maketitle

\begin{abstract}
We establish three families of Sobolev trace inequalities of orders two and four in the unit ball under higher order moments constraint, and are able to construct \emph{smooth} test functions to show all such inequalities are \emph{almost optimal}.  Some distinct feature  in \emph{almost sharpness} examples between the fourth order  and second order Sobolev trace inequalities is discovered. This has been neglected in higher order Sobolev inequality case in \cite{Hang}.  As a byproduct, the method of our construction can be used to show the sharpness of the generalized Lebedev-Milin inequality under constraints.

\medskip


\medskip

\end{abstract}

\tableofcontents

\section{Introduction}

The study of optimal constants in the Sobolev trace inequality has a long history. The purpose of this paper is to study three new families of Sobolev trace inequalities in the unit ball, under constraint of higher order moments with respect to the standard volume element on its boundary, and give examples to show these inequalities are \emph{almost optimal}. Such optimal constants are in connection with the \emph{cubature formulas} on spheres. We expect that such almost sharp  inequalities will bring us more interesting applications to geometric problems in the future.

We would like to give a brief survey on the history of Sobolev trace inequalities in the unit ball. For $ n\geqslant 2$ and $m  \in \mathbb{N}$, we define
   \begin{align*}
   	\mathcal{P}_{m}=&\left\{\mathrm{all~polynomials~on~} \mathbb{R}^{n} \mathrm {~with~degree~at~most~} m\right\},\\
   	\mathring{\mathcal{P}}_{m}=&\left\{p\in \mathcal{P}_{m}; \int_{\mathbb{S}^{n-1}} p \ud \mu_{\Sp^{n-1}}=0\right\}.
\end{align*}

For the second order Sobolev trace inequality in the unit disk, Lebedev-Milin \cite{Lebedev-Milin} established the following inequality in $1951$. 
\begin{theorem}[Lebedev-Milin]
For $f \in C^\infty(\Sp^1)$, let $u$ be a smooth extension of $f$ to the unit disk $\mathbb{B}^2$ , then
	\begin{equation}\label{ineq:Lebedev-Milin}
		\log\left(\frac{1}{2\pi} \int_{\Sp^1} e^f \ud \mu_{\Sp^1}\right) \leqslant \frac{1}{4\pi} \int_{\mathbb{B}^2} |\nabla u|^2 \ud x+\frac{1}{2\pi}\int_{\Sp^1} f \ud \mu_{\Sp^1}.
	\end{equation}
\end{theorem}

The generalization of Lebedev-Milin inequality was first proved in $1988$ by Osgood-Phillips-Sarnak \cite{OPS} for the first order moment and later extended by Widom \cite{Widom}  to all higher order moments.
\begin{theorem}[Osgood-Phillips-Sarnak $\&$ Widom]\label{Thm:Widom}
Let $f \in C^\infty(\Sp^1)$ satisfy $\int_{\Sp^1}e^f p \ud \mu_{\Sp^1}=0$ for all $p\in \mathring{\mathcal{P}}_m$, and $u$ be a smooth extension of $f$ to the unit disk $\mathbb{B}^2$ , then
	\begin{equation}\label{ieq:Widom}
		\log\left(\frac{1}{2\pi} \int_{\Sp^1} e^f \ud \mu_{\Sp^1}\right) \leqslant \frac{1}{4\pi(m+1)} \int_{\mathbb{B}^2} |\nabla u|^2 \ud x+\frac{1}{2\pi}\int_{\Sp^1} f \ud \mu_{\Sp^1}.
	\end{equation}
\end{theorem}

In higher dimension $n\geqslant 3$, Beckner \cite{Beckner} proved a sharp Sobolev inequality in the unit ball $\Bn$ with boundary $\Sp^{n-1}$.
\begin{theorem}[Beckner]
Let $n\geqslant 3$ and $1< q \leqslant n/(n-2)$, then for all $u \in H^1 (\mathbb{B}^{n})$, there holds 
   	\begin{align}\label{ineq:Beckner}
   		|\mathbb{S}^{n-1}|^\frac{q-1}{q+1}\left(\int_{\mathbb{S}^{n-1}}|u|^{q+1} \ud \mu_{\Sp^{n-1}}\right)^{\frac{2}{q+1}}\leqslant (q-1) \int_{\mathbb{B}^{n}}|\nabla u|^{2} \ud x+\int_{\mathbb{S}^{n-1}}u^{2} \ud \mu_{\Sp^{n-1}},
   	\end{align} 
	where $|\mathbb{S}^{n-1}|$ denotes the volume of the standard unit sphere $\Sp^{n-1}$.
\end{theorem}

For $n\geqslant 3$ and $1<q \leqslant n/(n-2)$, we define
	$$c(n, q)= \frac{q-1}{|\mathbb{S}^{n-1}|^{\frac{q-1}{q+1}}}$$
	and for all $u \in H^1(\Bn)$ with $u \not\equiv 0$ on $\Sp^{n-1}$,
	$$E_q[u]=\frac{(q-1) \int_{\mathbb{B}^{n}}|\nabla u|^{2} \ud x+\int_{\mathbb{S}^{n-1}}u^{2} \ud \mu_{\Sp^{n-1}}}{\left(\int_{\mathbb{S}^{n-1}}|u|^{q+1} \ud \mu_{\Sp^{n-1}}\right)^{\frac{2}{q+1}}}.$$
   
 In the critical case $q=n/(n-2)$, Escobar \cite{escobar7} employed the Obata type argument to characterize  all positive minimizers of the above functional  $E_q$; Y. Y. Li -Zhu \cite{li-zhu} applied the method of moving spheres to classify all positive critical points of the Euler-Lagrange equation of $E_q$; see also Y. Y. Li-Zhang \cite{Li-Zhang}.  In the subcritical case $q\in (1,n/(n-2))$, Guo-Wang \cite{Guo-Wang} recently proved a Liouville type theorem for the Euler-Lagrange equation of $E_q$.

	Following an argument in Aubin's book \cite[pp.61-63]{Aubin_book}, we can prove the following Sobolev trace inequality under the vanishing first order moment of the boundary volume element (see Chang-Xu-Yang \cite[Inequality (2.5) or Lemma 2.3]{Chang-Xu-Yang}):  Let $n\geqslant 3$ and $1< q \leqslant n/(n-2)$,  for any $0<\ve<1$ and any $u \in H^1(\Bn)$ with 
	$$\int_{\Sp^{n-1}}x_i |u|^\frac{2(n-1)}{n-2} \ud \mu_{\Sp^{n-1}}=0,  \qquad 1 \leqslant i \leqslant n,$$
	then there exists a positive constant $C_\ve$ such that
	\begin{equation}\label{ineq:Aubin}
	\left(\int_{\mathbb{S}^{n-1}}|u|^{q+1} \ud \mu_{\Sp^{n-1}}\right)^{\frac{2}{q+1}} \leqslant c(n,q)\left(2^{\frac{1-q}{q+1}}+\ve\right) \int_{\mathbb{B}^{n}}|\nabla u|^{2} \ud x+C_\ve\int_{\mathbb{S}^{n-1}}u^{2} \ud \mu_{\Sp^{n-1}}.
	\end{equation}
	The above inequality has applications to the prescribed boundary mean curvature problem in $\Bn$. For instance, given $0<h \in C^\infty(\Sp^{n-1})$ with $\max_{\Sp^{n-1}}f/\min_{\Sp^{n-1}}f<2^{1/(2(n-2))}$ and $2\leqslant q< n/(n-2)$, let $u_q$ be a positive smooth minimizer of 
	$$\inf_{u \in \mathscr{S}_h^q}\frac{\int_{\Bn} |\nabla u|^2 \ud x+\frac{n-2}{2}\int_{\Sp^{n-1}}u^2 \ud \mu_{\Sp^{n-1}}}{(\int_{\Sp^{n-1}}h |u|^{q+1} \ud \mu_{\Sp^{n-1}})^{\frac{2}{q+1}}},$$
	where 
	$$\mathscr{S}_h^q=\left\{u \in H^1(\Bn); \int_{\Sp^{n-1}}x |u|^{q+1} \ud \mu_{\Sp^{n-1}}=0 \quad\mathrm{and} \quad \int_{\Sp^{n-1}}h|u|^{q+1} \ud \mu_{\Sp^{n-1}}=1\right\}.$$ 
	The above inequality \eqref{ineq:Aubin} plays a central role in deriving a key estimate, the lower bound of $\|u_q\|_{L^2(\Sp^{n-1})}$, in the subcritical approximations method; see \cite{Chang-Xu-Yang}.

In a celebrated paper of Ache-Chang \cite{Ache-Chang}, authors proved the fourth order sharp Sobolev trace inequalities in $\mathbb{B}^4$ and $\Bn$ for $n\geqslant 5$, which are natural  counterparts of Lebedev-Milin and Beckner inequalities, respectively. For readers' convenience, we restate Ache-Chang sharp Sobolev trace inequalities as follows.

\begin{theorem}[Ache-Chang]\label{Thm:Ache-Chang}
Let  $f\in C^{\infty}(\mathbb{S}^{n-1})$ and $n \geqslant 5$. Suppose $u$ is a smooth extension of $f$ to $\mathbb{B}^n$ which satisfies Neumann boundary condition
\begin{equation*}
	\frac{\pa u}{\pa r} =-\frac{n-4}{2}f \qquad \mathrm{~~on~~} \Sp^{n-1},
\end{equation*}
where $\pa_r$ is the outward unit normal derivative on $\Sp^{n-1}$. Then
\begin{align}\label{ineq:Ache-Chang_n>4}
	&c_n |\mathbb{S}^{n-1}|^{\frac{3}{n-1}} \left(\int_{\mathbb{S}^{n-1}}|f|^{\frac{2(n-1)}{n-4}}\ud\mu_{\mathbb{S}^{n-1}}\right)^{\frac{n-4}{n-1}}\no\\
	\leqslant& \int_{\mathbb{B}^n} \left(\Delta u\right)^2 \ud x + 2\int_{\mathbb{S}^{n-1}}|\nabla f|_{\Sp^{n-1}}^2 \ud\mu_{\Sp^{n-1}} +b_n \int_{\mathbb{S}^{n-1}} f^2 \ud\mu_{\Sp^{n-1}},
\end{align}
where $c_n =n(n-2)(n-4)/4$ and $b_n =n(n-4)/2$. Moreover, equality holds if and only if $u$ is the biharmonic extension of $f_{z_0}(x)=c |1-z_0 \cdot x|^{(4-n)/2}$ on $\Sp^{n-1}$,  also satisfying the above Neumann boundary condition, where $c \in \R,  z_0 \in \Bn$. In particular, if $f=1$, then $u(x)=1+(n-4)(1-|x|^2)/4, x \in \Bn$.
\end{theorem}

\begin{theorem}[Ache-Chang]\label{Thm:Ache-Chang_n=4}
Given  $f\in C^{\infty}(\mathbb{S}^{3})$, let $u$ be a smooth extension of $f$ to  $\mathbb{B}^4$ coupled with zero Neumann boundary condition, that is
\begin{align*}
	\frac{\pa u}{\pa r}=0 \qquad \mathrm{~~on~~} \Sp^3.
\end{align*}
 Then with $\bar f:=\int_{\Sp^3}f \ud \mu_{\Sp^3}/(2\pi^2)$, there holds
\begin{align}\label{ineq:Ache-Chang_n=4}
	\log \left(\frac{1}{2 \pi^{2}} \int_{\Sp^{3}} e^{3(f-\bar{f})} \ud \mu_{\Sp^3}\right)
	 \leqslant \frac{3}{16 \pi^{2}}\left[ \int_{\mathbb{B}^{4}}\left(\Delta u\right)^{2} \ud x+2 \int_{\Sp^{3}}|\nabla f|_{\Sp^3}^{2} \ud \mu_{\Sp^3}\right].
\end{align}
Moreover, equality holds if and only if $u$ is a biharmonic extension of some function $f_{z_0}(x)=-\log|1-z_0\cdot x|+C$ on $\Sp^{3}$, and satisfies zero Neumann boundary condition, where $z_0 \in \mathbb{B}^4, C \in \R$. 
\end{theorem}

In $2019$, Chang-Hang \cite{Chang-Hang} initiated a study on Moser-Trudinger-Onofri inequalities under higher order moments constraint, which are similar improvements of \eqref{ineq:Lebedev-Milin} or \eqref{ieq:Widom}. Subsequently, Hang-Wang \cite{Hang-Wang} made an extension of Sobolev inequality for functions in $W^{1,p}(\Sn)$ with $1<p<n$ under the same constraint.

In this paper, we continue with an effort for Sobolev trace inequalities in the unit ball under higher order moments constraint and establish three families of \emph{almost sharp} Sobolev trace inequalities. 

To continue, we need to set up some notations. For $0 < \theta < 1$, as in \cite{Hang-Wang} we define 
\begin{align*}
\mathcal{M}_{m}^{c}\left(\mathbb{S}^{n-1}\right)=&\bigg\{\nu; ~\nu \mathrm{~is~a~probability~measure~on~} \mathbb{S}^{n-1}\mathrm {~supported~on~ } \\
&\quad~~\{\xi_i; i\in \mathbb{N}\}\mathrm{~s. t.~} \int_{\Sp^{n-1}} p \ud \nu=0,~~ \forall~ p \in \mathring{\mathcal{P}}_{m}\bigg\}
\end{align*}
and
\begin{align*}
 \Theta(m, \theta, n-1)=& \inf \bigg\{\sum_{i} \nu_{i}^{\theta}; ~\nu \in \mathcal{M}_{m}^{c}\left(\mathbb{S}^{n-1}\right) \mathrm {~is~supported~on~}\left\{\xi_{i}\right\} \subset \mathbb{S}^{n-1},\\
 &\qquad\qquad\quad~~ \nu_{i}=\nu\left(\left\{\xi_{i}\right\}\right)\bigg\}.
 \end{align*}
 Indeed, it has been proved by Putterman  \cite[Proposition 3.1]{Putterman} that  $\Theta(m, \theta, n-1)$ can only  be achieved by a Dirac probability measure supported on finitely many points. This directly implies that the infimum for $\Theta(m, \theta, n-1)$ is a minimum, which can follow from \cite[Corollary 3.2]{Putterman}) or the proof of Proposition \ref{prop:two_numbers} below implicitly.
 Moreover, some known exact values of $ \Theta(m, \theta, n-1)$ are $ \Theta(1, \theta, n-1)=2^{1-\theta}$ and $ \Theta(2,\theta,n-1)=(n+1)^{1-\theta}$  by Hang-Wang \cite{Hang-Wang} and $ \Theta(3, \theta, n-1)=(2n)^{1-\theta}$ by Putterman \cite[Theorem 5.1]{Putterman}, whose method is related to the idea of deriving cubature formulas, such as the technique of reproducing kernels on spheres, etc. 

We first state the second order Sobolev trace inequality in higher order moments case, which is a natural generalization of Beckner inequality. 

\begin{theorem}\label{Thm:Improved_Sobolev_Trace_ineq}
	Let $n \geqslant 3, m \in \mathbb{N}$ and $1< q \leqslant n/(n-2)$, then for any $u \in H^1 (\mathbb{B}^{n})$ with
	\begin{align}\label{ineq:2nd order_subcritical}
		\int_{\mathbb{S}^{n-1}} p |u|^{q+1} \ud \mu_{\Sp^{n-1}}=0
	\end{align}
	for all $p \in \mathring{\mathcal{P}}_{m}$, and  for any $\ve >0$ we distinguish it into two cases:
	
	\begin{enumerate}
	\item[(i)] when $1< q < n/(n-2)$, there exists a positive constant $C_\ve$ such that
	\begin{align*}
	\left(\int_{\mathbb{S}^{n-1}}|u|^{q+1} \ud \mu_{\Sp^{n-1}}\right)^{\frac{2}{q+1}}\leqslant\ve \int_{\mathbb{B}^{n}}|\nabla u|^{2} \ud x+C_\ve \int_{\mathbb{S}^{n-1}}u^{2} \ud \mu_{\Sp^{n-1}};
	\end{align*}
	 \item[(ii)] when $q=n/(n-2)$, there exists a positive constant $C_\ve$ such that
	\begin{align}\label{ineq:2nd order}
	&\left(\int_{\mathbb{S}^{n-1}}|u|^{\frac{2(n-1)}{n-2}} \ud \mu_{\Sp^{n-1}}\right)^{\frac{n-2}{n-1}}\no\\
	 \leqslant&\left( \frac{c(n, \frac{n}{n-2})}{\Theta(m, \frac{n-2}{n-1}, n-1)} +\ve \right)\int_{\mathbb{B}^{n}}|\nabla u|^{2} \ud x+C_\ve \int_{\mathbb{S}^{n-1}}u^{2} \ud \mu_{\Sp^{n-1}},
	\end{align}
	where 
	$$c(n, \frac{n}{n-2})= \frac{\frac{2}{n-2}}{|\mathbb{S}^{n-1}|^{\frac{1}{n-1}}}.$$
	\end{enumerate}
\end{theorem}
The proof of Theorem \ref{Thm:Improved_Sobolev_Trace_ineq}, as well as Theorems \ref{Thm1:Ache-Chang-type} and  \ref{Thm:Sobolev trace n=4} below, relies on a modified compactness and concentration argument of Lions. One of our main advances is the use of conic way to connect a Borel measure in the ball with a Borel measure on the sphere for the deduction of concentration compactness principle in this setting; see the proof of Lemma \ref{Concentration compactness bdry} for example.

Moreover, we can show that the number $c(n, \frac{n}{n-2})/\Theta(m, \frac{n-2}{n-1}, n-1)$ in \eqref{ineq:2nd order} is \emph{almost optimal}.

\begin{proposition}\label{example:n>2}
If $n \geqslant 3, m \in \N$ and there exist $a,b \in \R$ such that	
		\begin{align}\label{ineq:almost_optimal}
	\left(\int_{\mathbb{S}^{n-1}}|u|^{\frac{2(n-1)}{n-2}} \ud \mu_{\Sp^{n-1}}\right)^{\frac{n-2}{n-1}} \leqslant a\int_{\mathbb{B}^{n}}|\nabla u|^{2} \ud x+b \int_{\mathbb{S}^{n-1}}u^{2} \ud \mu_{\Sp^{n-1}}
	\end{align}
	for any $u\in H^1(\mathbb{B}^n)$ with 
\begin{align}\label{moment vanishing}
\int_{\mathbb{S}^{n-1}}p |u|^{\frac{2(n-1)}{n-2}}  \ud \mu_{\Sp^{n-1}}=0, \quad  \qquad \forall~ p \in \mathring{\mathcal{P}}_{m},
\end{align}
then 
$$a \geqslant  \frac{c(n, \frac{n}{n-2})}{\Theta(m, \frac{n-2}{n-1}, n-1)}.$$ 
\end{proposition}

\vskip 8pt

Next we transfer to the other two families of fourth order Sobolev trace inequalities under constraints.

\begin{theorem}\label{Thm1:Ache-Chang-type}
	Suppose  $n \geqslant 5, m \in \N$ and for any $\ve >0$ and any $f \in C^\infty (\Sp^{n-1})$ with
	\begin{align*}
		\int_{\mathbb{S}^{n-1}} p |f|^{\frac{2(n-1)}{n-4}} \ud \mu_{\Sp^{n-1}}=0
	\end{align*}
	for all $p \in \mathring{\mathcal{P}}_{m}$, then there exists a positive constant $C_\ve$ such that
	\begin{align}\label{ineq:n>5}
	&\left(\int_{\mathbb{S}^{n-1}}|f|^{\frac{2(n-1)}{n-4}} \ud \mu_{\Sp^{n-1}}\right)^{\frac{n-4}{n-1}} \no \\\leqslant&\left( \frac{\alpha(n)}{\Theta(m, \frac{n-4}{n-1}, n-1)} +\ve \right)\int_{\mathbb{B}^{n}}\left(\Delta u\right)^{2} \ud x+C_\ve\int_{\mathbb{S}^{n-1}} \left( |\nabla f|_{\Sp^{n-1}}^2  +b_n  f^2 \right)\ud\mu_{\Sp^{n-1}},	\end{align}
	where $u$ is an $H^2(\Bn)$ norm extension of $f$ and satisfies the Neumann boundary condition
\begin{align}\label{Neumann_bdry_cond_n>4}
	 \frac{\pa u}{\pa r} =-\frac{n-4}{2}f \qquad \mathrm{~~on~~} \Sp^{n-1} 
\end{align}
and 
	$$ b_n =\frac{n(n-4)}{2}, \quad \alpha(n)= \frac{4}{n(n-2)(n-4) |\mathbb{S}^{n-1}|^{\frac{3}{n-1}} }.$$
\end{theorem} 

For non-vanishing Neumann (precisely, Robin) boundary condition \eqref{Neumann_bdry_cond_n>4}, our conic proof of Theorem \ref{Thm1:Ache-Chang-type} has some advantage and thus sounds interesting to readers.

As in \cite{Chang-Hang} we define
\begin{align*}
&\mathcal{N}_{m}(\Sp^{n-1})   \\
=&\bigg\{ N\in \mathbb{N}; \exists ~x_{1},\cdots ,x_{N}\in \mathbb{S}^{n-1}%
\mathrm{~~and~~}\nu 
_{1},\cdots ,\nu _{N}\in \left[ 0,\infty \right) \\
&\qquad \qquad  \mathrm{~~with~~
}\sum_{i=1}^N \nu _i=1 \mathrm{~~and~~}%
\sum_{i=1}^N \nu _{i}p( x_{i})=0,~~\forall~
p\in \mathring{\mathcal{P}}_{m}\bigg\} .
\end{align*}
The smallest number in $\mathcal{N}_{m}(\Sp^{n-1})$ is denoted as $N_{m}(\Sp^{n-1})$, i.e. $%
N_{m}(\Sp^{n-1})=\min \mathcal{N}_{m}(\Sp^{n-1})$. Moreover, Chang-Hang gave an elementary proof in \cite{Chang-Hang} to show that $N_m(\Sp^1)=m+1$ for all $m \in \N$ and $N_1(\Sp^2)=2, N_2(\Sp^2)=4$. Later Hang \cite{Hang} extended to prove $N_1(\Sp^{n-1})=2$ and $N_2(\Sp^{n-1})=n+1$ for all $n\geqslant 2$. Lower bounds of $N_m(\Sp^{n-1})$ can be found in Dai-Xu's book \cite[Theorems 6.1.2 and 6.1.4]{Dai-Xu} (see also Delsarte-Goethals-Seidel \cite{DGS} in $1977$):
\begin{equation}\label{lbd:DGS}
\begin{split}
N_m(\Sp^{n-1})\geqslant C_{n-1+\frac{m}{2}}^{n-1}+C_{n+\frac{m}{2}-2}^{n-1}  \qquad&\mathrm{if~} m \mathrm{~is~even}; \\
N_m(\Sp^{n-1})\geqslant 2C_{n-1+\frac{m-1}{2}}^{n-1}  \qquad &\mathrm{if~} m \mathrm{~is~odd}.
\end{split}
\end{equation}
Besides \cite{Hang}, some examples meeting the above lower bounds, called ``\emph{tight $m$-spherical designs}" in  \cite{DGS}, can help us  to derive other exact values of $N_m(\Sp^{n-1})$;  see also  \cite{Bannai1,Bannai2,Bannai3}. In particular, $N_3(\Sp^{n-1})=2n$ by virtue of Lemma \ref{lem:N_3}. The exact value of  $N_m(\Sp^{n-1})$ is intimately related to  \emph{cubature formulas} on spheres; see \cite[Chapter 6]{Dai-Xu}, \cite{Mysovskih} etc.


\begin{theorem}\label{Thm:Sobolev trace n=4}
	Given  $f\in C^{\infty}(\mathbb{S}^{3})$, suppose $u$ is a smooth extension of $f$ to $\mathbb{B}^4$ which also satisfies zero Neumann boundary condition.  
 Assume that $f$ satisfies  $
\int_{\mathbb{S}^{3}}pe^{3f}\ud\mu_{\Sp^3} =0$ for all $p\in \mathring{\mathcal{P}}_{m}$  with some $m \in \N$, then for any $\varepsilon >0$, there exists a positive constant $C_\ve$ such that
\begin{align}\label{ineq:n=4}
	&\log \left(\frac{1}{2 \pi^{2}} \int_{\Sp^{3}} e^{3(f-\bar f)} \ud \mu_{\Sp^3}\right)\no\\
	 \leqslant& \left(\frac{3}{16 \pi^{2}N_m(\Sp^3)}+\ve \right)\left[\int_{\mathbb{B}^{4}}\left(\Delta u\right)^{2} \ud x+2 \int_{\Sp^{3}}|\nabla f|_{\Sp^3}^{2} \ud \mu_{\Sp^3}\right]+C_\ve,
	\end{align}
	where $\bar f=\int_{\Sp^3} f \ud \mu_{\Sp^3}/(2\pi^2)$.
	\end{theorem}

Since there exist several more challenging obstructions in addition to the second order case, \emph{almost sharpness} examples for the fourth order Sobolev trace inequalities become more subtle and interesting. At the same time, without the second order example, the ones for fourth order case are impossible to appear.

\vskip 8pt

\begin{proposition}\label{example:n>4}
Suppose $n \geqslant 5, m \in \N$  and there exist $a,b \in \R$ such that 
\begin{align*}\label{Fourth Trace inequality}
	&\left(\int_{\mathbb{S}^{n-1}}|f|^{\frac{2(n-1)}{n-4}}\ud\mu_{\mathbb{S}^{n-1}}\right)^{\frac{n-4}{n-1}}\no\\
	\leqslant& a\int_{\mathbb{B}^n} \left(\Delta u\right)^2 \ud x + b \int_{\mathbb{S}^{n-1}}\left(|\nabla f|_{\Sp^{n-1}}^2 \ud\mu_{\Sp^{n-1}}+b_n  f^2\right) \ud\mu_{\Sp^{n-1}} 
\end{align*}
for any $u \in H^2(\Bn)$ satisfying 
\begin{align*}
		\int_{\mathbb{S}^{n-1}} p |f|^{\frac{2(n-1)}{n-4}} \ud \mu_{\Sp^{n-1}}=0 \qquad \forall~p \in \mathring{\mathcal{P}}_{m},
\end{align*}
coupled with boundary conditions: 
\begin{equation*}
	u=f, \quad \frac{\pa u}{\pa r} =-\frac{n-4}{2}f \qquad \mathrm{~~on~~} \Sp^{n-1}.
\end{equation*}
Then 
$$a \geqslant  \frac{\alpha(n)}{\Theta(m, \frac{n-4}{n-1}, n-1)}$$
with
$$\alpha(n)= \frac{4}{n(n-2)(n-4) |\mathbb{S}^{n-1}|^{\frac{3}{n-1}} }.$$
\end{proposition}

For the construction of examples in Propositions \ref{example:n>2} and \ref{example:n>4}, the following fact is enough to our use: The number $\Theta(m, \theta, n-1)$ is achieved by some $\nu=\sum_{i=1}^{N}\nu_i \delta_{x_i} \in \mathcal{M}_{m}^{c}\left(\mathbb{S}^{n-1}\right)$ for some $N\geqslant N_m(\Sp^{n-1})$.

\begin{proposition}\label{example:n=4}
Suppose there exist $a,b \in \R$ such that for all $u \in H^2(\mathbb{B}^4)$ satisfying boundary conditions:
	$$u=f, \quad \frac{\pa u}{\pa r}=0 \qquad \mathrm{on~~} \Sp^3,$$
	and $
\int_{\mathbb{S}^{3}}pe^{3f}\ud\mu_{\Sp^3} =0$ for all $p\in \mathring{\mathcal{P}}_{m}$ with $m \in  \mathbb{N}$,  we have
\begin{align*}
	\log \left(\frac{1}{2 \pi^{2}} \int_{\Sp^{3}} e^{3(f-\bar f)} \ud \mu_{\Sp^3}\right) \leqslant a\left[\int_{\mathbb{B}^{4}}\left(\Delta u\right)^{2} \ud x+2 \int_{\Sp^{3}}|\nabla f|_{\Sp^3}^{2} \ud \mu_{\Sp^3}\right]+b,
	\end{align*}
	where $\bar f=(2\pi^2)^{-1}\int_{\mathbb{S}^3} f\ud \mu_{\Sp^3}$.
	Then
	$$a \geqslant\frac{3}{16 \pi^{2}N_m(\Sp^3)}.$$
\end{proposition}

As have shown before, optimal constants are related to two numbers $\Theta(m, \theta, n-1)$ and $N_m(\Sp^{n-1})$. Furthermore, it is of independent interest that a natural relationship between these two numbers has been discovered.
 \begin{proposition}\label{prop:two_numbers}
For $n\geqslant 2$ and all $m \in  \mathbb{N}$, there holds
$$\lim_{\theta \searrow  0} \Theta(m, \theta, n-1)=N_m(\Sp^{n-1}).$$
\end{proposition}

Compared with Chang-Hang \cite{Chang-Hang} and Hang-Wang \cite{Hang-Wang}, some additional difficulties arise from compact manifolds with boundary and high order of conformally covariant operators.  The following are the novelties of our constructions. First, we introduce a union of finitely many conic annuli $\cup_{i=1}^{N}\mathcal{A}_\delta(x_i)$ (see \eqref{def:Conic_Annulus} for the definition  and Figure \ref{fig:Conic Annuli} (a) for example) to isolate the dominated terms, which contribute to the sharp constant, and its complement in $\Bn$ controls higher order terms with delicate computations and observations. Second, the Neumann boundary conditions are of geometric favor, arising from conformally covariant boundary operator $P_3^g$, \emph{GJMS operator} of order three. As the involved operator is fourth order, precisely, the \emph{Paneitz operator}, we need to make an appropriate correction to the test function like the second order one such that the new test function satisfies Neumann boundary condition.   Third, since the extremal metric in Ache-Chang's sharp Sobolev trace inequalities of order four in \cite{Ache-Chang} is not the flat metric, but ``\emph{the adapted metric}",  first introduced by Case-Chang \cite{Case-Chang}, this forces us to know the exact \emph{geometric local bubbles}, which are resolved in Section \ref{Sect2}. In dimension $n\geqslant 5$, a further correction is still needed to the  \emph{geometric  local  bubble} to ensure that \emph{the new local bubble} satisfies the Neumann boundary condition and controls higher order terms.  By the way, the local bubble for the second order case is well-known in the study of boundary Yamabe problem. Geometric intuitions play an important role in all constructions. Finally,  in contrast to the example for $n\geqslant 5$, we soon realize that  the Chang-Hang type estimate  \cite[(3.12)]{Chang-Hang}\footnote{The estimate (3.12) is also crucial in the  Chang-Hang's example. At first glance, it sounds very tough to be improved. Fortunately, a clever observation on $\mathring{\mathcal{P}_1}$ in Proposition \ref{prop:3-dim_m=1} motivates us to achieve this goal for all $m\in \N$ in Proposition \ref{prop:n-dim_N_2}.} is not enough in dimension four. Whenever we struggled in this optimal constant, a geometric intuition/Branson's intuitive proof always indicates that such an example should be there. That is exactly our motivation for improving the Chang-Hang type estimate. Such an essential improvement  rescues us from above dilemma. We eventually achieve this goal in Section \ref{Subsect:improved_estimate}. This also demonstrates one of main differences of higher order Sobolev trace inequality from the second order one.\footnote{After the completion of our examples, we look back to the example for higher order Sobolev inequality in  \cite[p.21]{Hang}. We do not know how to construct the example as the author advised in \cite{Hang} without our new type estimate as in Proposition \ref{prop:n-dim_N_2}.}

To the best knowledge of authors, \emph{almost sharp} examples in this paper seem to be \emph{the first ones}  of fourth order Sobolev trace inequalities.

To demonstrate the relationships among these Sobolev trace inequalities in the unit ball, we draw a diagram for readers' convenience.
 \vskip 8pt
 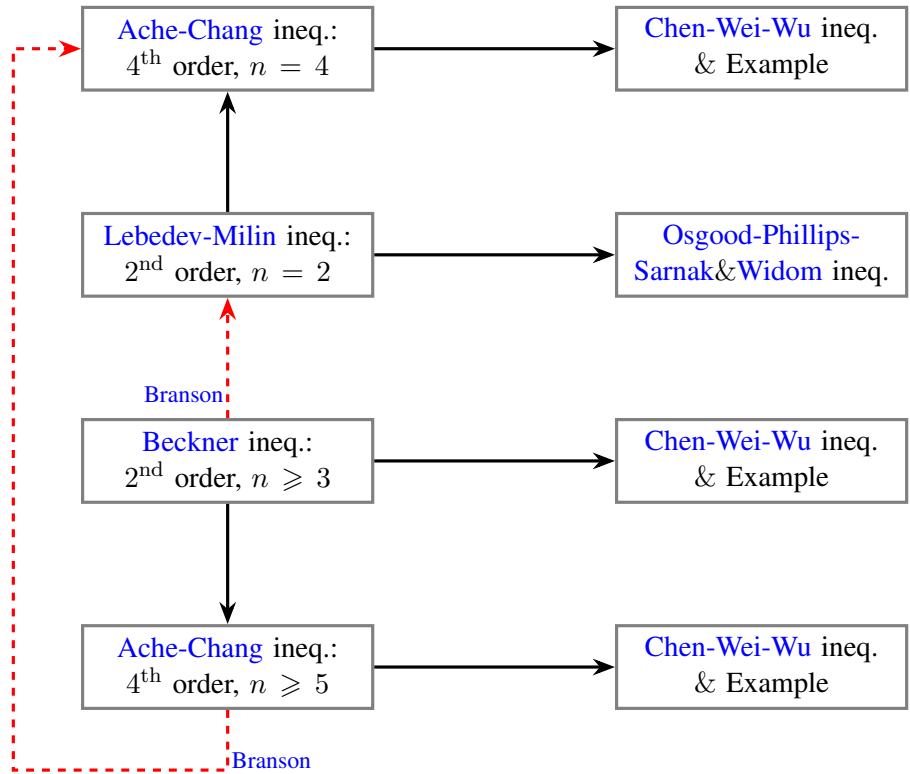
\begin{figure}[!h]
    \centering
        \begin{tikzpicture}[auto,
                       > = Stealth, 
           node distance = 16mm and 32mm,
              box/.style = {draw=gray, very thick,
                            minimum height=11mm, text width=36mm, 
                            align=center},
       every edge/.style = {draw, ->, very thick},
every edge quotes/.style = {font=\footnotesize, align=center, inner sep=1pt}
                            ]
    \node (n11) [box]               {\textcolor{blue}{Ache-Chang} ineq.: $4^{\mathrm{th}}$ order, $n \geqslant 5$};
    \node (n12) [box, right=of n11] { \textcolor{blue}{Chen-Wei-Wu} ineq. \\ $\&$ Example};
    \node (n21) [box, above=of n11] { \textcolor{blue}{Beckner} ineq.: $2^{\mathrm{nd}}$ order, $n\geqslant 3$};
    \node (n22) [box, above=of n12] {\textcolor{blue}{Chen-Wei-Wu} ineq. \\ $\&$ Example};
    \node (n31) [box, above=of n21] {  \textcolor{blue}{Lebedev-Milin} ineq.: $2^{\mathrm{nd}}$ order, $n=2$};
    \node (n32) [box, above=of n22] {\textcolor{blue}{Osgood-Phillips-Sarnak}$\&$\textcolor{blue}{Widom} ineq.};
    \node (n41) [box, above=of n31] {\textcolor{blue}{Ache-Chang} ineq.: $4^{\mathrm{th}}$ order, $n=4$};
    \node (n42) [box, above=of n32] {\textcolor{blue}{Chen-Wei-Wu} ineq. \\ $\&$ Example};

\draw  
        (n21) edge [" "]         (n11)
        (n21) edge [" "]         (n22)
        (n21) edge [pos=0.2, red, dashed, "\textcolor{blue}{Branson}"]     (n31)
        (n31) edge [" "]         (n32)
        (n31) edge [" "]         (n41)
        (n41) edge [" "]         (n42)
        (n11) edge [" "]          (n12);
\path[draw=red,dashed, very thick, ->]
    (n11.south) -+ (0,-11mm) to ["\textcolor{blue}{Branson}"] 
    ([yshift=-8mm] n11.south) -|   ([xshift=-9mm] n41.west) -- (n41);
        \end{tikzpicture}
    \caption{\textcolor{red}{Sharpness} vs \textcolor{red}{Almost Sharpness}}
    \label{fig:Sobolev trace inequalities}
    \end{figure}

Here, ``\textcolor{blue}{Branson}" means Thomas P. Branson's intuitive proof: \emph{dimension continuity argument} and the \textcolor{blue}{Widom} (called \textcolor{blue}{Osgood-Phillips-Sarnak} for $m=1$) inequality is also \textcolor{red}{sharp}.
\vskip 8pt

We outline a unified approach to our constructions for fourth order Sobolev trace inequalities in three steps. \emph{Step 1}, a localized analysis in the conic annulus $\mathcal{A}_\delta(x_i)$ together with ``\emph{a local bubble}" $\phi_{\ve,i}$ is used to find the restriction of the test function to the boundary, say, $v$, satisfying higher order moments constraint. \emph{Step 2}, find a natural way to extend $v$ to a global function $u$ in $\Bn$ for $n\geqslant 4$, satisfying the Neumann boundary condition in the following ways.
\begin{itemize}
\item For $n\geqslant 5$, we denote the test function by
$$u^{\frac{2(n-1)}{n-4}}=v^{\frac{2(n-1)}{n-4}}+(1-r)\tilde g(x), \qquad r=|x|, x\in \Bn,$$
where $1-r$ is exactly the boundary defining function. Again the Neumann boundary condition \eqref{Neumann_bdry_cond_n>4} enables us to obtain the exact expression of $\tilde g$.
\item For $n=4$, with delicate selections of suitable cut-off functions and the `\emph{nice}' properties of the ``\emph{local bubble}" $\phi_{\ve,i}$, we find the scheme used in the second order case still valid in this case except for an improved Chang-Hang type estimate.
\end{itemize}
\emph{Step 3}, delicate calculations and  deep insights are employed to capture optimal constants and to control higher order terms. We shall convince geometric intuitions through an analytic way.

\vskip 8pt

The following is the organization of this paper. In Section \ref{Sect2}, we revisit Ache-Chang's elegant proof   to understand ``\emph{local bubbles}" and their `\emph{nice}' geometric properties in the viewpoint of conformal geometry. In Section \ref{Sect3}, we prove the second order Sobolev trace inequality and give an example of precise test functions to show the \emph{almost sharpness}. A proof of Proposition \ref{prop:two_numbers} is also presented there. In Section \ref{subsect:4th order ineq}, we adapt our conic proof to the fourth order Sobolev trace inequality for $n\geqslant 5$.
In Section \ref{Sect:n>4}, we complete the construction of our example for $n\geqslant 5$ and thus finish the proof of Proposition \ref{example:n>4}. Section \ref{Sect5} is devoted to the fourth order Sobolev trace inequality established in Section \ref{Subsect:4th_Sobolev_trace_ineq} and the four dimensional example. In Section \ref{Subsect:improved_estimate}, we give an elementary proof of the exact value of $N_3(\Sp^{n-1})$ and establish an improvement of the Chang-Hang type estimate, which enables us to construct an example to complete the proof of Proposition \ref{example:n=4} in Section \ref{Subsect:n=4}. In Appendix \ref{Append}, we employ an example to show the \emph{sharpness} of Widom inequality as a warm-up of the four dimensional example, which is a subtle case as we know. 

\section{Geometric interpretations of local bubbles}\label{Sect2}

Before presenting our concrete constructions of test functions, we think it important to describe geometric ideas behind them, as it will be a long journey.

As we have pointed out before, the local bubble in the second order is originated from the study of boundary Yamabe problem. Precisely, the bubble function is the conformal factor of a conformal metric in the class of the flat metric in upper half-space, which is scalar-flat with positive constant boundary mean curvature. So, we only focus on the ones associated to the fourth order Sobolev trace inequalities.

The \emph{Paneitz operator} on a smooth Riemannian manifold $(M,g)$ of dimension $n \geqslant 3$ is defined by
$$P_4^g=\Delta_g^2+\delta_g (4A_g-(n-2)J_g g)(\nabla \cdot,\cdot)+\frac{n-4}{2}Q_4^g,$$
where $\delta_g$ is the divergent operator, $A_g=\frac{1}{n-2}[\mathrm{Ric}_g-\frac{R_g}{2(n-1)}g]$, $R_g$ is the scalar curvature, $J_g=\mathrm{tr}_g(A_g)$ and $Q_4^g$ is the $Q$-curvature
$$Q_4^g=-\Delta_g J_g+\frac{n}{4}J_g^2-2|A_g|_g^2.$$
It is well-known that $P_4^g$ is conformally covariant in the sense that if we let $\tilde g=u^{4/(n-4)}g \in [g]$ for $n\geqslant 5$, then for all $\psi \in C^\infty(\overline{M})$,
$$u^{\frac{n+4}{n-4}}P_4^{\tilde g}\psi=P_4^{g} (u\psi)=\Delta (u\psi)$$
and if we let $\tilde g=e^{2u}g$ for $n=4$, then $P_4^{\tilde g}=e^{-4u}P_4^{g}$ and
$$P_4^g u+Q_4^g=Q_4^{\tilde g} e^{4u}.$$

Since our construction is of not only analytic but also geometric favors, it sounds important to grasp some key ideas in Ache-Chang's elegant proof. It will enable us to discover  \emph{``a local bubble"}  naturally associated to \emph{``the adapted metric"} in $\Bn$: 
$$g^\ast:=e^{2\tau(x)}|\ud x|^2 \qquad\quad \mathrm{with~~}\quad \tau(x)=\frac{1-|x|^2}{2}\quad\mathrm{~~and~~} n = 4 $$ 
and
$$g^\ast:=\psi(x)^{\frac{4}{n-4}}|\ud x|^2 \qquad\quad \mathrm{with~~}\quad \psi(x)=1+\frac{n-4}{2}\tau(x)\quad\mathrm{~~and~~} n \geqslant 5.$$
 These correspond to the extremal metrics in Ache-Chang's Sobolev trace inequality of order four in \cite[Theorem B]{Ache-Chang} and \cite[Theorem A]{Ache-Chang}, respectively. Moreover, \emph{the adapted metric} $g^\ast$ satisfies the following properties (cf. a particular case of \cite[Proposition 2.1]{Ache-Chang}, and first shown by Case-Chang \cite{Case-Chang} on a Poincar\'e-Einstein manfiold):
\begin{itemize}
\item $Q_4^{g^\ast}=0$;
\item $\Sp^{n-1}=\pa \Bn$ is totally geodesic with respect to $g^\ast$ and $g^\ast=g_{\Sp^{n-1}}$ on $\Sp^{n-1}$;
\item Let $P_3=P_3^{\Sp^{n-1}}=(B-1)B(B+1)$ be the fractional GJMS operator defined in \cite[(4.1)]{Ache-Chang} with 
 $$B=\sqrt{-\Delta_{\Sp^{n-1}}+\frac{(n-2)^2}{4}}$$
 and its associated $Q$-curvature 
 $$Q_3^{\Sp^{n-1}}:=\frac{2}{n-4}P_3(1)=\frac{n(n-2)}{4}.$$
 Then
\begin{align*}
\frac{1}{2}E_4(g^\ast)[U_f]=&\int_{\Sp^{n-1}}f P_3 f \ud \sigma_{g^\ast}-\frac{n-4}{2}\int_{\Sp^{n-1}}Q_3^{g^\ast} f^2 \ud \sigma_{g^\ast}\\
=&\int_{\Sp^{n-1}}f P_3 f \ud \mu_{\Sp^{n-1}}-\frac{n(n-2)(n-4)}{8}\int_{\Sp^{n-1}}f^2 \ud \mu_{\Sp^{n-1}}.
\end{align*}
\end{itemize}

  We give a brief summary on Ache-Chang's proof: Authors started with the following energy identity
$$0=\int_{\Bn}U_f P_4^{g^\ast} U_f \ud \mu_{g^\ast}=E_4(g^\ast)[U_f]+J(g^\ast)[U_f]$$
for all $U_f$ satisfying
\begin{align}\label{PDE:U_f}
\begin{cases}
\displaystyle P_4^{g^\ast} U_f=0 &\qquad \mathrm{~~in~~} \Bn;\\
\displaystyle U_f=f &\qquad \mathrm{~~on~~} \Sp^{n-1};\\
\displaystyle \frac{\pa U_f}{\pa \nu_{g^\ast}}=0 &\qquad \mathrm{~~on~~} \Sp^{n-1},
\end{cases}
\end{align}
where $\nu_{g^\ast}$ is the outward unit normal with respect to the metric $g^\ast$ on $\Sp^{n-1}$. Here
$$E_4(g^\ast)[U_f]=\int_{\Bn}\left[(\Delta_{g^\ast} U_f)^2-\left(4A_{g^\ast}-(n-2)J_{g^\ast}g^\ast\right)(\nabla U_f,\nabla U_f)\right]\ud \mu_{g^\ast}$$
and $J(g^\ast)[U_f]$ is the boundary terms arising from the integration by parts.

It follows from conformal invariance property of $P_4^{g^\ast}$ that if $n\geqslant 5$, then for all $u \in C^\infty(\overline{\Bn})$,
$$\psi^{\frac{n+4}{n-4}}P_4^{g^\ast}u=P_4^{|\ud x|^2} (u\psi)=\Delta^2 (u\psi)$$
and if $n=4$, then $P_4^{g^\ast}=e^{-4\tau}P_4^{|\ud x|^2}=e^{-4\tau}\Delta^2$. This together with \eqref{PDE:U_f} indicates
\begin{itemize}
\item If $n\geqslant 5$, then
\begin{align}\label{PDE:psi U_f_n>4}
\begin{cases}
\displaystyle \Delta^2(\psi U_f)=0 &\qquad \mathrm{~~in~~} \Bn;\\
\displaystyle \psi U_f=f &\qquad \mathrm{~~on~~} \Sp^{n-1};\\
\displaystyle \frac{\pa (\psi U_f)}{\pa r}=-\frac{n-4}{2}f &\qquad \mathrm{~~on~~} \Sp^{n-1}.
\end{cases}
\end{align}
\item If $n=4$, then
\begin{align}\label{PDE:U_f_n=4}
\begin{cases}
\displaystyle  \Delta^2 U_f=0 &\qquad \mathrm{~~in~~} \Bn;\\
\displaystyle  U_f=f &\qquad \mathrm{~~on~~} \Sp^{n-1};\\
\displaystyle  \frac{\pa U_f}{\pa r}=0 &\qquad \mathrm{~~on~~} \Sp^{n-1}.
\end{cases}
\end{align}
\end{itemize}
Direct consequences of \eqref{PDE:psi U_f_n>4} and \eqref{PDE:U_f_n=4} are
\begin{align*}
0=&\int_{\Bn}|\Delta (\psi U_f)|^2 \ud x+\int_{\Sp^{n-1}}f \frac{\pa}{\pa r}\Delta(\psi U_f) \ud \mu_{\Sp^{n-1}}-\int_{\Sp^{n-1}}\frac{\pa (\psi U_f)}{\pa r} \Delta(\psi U_f) \ud \mu_{\Sp^{n-1}}\\
:=& \int_{\Bn}|\Delta (\psi U_f)|^2 \ud x+J(g_0)[\psi U_f]
\end{align*}
for $n\geqslant 5$; and 
\begin{align*}
0=&\int_{\Bn}(\Delta U_f)^2 \ud x+\int_{\Sp^{n-1}}f \frac{\pa}{\pa r}\Delta U_f \ud \mu_{\Sp^{n-1}}:=\int_{\Bn}(\Delta U_f)^2 \ud x+J(g_0)[U_f]
\end{align*}
for $n=4$.

One of key ingredients in Ache-Chang's proof is to calculate the exact expression of boundary term $J(g^\ast)[U_f]$ at the cost of lengthy computations. Finally, a bridge to Ache-Chang's sharp Sobolev trace inequalities is the following Beckner's sharp Sobolev inequalities (see \cite{Beckner} or \cite[Theorem 4.2]{Ache-Chang}):
\begin{itemize}
\item If $n\geqslant 5$, then for all $f \in C^\infty(\Sp^{n-1})$,
$$\frac{n(n-2)(n-4)}{8}\left(\frac{1}{|\Sp^{n-1}|}\int_{\Sp^{n-1}}|f|^{\frac{2(n-1)}{n-4}}\ud \mu_{\Sp^{n-1}}\right)^{\frac{n-4}{n-1}}\leqslant \frac{1}{|\Sp^{n-1}|}\int_{\Sp^{n-1}}f P_3 f \ud \mu_{\Sp^{n-1}}.$$
\item If $n=4$, then for all $f \in C^\infty(\Sp^3)$,
$$\log \left(\frac{1}{2\pi^2}\int_{\Sp^3}e^{3(f-\bar f)} \ud \mu_{\Sp^3}\right)\leqslant \frac{3}{8\pi^2}\int_{\Sp^3}f P_3 f \ud \mu_{\Sp^3},$$
where $\bar f=(2\pi^2)^{-1}\int_{\Sp^3} f\ud \mu_{\Sp^3}$.
\end{itemize}
With the help of differences of $J(g^\ast)[U_f]-J(g_0)[\psi U_f]$ for $n\geqslant 5$ and $J(g^\ast)[U_f]-J(g_0)[U_f]$ for $n=4$, respectively, authors obtained the desired assertions, see \cite[Theorems A and B]{Ache-Chang}. 

\vskip 8pt

From the viewpoint of conformal geometry, we collect some elementary facts together, which stimulate us to find aforementioned  \emph{``local bubbles"}.
\begin{itemize}
	\item For $n\geqslant 2$, let $F: (\Bn,|\ud x|^2) \to (\Rn_+, |\ud z|^2)$ be an inversion with respect to the sphere $\pa\Bn_{\sqrt{2}} (-e_n)$ with radius $\sqrt 2$ and center at $-e_n$, that is, 
		$$z=F(x)=-e_n+\frac{2(x+e_n)}{|x+e_n|^2},$$
	then $F$ is a conformal map with property that
	$$F_\ast (|\ud x|^2)=\left(\frac{2}{(1+z_n)^2+|z'|^2}\right)^2 |\ud z|^2:=U(z)^{\frac{n-4}{2}}|\ud z|^2,\quad z=(z',z_n)\in \Rn_+.$$
	\item Notice that
		$$F_\ast(\psi^{\frac{4}{n-4}}|\ud z|^2)=\hat \psi(z)^{\frac{4}{n-4}}|\ud z|^2$$
	with
	$$\hat \psi(z)=\left[1+\frac{n-4}{2}\frac{2z_n}{(1+z_n)^2+|z'|^2}\right]U(z).$$
	For $\ve>0$, by scalings we define
	$$\hat \psi_\ve(z)=\ve^{\frac{4-n}{2}}\hat \psi(\frac{z}{\ve}).$$
	\item Combining Case \cite[Corollary 5.2]{Case} and Ache-Chang \cite[Theorem 4.2 (b)]{Ache-Chang}, we can  immediately obtain the following sharp Sobolev trace inequality $H^2(\Rn_+)\hookrightarrow H^{3/2}(\R^{n-1})\hookrightarrow L^{2(n-1)/(n-4)}(\R^{n-1})$ : For all $u \in C_c^\infty(\overline{\Rn_+})$ with $\pa_{z_n}u=0$ on $\R^{n-1}$, there holds
	$$\frac{n(n-2)(n-4)}{4}|\Sp^{n-1}|^{\frac{3}{n-1}} \left( \int_{\R^{n-1} }|u|^{\frac{2(n-1)}{n-4}} \ud z'\right)^{\frac{n-4}{n-1}}\leqslant \int_{\Rn_+}(\Delta u)^2 \ud z.$$
	Moreover, when equality holds, the extremal metric is exactly $\hat \psi^{4/(n-2)} |\ud z|^2$. In particular, 
	 $\hat \psi$ satisfies
	\begin{align*}
	\begin{cases}
	\displaystyle \Delta^2 \hat \psi=0 &\qquad \mathrm{in~~} \Rn_+,\\
	\displaystyle \hat \psi=f &\qquad \mathrm{on~~} \R^{n-1},\\
	\displaystyle \pa_{z_n} \hat \psi=0 &\qquad \mathrm{on~~} \R^{n-1},
	\end{cases}
	\end{align*}
	where 
	$$f(z')=\left(\frac{2}{1+|z'|^2}\right)^{\frac{n-4}{2}}.$$
	As a byproduct, we obtain
	$$\int_{\Rn_+}(\Delta \hat \psi)^2 \ud z=\frac{n(n-2)(n-4)}{4}|\Sp^{n-1}|.$$
	In this specific case, we give an explicit  extremal functions in \cite[Corollary 5.2]{Case}.
	\item This arrives at the construction of the aforementioned \emph{``local bubble"} in $\Bn$. Based on successful experience of the second order example, we replace $z_n$ by $1-r$ and $|z'|$ by $\rho$ in $\psi_\ve=\hat \psi_\ve(z_n,|z'|)$, which is called ``\emph{a local geometric bubble}". However, in order to satisfy the Neumann boundary condition, with a further modification on $\psi_\ve(r,\rho):=\hat \psi_\ve(1-r,\rho)$ we define
		\begin{align*}
	\phi_\ve(r,\rho)=& \psi(r) \psi_\ve(r,\rho)\\
	=&\psi(r)\left[1+\frac{n-4}{2}\frac{2\ve(1-r)}{(\ve+1-r)^2+\rho^2}\right]\left(\frac{2}{(\ve+1-r)^2+\rho^2}\right)^{\frac{n-4}{2}}
	\end{align*}
	for the flat metric
	$$|\ud x|^2=\ud r^2 +r^2 \left(\ud\rho^2 + \sin^2\rho g_{\mathbb{S}^{n-2}}\right)$$
for $x=r\xi \in \Bn$ and $\rho(\xi)=\ud_{\mathbb{S}^{n-1}}(\xi, x_i)=\stackrel{\frown}{\xi x_i}$. For convenience, we call $\phi_\ve$ \emph{``a local bubble"}, which will save us a lot of  energy on calculations.

Then it is not hard to verify that 
\begin{equation}\label{Neumann_bdry_local_bubble:n>4}
\pa_r \phi_\ve= \hat \psi_\ve \pa_r \psi+\psi \pa_r \hat \psi_\e=-\frac{n-4}{2}\hat \psi_\ve=-\frac{n-4}{2}\phi_\ve \qquad \mathrm{on~~} \Sp^{n-1}.
\end{equation}
The above property is fundamental in the positivity of our test function $u$ in $\overline\Bn$, as well as the control of higher order terms.

\end{itemize}

As above, we can apply the same scheme to find the  ``\emph{local geometric bubble}" in dimension $n=4$ as follows.
\begin{itemize}
\item Notice that
$$F_\ast (|\ud x|^2)=\left(\frac{2}{(1+z_n)^2+|z'|^2}\right)^2 |\ud z|^2=e^{2U_1(z)}|\ud z|^2$$
and then
\begin{equation}\label{adapted_metric_n=4}
F_\ast(g^\ast)=F_\ast(e^{2\tau(x)} |\ud z|^2)=e^{2\tau\circ F^{-1}(z)}e^{2U_1(z)} |\ud z|^2:=e^{2\hat \tau(z)}|\ud z|^2.
\end{equation}
This implies 
$$U_1(z)=\log\frac{2}{(1+z_n)^2+|z'|^2}$$
and
$$\hat \tau(z)=\tau\circ F^{-1}(z)+U_1(z)=\frac{2z_n}{(1+z_n)^2+|z'|^2}+\log\frac{2}{(1+z_n)^2+|z'|^2}$$
for $z=(z',z_n)\in \Rn_+$.

Since
$$Q_4^{g^\ast}=0 \quad \Longrightarrow \quad Q_4^{F_\ast(g^\ast)}=Q_4^{g^\ast}\circ F^{-1}=0,$$
the $Q_4$-curvature equation together with \eqref{adapted_metric_n=4} yields
$$P_4^{|\ud z|^2}\hat \tau=Q_4^{F_\ast(g^\ast)} e^{n\hat \tau}=0 \quad \Longrightarrow \quad \Delta^2 \hat \tau=0 \qquad \mathrm{in~~} \Rn_+.$$
Thus, we obtain
\begin{align*}
	\begin{cases}
	\displaystyle \Delta^2 \hat \tau=0 &\qquad \mathrm{in~~} \Rn_+,\\
	\displaystyle \hat \tau=f_1 &\qquad \mathrm{on~~} \R^{n-1},\\
	\displaystyle \pa_{z_n} \hat \tau=0 &\qquad \mathrm{on~~} \R^{n-1},
	\end{cases}
	\end{align*}
	where 
	$$f_1(z')=\log \frac{2}{1+|z'|^2}.$$

For $\ve>0$, by scalings we define
\begin{align*}
	\hat \tau_\ve(z)=&\tau(\frac{z}{\ve})-\log \ve\\
	=&\frac{2\ve z_n}{(\ve+z_n)^2+|z'|^2}+\log\frac{2\ve}{(\ve+z_n)^2+|z'|^2}.
	\end{align*}

	\item For each $\mathcal{A}_\delta(x_i)\subset \Bn$, under the local coordinates near $x_i$, the flat metric can be expressed as
	$$|\ud x|^2=\ud r^2 +r^2 \left(\ud\rho^2 + \sin^2\rho g_{\mathbb{S}^{2}}\right)$$
for $x=r\xi \in \mathbb{B}^n$ and $\rho(\xi)=\ud_{\mathbb{S}^{n-1}}(\xi, x_i)=\stackrel{\frown}{\xi x_i}$.
	
	As before, replacing $z_n$ and $|z'|$ in $\hat \rho(z)=\hat \rho(z_n, |z'|)$ by $1-r$ and $\rho$, respectively, up to a constant $\log(2\ve)$, we define \emph{``a local geometric bubble"} in dimension four  by
	$$\phi_\ve(r,\rho)=-\log\left((\ve+1-r)^2+\rho^2\right)+\frac{2\ve (1-r)}{(\ve+1-r)^2+\rho^2}.$$
	
	A direct computation yields
	\begin{equation}\label{bubble:Neumann_bdry_cond}
	\pa_r \phi_\ve=0\qquad \mathrm{on~~} \Sp^{n-1}.
	\end{equation}
	This implies that $\phi_\ve$ satisfies zero Neumann boundary condition, which is crucial to control higher order terms.
	\end{itemize}
	
	\vskip 8pt
	
	With the above geometric intuitions at hand, we are confident that our strategy is feasible. A confirmation unavoidably involves  lengthy calculations and deep insights.

\section{Almost sharp Sobolev trace inequality of order two under constraints in dimension three and higher}\label{Sect3}

 The second order example is a basis for the ones of the fourth order case. The successful experience on the second order example opens the doors to the construction of the fourth order examples. It is our first time to introduce the conic proof and the conic annulus $\mathcal{A}_\delta(x_i)$ in the process of our example.

\subsection{Second order Sobolev trace inequality}
We first generalize the concentration and compactness principle in $\Rn$ or closed manifolds (e.g., Lions \cite{Lions1,Lions2}, Struwe \cite[Section 4.8]{Struwe_book} and Chang-Hang \cite{Chang-Hang} etc.) to the unit ball.

 \begin{lemma}\label{Concentration compactness bdry}
	Let $n\geqslant 3$ and $1< q \leqslant n/(n-2)$. Suppose that as $k \to \infty$,  $u_k \rightharpoonup u$ weakly in $H^1(\mathbb{B}^n)$ and
	\begin{align*}
		|\nabla u_{k}|^{2} \ud x \rightharpoonup&~~ |\nabla u|^{2} \ud x+\lambda \qquad\quad~~~\mathrm{as ~~measure~~in~~} \mathbb{B}^n,\\
		\left|u_{k}\right|^{q+1} \ud \mu_{\Sp^{n-1}} \rightharpoonup& ~~|u|^{q+1} \ud \mu_{\Sp^{n-1}}+\nu \qquad\mathrm{as~~measure~~in~~} \mathbb{S}^{n-1},
	\end{align*}
	where $\lambda$ and $\nu$ are bounded nonnegative Borel measures on $\Bn$ and $\Sp^{n-1}$, respectively. For any Borel set $\Omega \subset \Sp^{n-1}$, we define a nonnegative Borel measure on $\Sp^{n-1}$ associated to $\lambda$ through 
	$$\hat{\lambda}(\Omega):= \lambda(\mathcal{C}(\Omega)),$$
	where $\mathcal{C}(\Omega)$ is a cone with vertex at the origin and base $\Omega$, i.e.,
	\begin{align}\label{def:conic_Borel_set}
		\mathcal{C}(\Omega)= \left\{(r, \xi)\in \mathbb{B}^n; 0 \leqslant r \leqslant 1, \xi \in \Omega \right\}.
	\end{align} 
	Then there exist countably many points $x_i \in \mathbb{S}^{n-1}$ with $\nu_i=\nu(\{x_i\})>0$ such that 
	\begin{align*}
		\nu=\sum_{i} \nu_{i} \delta_{x_{i}}, \qquad \hat{\lambda} \geqslant \frac{1}{c(n, q)}\sum_{i}\nu_i^{\frac{2}{q+1}} \delta_{x_i} \qquad \mathrm{~~for~~}q=\frac{n}{n-2}; 
	\end{align*}
	and
	\begin{align*}
		\nu =0, \qquad \qquad \mathrm{~~for~~} 1< q < \frac{n}{n-2}. 
	\end{align*}
\end{lemma}
\begin{proof}
	Let $v_k:= u_k -u$, then up to a subsequence, as $k \to \infty$, $v_k \rightharpoonup 0$ weakly in $H^1(\mathbb{B}^n)$ and $v_k \to 0$ in $L^2 (\mathbb{B}^n)$,  and $v_k \to 0$ in $L^{s+1} (\mathbb{S}^{n-1})$ for any $0<s<n/(n-2)$; also $\nu_{k}:=|u_k|^{q+1} \ud\mu_{\Sp^{n-1}} -|u|^{q+1} \ud\mu_{\Sp^{n-1}} \rightharpoonup \nu$ and $\lambda_{k}:=\left|\nabla v_{k}\right|^{2} \ud x \rightharpoonup \lambda$ in the weak sense of measures.
	
		\vskip 8pt
		
	If $1< q < n/(n-2)$, then the compact embedding of Sobolev trace inequality  from $H^1(\Bn)$ to $L^{q+1}(\Sp^{n-1})$ forces $\nu=0$.
	
	\vskip 8pt
	
	In the following, it suffices to consider $q= n/(n-2)$.
	
	For any $\varphi \in C^1(\overline{\mathbb{B}^n})$, by the sharp Sobolev trace inequality in \eqref{ineq:Beckner} we have 
	\begin{align*}
\int_{\mathbb{S}^{n-1}}|\varphi|^{q+1} \ud \nu =&\lim _{k \rightarrow \infty} \int_{\mathbb{S}^{n-1}}|\varphi|^{q+1} \ud \nu_{k}=\lim _{k \rightarrow \infty} \int_{\mathbb{S}^{n-1}}\left|v_{k} \varphi\right|^{q+1} \ud x \\
 \leqslant& c(n, q)^{\frac{q+1}{2}} \lim _{k \rightarrow \infty}\left[\int_{\mathbb{B}^{n}}\left|\nabla \left(v_{k} \varphi\right)\right|^{2}+( v_k\varphi)^2 \ud x\right]^{\frac{q+1}{2}} \\
=& c(n, q)^{\frac{q+1}{2}} \lim_{k \rightarrow \infty}\left(\int_{\mathbb{B}^{n}}\varphi^{2}\left|\nabla v_{k}\right|^{2} \ud x\right)^{\frac{q+1}{2}} \\
=& c(n, q)^{\frac{q+1}{2}}\left(\int_{\mathbb{B}^{n}}\varphi^{2} \ud \lambda\right)^{\frac{q+1}{2}},
	\end{align*}
	that is,
	\begin{align}\label{ineq:trace_measures}
		\int_{\mathbb{S}^{n-1}}|\varphi|^{q+1} \ud \nu \leqslant c(n, q)^{\frac{q+1}{2}}\left(\int_{\mathbb{B}^{n}}\varphi^{2} \ud \lambda\right)^{\frac{q+1}{2}}
	\end{align}
	for all $\varphi \in C^1(\overline{\mathbb{B}^n})$. 
	
	Now let $J \subset \Sp^{n-1}$ be the set of the atoms of the measure $\nu$. From the assumption and the Sobolev trace inequality of $H^1(\Bn) \hookrightarrow L^{q+1}(\Sp^{n-1})$ that
	 $$\int_{\mathbb{S}^{n-1}}\ud\nu\leqslant \lim_{k\to \infty}\int_{\Sp^{n-1}}|u_k|^{q+1} \ud \mu_{\Sp^{n-1}} \leqslant C\liminf_{k\to \infty} \|u_k\|_{H^1(\Bn)}^{q+1}< \infty,$$
	 we know that $J$ is an  at most countable set denoted by $J=\{x_i; i \in \N\}$.  We decompose 
	 $$\nu = \nu_0 + \sum_{i}\nu_i \delta_{x_i},$$
where $\nu_0 $ is the singular continuous part of the measure $\nu$ and also a nonnegative Borel measure. 

We first claim that $\nu_0=0$. 

To this end, for any open set $\mathcal{O} \subset \Sp^{n-1}$, choosing a sequence $\{\varphi_k\} \subset C^1 (\overline{\mathbb{B}^n})$  in \eqref{ineq:trace_measures}  such that $\varphi_k$ converges to the characteristic function of $\mathcal C(\mathcal{O})$ and letting $k \to \infty$, we obtain
	 \begin{align*}
	 	\int_{\mathcal{O}} \ud\nu \leqslant c(n, q)^{\frac{q+1}{2}}\left(\int_{\mathcal C(\mathcal{O})} \ud \lambda\right)^{\frac{q+1}{2}} =c(n, q)^{\frac{q+1}{2}}\left(\int_{\mathcal{O}} \ud \hat \lambda\right)^{\frac{q+1}{2}} <\infty.
	 \end{align*}
	This obviously implies that $\nu$ is  absolutely continuous with respect to $\hat{\lambda}$ and thus from the Radon-Nikodym theorem that there exists a nonnegative $f \in L^{1}(\Sp^{n-1} , \hat{\lambda})$ such that $\ud \nu=f \ud \hat{\lambda}$.  Moreover, for $\hat{\lambda}$-almost every $x\in \mathbb{S}^{n-1}$ and  $B_{\rho}(x) \subset \mathbb{S}^{n-1}$, we have 
	 \begin{align*}
	 	f(x)=\lim _{\rho \searrow 0} \frac{\int_{B_{\rho}(x)} \ud \nu}{\int_{B_{\rho}(x)} \ud \hat{\lambda}}.
	 \end{align*}
	 In particular, if  we choose $x \in \Sp^{n-1}$ such that the segment $\overline{ox}$ does not carry any atom of $\lambda$, then 
	 \begin{align*}
	 	 f(x)^{\frac{2}{q+1}}=\lim _{\rho \searrow 0}\frac{\left(\int_{B_{\rho}(x)} \ud \nu\right)^{\frac{2}{q+1}}}{\left(\int_{B_{\rho}(x)} \ud \hat{\lambda}\right)^{\frac{2}{q+1}}} \leqslant c(n, q)\lim _{\rho \searrow 0}\left(\int_{\mathcal C(B_{\rho}(x))} \ud \lambda\right)^{\frac{q-1}{q+1}}=0.
	 \end{align*}
Furthermore, $\lambda$ has only countably many atoms, so does $\hat{\lambda}$ by definition. Hence, we conclude that for $\hat{\lambda}$-almost everywhere on $\Sp^{n-1}$, $f(x)=0$. This together with the fact that $\nu_0$ contains no atoms implies that $\nu_0=0$.
	
	Next we go to establish the desired inequality.

	For each $x_i \in J$, we choose $\varphi \in C^1(\overline \Bn)$ such that $0\leqslant \varphi \leqslant 1$ and $\varphi=1$ on the segment $\overline{ox_i}$. Hence,  with $\hat \lambda_i=\hat \lambda(\{x_i\})$, we can apply \eqref{ineq:trace_measures} with the above $\varphi$ to conclude that 
 \begin{align*}
		c(n, q) \hat{\lambda}_i\geqslant  \nu_i^{\frac{2}{q+1}} 
	\end{align*}
	and thus
	\begin{align*}
		\hat{\lambda} \geqslant \frac{1}{c(n, q)}\sum_{i}\nu_i^{\frac{2}{q+1}} \delta_{x_i}.
	\end{align*}
	This completes the proof.
\end{proof}

We are now in a position to prove Theorem \ref{Thm:Improved_Sobolev_Trace_ineq}.

\begin{proof}[Proof of Theorem \ref{Thm:Improved_Sobolev_Trace_ineq}] We shall prove these results by contradiction. If either \eqref{ineq:2nd order_subcritical} or \eqref{ineq:2nd order} is not true, then we can find a sequence $\{u_k\} \subset H^1(\mathbb{B}^n)$ and some $\alpha>0$ for $1<q<n/(n-2)$ or 
	$$\alpha =  \frac{c(n, \frac{n}{n-2})}{\Theta(m, \frac{n-2}{n-1}, n-1)} +\ve \qquad \mathrm{~~for~~}q=\frac{n}{n-2},$$
	respectively,
	such that
	\begin{align*} 
		\left(\int_{\mathbb{S}^{n-1}}|u_k|^{q+1} \ud \mu_{\Sp^{n-1}}\right)^{\frac{2}{q+1}} >\alpha\int_{\mathbb{B}^{n}}|\nabla u_k|^{2} \ud x+k \int_{\mathbb{S}^{n-1}}|u_k|^{2} \ud \mu_{\Sp^{n-1}}
	\end{align*} 
	and
	\begin{align*}
		\int_{\mathbb{S}^{n-1}} p |u_k|^{q+1} \ud \mu_{\Sp^{n-1}}=0, \qquad \forall~ p \in \mathring{\mathcal{P}}_{m}.
	\end{align*}
	We may normalize $\int_{\mathbb{S}^{n-1}}|u_k|^{q+1} \ud \mu_{\Sp^{n-1}}=1$ such that
	\begin{align*}
		\int_{\mathbb{B}^{n}}|\nabla u_k|^{2} \ud x \leqslant \frac{1}{\alpha}
	\end{align*}
	and
	\begin{align*}
		\int_{\mathbb{S}^{n-1}}|u_k|^{2} \ud \mu_{\Sp^{n-1}}\leqslant \frac{1}{k}.
	\end{align*}
	
	The Sobolev inequality
	\begin{align*}
		Y(\Bn,\pa \Bn) \|u_k\|_{L^{\frac{2n}{n-2}}(\mathbb{B}^n)}^2\leqslant \int_{\mathbb{B}^n} |\nabla u_k|^2 \ud x +2(n-1)\int_{\mathbb{S}^{n-1}}u_k^2 \ud \mu_{\Sp^{n-1}}
	\end{align*}
	indicates that $\|u_k\|_{L^{2n/(n-2)}(\mathbb{B}^n)}$ is uniformly bounded, so is $\|u_k\|_{L^{2}(\mathbb{B}^n)}$. Readers are referred to \cite{escobar4} for the definition of the Yamabe constant $Y( \Bn,\pa \Bn)$, which is positive. Thus,  $u_k$ is uniformly bounded  in $H^1(\mathbb{B}^n)$ and then up to a subsequence, $u_k  \rightharpoonup 0$ weakly in $H^1(\mathbb{B}^n)$ as $k \to \infty$.
	
	\vskip 8pt
	
	(i) If $1< q < n/(n-2)$, then $u_k \to 0$ in $L^{q+1}(\mathbb{S}^{n-1})$ as $k \to \infty$. A contradiction! 
	
	\vskip 8pt
	
	(ii) For $q=n/(n-2)$, there exist nonnegative Borel measures $\lambda$ and $\nu$ on the $\sigma$-algebras of $\Bn$ and $\Sp^{n-1}$, respectively, such that, up to a subsequence, as $k \to \infty$
	\begin{align*}
		|\nabla u_k(x)|^{2} \ud x \rightharpoonup  \lambda \qquad\quad~~~\mathrm{as ~~measures~~in~~} \mathbb{B}^n
	\end{align*}
	and 
	\begin{align*}
	|u_k|^{q+1} \ud \mu_{\Sp^{n-1}}   \rightharpoonup \nu \qquad \qquad\mathrm{as~~measures~~in~~} \mathbb{S}^{n-1}.
	\end{align*}
	
Let $\hat \lambda$ be the Borel measure on the $\sigma$-algebra of $\Sp^{n-1}$ associated to $\lambda$ as in Lemma \ref{Concentration compactness bdry}. Again by Lemma \ref{Concentration compactness bdry}, we can find countably many points $x_j \in \mathbb{S}^{n-1}$ such that 
\begin{align*}
	\nu=\sum_{j} \nu_{j} \delta_{x_{j}}
\end{align*}
with $\nu_{j}=\nu\left(\left\{x_{j}\right\}\right)$ and
\begin{align*}
	\hat{\lambda} \geqslant \frac{1}{c(n, \frac{n}{n-2})}\sum_{j}\nu_j^{\frac{n-2}{n-1}} \delta_{x_j}.
\end{align*}

Notice that
\begin{align*}
	\nu\left(\mathbb{S}^{n-1}\right)=1 \qquad \mathrm{~~and~~}  \qquad \hat{\lambda} \left(\mathbb{S}^{n-1}\right)=\lambda \left(\mathbb{B}^{n}\right) \leqslant \frac{1}{\alpha}.
\end{align*}
By definition of the weak convergence for measures, we know that
\begin{align*}
	\int_{\mathbb{S}^{n-1}} p \ud \nu=0 \qquad \mathrm{~~for~~} p \in \mathring{\mathcal{P}}_{m}.
\end{align*}

By definition of $\Theta(m, (n-2)/(n-1), n)$, with $\hat \lambda_{j}=\hat \lambda\left(\left\{x_{j}\right\}\right)$ we have 
\begin{align*}
	\Theta(m, \frac{n-2}{n-1}, n-1) \leqslant \sum_{j} \nu_{j}^{\frac{n-2}{n-1}} \leqslant& \sum_{j} c(n, \frac{n}{n-2})\hat \lambda_{j} \\
	\leqslant& c(n, \frac{n}{n-2}) \hat\lambda\left(\mathbb{S}^{n-1}\right) \leqslant \frac{c(n, \frac{n}{n-2})}{\alpha}.
\end{align*}
Hence, 
\begin{align*}
	\alpha \leqslant \frac{c(n, \frac{n}{n-2})}{\Theta(m, \frac{n-2}{n-1}, n-1)}.
\end{align*}
However, this contradicts the choice of $\alpha$.
\end{proof}

\subsection{Almost sharpness}

We start with a brief discussion on the number $\Theta(m, \theta, n-1)$. Using the idea of proof in Putterman \cite[Proposition 3.1]{Putterman}, one  can prove that if $\nu=\sum_{i=1}^{N}\nu_i \delta_{x_i}\in \mathcal{M}_{m}^{c}\left(\mathbb{S}^{n-1}\right)$ for any $N>\bar N:=\mathrm{dim}(\mathring{\mathcal{P}}_m)$ , then $\nu$ can not be an extremal element of $\Theta(m, \theta, n-1)$. A direct consequence is that the infimum in $\Theta(m, \theta, n-1)$ is a minimum by virtue of \cite[Corollary 3.2]{Putterman}. 

Some exact values have been known. 
\begin{itemize}
\item For  $m=1$, the \cite[Proposition 3.1]{Hang-Wang} states that
$$\Theta(1,\theta,n-1)=2^{1-\theta}$$
is achieved by $\nu_1 \in \mathcal{M}_{1}^{c}\left(\mathbb{S}^{n-1}\right)$ if and only if $\nu_1=\frac{1}{2}(\delta_\xi+\delta_{-\xi})$ for any $\xi \in \Sp^{n-1}$.
\item For $m=2$, the \cite[Proposition 3.2]{Hang-Wang} states that 
	$$\Theta(2, \theta, n-1)=(n+1)^{1-\theta}$$
	is achieved by $\nu_2\in \mathcal{M}_{2}^{c}\left(\mathbb{S}^{n-1}\right)$ if and only if $\nu_2= (\sum_{i=1}^{n+1} \delta_{x_{i}})/(n+1) \in \mathcal{M}_{2}^{c}\left(\mathbb{S}^{n-1}\right)$, where  $x_{1}, \cdots, x_{n+1} \in \mathbb{S}^{n-1}$ are the vertices of a regular $(n+1)$-simplex embedded in $\Bn$.
\item For $m=3$, the \cite[Theorem 1.2]{Putterman} states that 
$$\Theta(3,\theta,n-1)=(2n)^{1-\theta}$$
is achieved by $\nu_3 \in \mathcal{M}_{3}^{c}\left(\mathbb{S}^{n-1}\right)$ if and only if 
$$\nu_3=\frac{1}{2n}\sum_{i=1}^n \left(\delta_{e_i}+\delta_{-e_i}\right),$$
up to an isometry on $\Sp^{n-1}$, where $\{e_i; 1 \leqslant i \leqslant n\}$ is  the standard basis in $\Rn$.
\end{itemize}

  We would like to point out that a combination of the cubature formulas and examples meeting the lower bounds \eqref{lbd:DGS} in \cite{DGS} is helpful to know more exact values of $\Theta(m, \theta, n-1)$. 

Based on the above known results, it is expected that as $\theta \searrow 0$, the limit of $\Theta(m, \theta, n-1)$ should be $N_m(\Sp^{n-1})$. 
\begin{proof}[Proof of Proposition \protect{\ref{prop:two_numbers}}]
First, we show that
$$\limsup_{\theta\searrow 0}\Theta(m, \theta, n-1)\leqslant N_m(\Sp^{n-1}):=N_m.$$

To this end, it follows from the definition of $N_m(\Sp^{n-1})$ that there exist $\{\xi_i; 1 \leqslant i \leqslant N_m\} \subset \Sp^{n-1}$ and $\nu_i>0, \sum_{i=1}^{N_m}\nu_i=1$, such that $\nu=\sum_{i=1}^{N_m}\nu_i \delta_{\xi_i}\in \mathcal{M}_3^c(\Sp^{n-1})$. By definition of $\Theta(m, \theta, n-1)$, we have
\begin{align*}
\Theta(m, \theta, n-1)\leqslant \sum_{i=1}^{N_m} \nu_i^\theta.
\end{align*}
Letting $\theta \searrow 0$ in the above inequality, the desired assertion follows.

Next, it suffices to show
$$ \underline{\Theta}:=\liminf_{\theta\searrow 0}\Theta(m, \theta, n-1)\geqslant N_m.$$

To that end, it follows from \cite[Proposition 3.1 and Corollary 3.2]{Putterman} that $\forall~\theta \in (0,1)$, there exist $\xi_i \in \Sn$ and $\nu_i \geqslant 0$ for $1 \leqslant i \leqslant \bar N$, such that $\sum_{i=1}^{\bar N}\nu_i=1, \sum_{i=1}^{\bar N}\nu_i p_j(\xi_i)=0$ for a basis $\{p_j\} \subset \mathring{\mathcal{P}}_m$ and
$$\Theta(m, \theta, n-1)=\sum_{i=1}^{\bar N}\nu_i^\theta.$$

We can find a sequence of real numbers $\{\theta_k\}$ such that $\theta_k \to 0$ and
$$\Theta(m, \theta_k, n-1) \to \underline{\Theta} \qquad \mathrm{as~~} k \to \infty.$$
For any fixed $\theta_k$, there exist  $\xi_i^{(k)} \in \Sn$ and $\nu_i^{(k)} \geqslant 0$ for $1 \leqslant i \leqslant \bar N$, such that 
\begin{equation}\label{cubature_cond}
\sum_{i=1}^{\bar N}\nu_i^{(k)}=1,\quad \sum_{i=1}^{\bar N}\nu_i^{(k)} p_j(\xi_i^{(k)})=0 \qquad\mathrm{for~all~~} p_j \in \mathring{\mathcal{P}}_m
\end{equation}
and
$$\Theta(m, \theta_k, n-1)=\sum_{i=1}^{\bar N}\left(\nu_i^{(k)}\right)^{\theta_k}.$$
Then up to a subsequence, there exist $\nu_i^{(\infty)}\geqslant 0$ and  $\xi_i^{(\infty)}\in \Sp^{n-1}, 1 \leqslant i \leqslant \bar N$, such that for each $i$,
$$\nu_i^{(k)} \to \nu_i^{(\infty)}, \quad  \xi_i^{(k)}\to \xi_i^{(\infty)}\qquad \mathrm{as~~} k \to \infty.$$
Letting $k \to \infty$ in \eqref{cubature_cond} we obtain
$$\sum_{i=1}^{\bar N}\nu_i^{(\infty)}=1,\quad \sum_{i=1}^{\bar N}\nu_i^{(\infty)} p_j(\xi_i^{(\infty)})=0 \qquad\mathrm{for~all~~} p_j \in \mathring{\mathcal{P}}_m,$$
whence
$$\nu^{(\infty)}:=\sum_{i=1}^{\bar N}\nu_i^{(\infty)}\delta_{\xi_i^{(\infty)}}\in\mathcal{M}_m^c(\Sp^{n-1}).$$
It follows from the definition of $N_m(\Sp^{n-1})$ that 
$$N_m \leqslant \sharp\left\{\nu_i^{(\infty)}>0,~~1 \leqslant i \leqslant \bar N\right\}.$$

On the other hand, we notice that if $\nu_i^{(\infty)}>0$, then $\left(\nu_i^{(k)}\right)^{\theta_k} \to 1$ as $k \to \infty$; if $\nu_i^{(\infty)}>0$, then $\liminf_{k \to \infty}\left(\nu_i^{(k)}\right)^{\theta_k}\geqslant 0$. Hence, putting these facts together we obtain
\begin{align*}
\underline{\Theta}=&\lim_{k \to \infty}\Theta(m, \theta_k, n-1)\\
=&\lim_{k \to \infty} \sum_{i=1}^{\bar N}\left(\nu_i^{(k)}\right)^{\theta_k}\\
\geqslant& \lim_{k \to \infty} \sum_{\substack{1\leqslant i \leqslant \bar N\\\nu_i^{(\infty)}>0}}\left(\nu_i^{(k)}\right)^{\theta_k}=\sharp\left\{\nu_i^{(\infty)}>0,~~1 \leqslant i \leqslant \bar N\right\}\geqslant N_m
\end{align*}
as desired.
\end{proof}
\vskip 8pt

We now transfer to the construction of precise test functions in order to show that the constant $c(n, \frac{n}{n-2})/\Theta(m, \frac{n-2}{n-1}, n-1)$ in the inequality \eqref{ineq:2nd order} is \emph{almost optimal}.

\begin{proof}[Proof of Proposition \protect\ref{example:n>2}]
For each $m$, there exist some natural number $N\geqslant N_m(\Sp^{n-1})$ and   $\nu=\sum_{i=1}^N \nu_i \delta_{x_i} \in  \mathcal{M}_{m}^{c}$ for $1 \leqslant i \leqslant N$, such that
$$\Theta(m, \frac{n-2}{n-1}, n-1)=\sum_{i=1}^N \nu_i^{\frac{n-2}{n-1}}.$$	
 
	 We denote by $\stackrel{\frown}{xy}$  the geodesic distance between $x$ and $y$ in $\mathbb{S}^{n-1}$. Fix $\delta>0$ small enough such that $\overline{\mathcal{A}_{2 \delta}\left(x_{i}\right)} \cap \overline{\mathcal{A}_{2 \delta}\left(x_{j}\right)}=\emptyset$ for $1 \leqslant i<j \leqslant N$, where
	 \begin{equation}\label{def:Conic_Annulus}
	\mathcal{A}_\delta (x_i):=\left\{x=r\xi \in \mathbb{B}^n; \xi \in \Sp^{n-1}, 1-r< \delta, ~~ \stackrel{\frown}{\xi x_i} < \delta \right\}
	\end{equation}
	for each $x_i \in \mathbb{S}^{n-1}$. In other words, $\mathcal{A}_\delta (x_i)$ is a conic annulus, where the cone has  the origin as its vertex and a geodesic ball $B_{\delta}(x_i) \subset \Sp^{n-1}$ as its base. 
	
	 For $0 < \ve <\delta$ and each $1 \leqslant i \leqslant N$, under the above coordinates we define
	\begin{align*}
		\phi_{\ve, i}(x)= \chi_i(x)\left((\ve+1-r)^2+\stackrel{\frown}{\xi x_i}^2\right)^{\frac{2-n}{2}},
		\end{align*}
		where $\chi_i (x)$ is a smooth cut-off function, $\chi_i(x)=1$ in $\mathcal{A}_{\delta}(x_i)$ and $\chi_i(x)=0$ outside $\mathcal{A}_{2\delta}(x_i)$.
		
Define 
$$ v(x)=\sum_{i=1}^{N}\nu_i^{\frac{n-2}{2(n-1)}}\phi_{\ve, i} (x),$$
then a direct computation yields
\begin{align}\label{bubble esti1_v}
	&\int_{\mathbb{S}^{n-1}}v^{\frac{2(n-1)}{n-2}}\ud \mu_{\Sp^{n-1}}\no\\
	=&\sum_{i=1}^{N}\nu_i\int_{B_{2\delta}(x_i) }\phi_{\ve, i}^{\frac{2(n-1)}{n-2}}\ud \mu_{\Sp^{n-1}}\no\\
	=&\left(\sum_{i=1}^{N}\nu_i \right)|\Sp^{n-2}|\int_0^{\delta}(\ve^2+\rho^2)^{1-n}\sin^{n-2}\rho \ud \rho+O(\ve^{1-n})\int_{\frac{\delta}{\ve}}^{\frac{2\delta}{\ve}}(1+\rho^2)^{1-n}\rho^{n-2}\ud \rho\no\\
	=&(2\ve)^{1-n}\left|\mathbb{S}^{n-1}\right|+ \begin{cases} 
	O(\log \ve^{-1}) \quad \mathrm{if~~} n=3;\\
	O(\ve^{3-n}) \qquad \mathrm{if~~} n\geqslant 4;
	\end{cases}\no\\
	=&(2\ve)^{1-n}\left|\mathbb{S}^{n-1}\right|+O(\ve^{3-n}\log \ve^{-1}) \quad \mathrm{~~as~~} \varepsilon \rightarrow 0.
\end{align}

 For any $p \in \mathring{\mathcal{P}}_{m}$, we have 
\begin{align}\label{correction esti:}
& \int_{\mathbb{S}^{n-1}} v^{\frac{2(n-1)}{n-2}} p \ud \mu_{\Sp^{n-1}} \no\\
=& \sum_{i=1}^{N}\nu_i \int_{B_{2 \delta}\left(x_{i}\right)} \phi_{\ve, i}^{\frac{2(n-1)}{n-2}}(x) p(x) \ud \mu_{\Sp^{n-1}} \no\\
=& \sum_{i=1}^{N}\int_{B_{2 \delta}\left(x_{i}\right)} \left[\phi_{\ve, i}^{\frac{2(n-1)}{n-2}}(\stackrel{\frown}{x x_{i}}) \nu_i p\left(x_{i}\right) +\phi_{\ve, i}^{\frac{2(n-1)}{n-2}}(\stackrel{\frown}{x x_{i}}) O\left(\stackrel{\frown}{x x_{i}}^{2}\right)\right] \ud \mu_{\Sp^{n-1}} \no \\
=& O\left(\ve^{3-n}\log \ve^{-1}\right),
\end{align}
where the last equality follows from
\begin{align*}
	\nu \in \mathcal{M}_{m}^{c}\left(\mathbb{S}^{n-1}\right) \quad \Longrightarrow \quad \sum_{i=1}^{N} \nu_i p\left(x_{i}\right)=0.
\end{align*}

Obviously, we shall make a further correction for the above $v$ to fulfill the condition \eqref{moment vanishing}. It is shown in \cite[Theorem 2.1 in Chapter IV]{Stein-Weiss} that there exists a basis 
$\{P_{1}, \cdots, P_{L}\}$ of $\mathring{\mathcal{P}}_{m},$
 such that 
 $$\left.p_{1}=P_1\right|_{\mathbb{S}^{n-1}}, \cdots,\left.p_{L}=P_L\right|_{\mathbb{S}^{n-1}}$$
  are spherical harmonics, where $L=n+\sum_{i=2}^{m}(C_{n+i-1}^{n-1}-C_{n+i-3}^{n-1})$. Then  for each $1\leqslant i \leqslant N$, we claim that there exist $\psi_{1}, \cdots, \psi_{L} \in C_{c}^{\infty}\left(\mathbb{B}^{n} \backslash \bigcup_{i=1}^{N} \overline{\mathcal{A}_{2 \delta}\left(x_{i}\right)}\right)$ such that the determinant 
\begin{align}\label{nonsingular}
	\det\left[\int_{\mathbb{S}^{n-1}} \psi_{j} p_{k} \ud \mu_{\Sp^{n-1}}\right]_{1 \leqslant j, k \leqslant L} \neq 0.
\end{align}
To this end, we can choose a nonzero smooth function $\eta \in C_{c}^{\infty}\left(\mathbb{B}^{n} \backslash \bigcup_{i=1}^{N} \overline{\mathcal{A}_{2 \delta}\left(x_{i}\right)}\right)$ such that $\eta P_{1}, \cdots, \eta P_{L}$ are linearly independent. It follows that the Gram matrix
\begin{align*}
	\left[\int_{\mathbb{S}^{n-1}} \eta^{2} p_{j} p_{k} \ud \mu_{\Sp^{n-1}}\right]_{1 \leqslant j, k \leqslant L}
\end{align*}
is positive definite, then $\psi_{j}=\eta^{2} P_{j}$ satisfies \eqref{nonsingular}.

The fact \eqref{nonsingular} enables us to find $\beta_{1}, \cdots, \beta_{L} \in \mathbb{R}$ such that 
\begin{align}\label{moment vanishing construction}
	\int_{\mathbb{S}^{n-1}}\left(v^{\frac{2(n-1)}{n-2}}+\sum_{j=1}^{L} \beta_{j} \psi_{j}\right) p_{k} \ud \mu_{\Sp^{n-1}}=0 \qquad \forall~ 1 \leqslant k \leqslant L.
\end{align}
Moreover, it follows from \eqref{correction esti:} that for all $1 \leqslant j \leqslant L$, $\beta_{j}=O\left(\ve^{3-n}\log \ve^{-1}\right)$ as $\varepsilon \rightarrow 0$. As a consequence we can find a constant $c_1 >0$ such that 
\begin{align*}
	\sum_{j=1}^{L} \beta_{j} \psi_{j}+c_{1} \ve^{3-n}\log \ve^{-1} \geqslant \ve^{3-n}\log \ve^{-1}.
\end{align*} 

We define the test function by
\begin{align}\label{def:test_fcn}
	u^{\frac{2(n-1)}{n-2}}=v^{\frac{2(n-1)}{n-2}}+\sum_{j=1}^{L} \beta_{j} \psi_{j}+c_{1} \ve^{3-n}\log \ve^{-1}.
\end{align}
Clearly, it follows from \eqref{moment vanishing construction} and \eqref{def:test_fcn} that 
\begin{align*}
	\int_{\mathbb{S}^{n-1}}p u^{\frac{2(n-1)}{n-2}}\ud \mu_{\Sp^{n-1}}=0, \qquad \forall~ p  \in \mathring{\mathcal{P}}_{m}.
\end{align*}

Hence, as $\varepsilon \rightarrow 0$, by \eqref{bubble esti1_v} and \eqref{def:test_fcn} we have
\begin{align}\label{bubble esti1}
	&\|u\|_{L^{\frac{2(n-1)}{n-2}}(\mathbb{S}^{n-1})}^2 \no\\
	=&\left[\int_{\Sp^{n-1}}\left(v^{\frac{2(n-1)}{n-2}}+\sum_{j=1}^{L} \beta_{j} \psi_{j}+c_{1} \ve^{3-n}\log \ve^{-1}\right) \ud \mu_{\Sp^{n-1}}\right]^{\frac{n-2}{n-1}}\no\\
	=&(2\ve)^{2-n}\left|\mathbb{S}^{n-1}\right|^{\frac{n-2}{n-1}} \left(1 + O(\ve^2\log \ve^{-1})\right) .
\end{align}

Next, we estimate the term $\|u\|_{L^2(\Sp^{n-1})}$. To this end,  we can apply \eqref{def:test_fcn} to show
\begin{align*}
		u^{\frac{2(n-1)}{n-2}} \leqslant v^{\frac{2(n-1)}{n-2}} +  C\ve^{3-n}\log \ve^{-1},
\end{align*}
which directly yields
\begin{align*}
	u^2 \leqslant \left(v^{\frac{2(n-1)}{n-2}} + C \ve^{3-n}\log \ve^{-1} \right)^{\frac{n-2}{n-1}} \leqslant v^2 + C  (\ve^{3-n}\log \ve^{-1})^{\frac{n-2}{n-1}}. 
\end{align*}
Hence, as $\varepsilon \rightarrow 0$ we obtain
\begin{align}\label{bubble esti2}
	&\int_{\mathbb{S}^{n-1}}u^2 \ud \mu_{\Sp^{n-1}}\no\\
	\leqslant &\int_{\mathbb{S}^{n-1}}\left(v^2 + C(\ve^{3-n}\log \ve^{-1})^{\frac{n-2}{n-1}}\right) \ud \mu_{\Sp^{n-1}}\no\\
	= & \sum_{i=1}^{N}\nu_i^{\frac{n-2}{n-1}}\int_{B_{2\delta}(x_i)} \phi_{\ve, i}^2 \ud \mu_{\Sp^{n-1}} + C (\ve^{3-n}\log \ve^{-1})^{\frac{n-2}{n-1}}\no\\
	= & O (\ve^{3-n}\log \ve^{-1})^{\frac{n-2}{n-1}}=O(\ve^{2-n})(\ve^2 \log \ve^{-1})^{\frac{n-2}{n-1}}.
\end{align}

To estimate $\|\nabla u\|_{L^{2}(\mathbb{B}^n)}^{2}$. Near each $x_i, 1 \leqslant i \leqslant N$ the flat metric $|\ud x|^2$ in $\overline{\Bn}$ is expressed as
$$|\ud x|^2=\ud r^2 +r^2 \left(\ud\rho^2 + \sin^2\rho g_{\mathbb{S}^{n-2}}\right)$$
for $x=r\xi \in \Bn$ and $\rho=\ud_{\mathbb{S}^{n-1}}(\xi, x_i)=\stackrel{\frown}{\xi x_i}$. Under this coordinate system, we rewrite
$$\phi_{\ve, i}(x)=\chi_i(x)\left((\ve+1-r)^2 +\rho^2\right)^{\frac{2-n}{2}}.$$

Recall that the Beta function is defined by
$$\int_0^\infty\frac{x^{\alpha-1}}{(1+x)^{\alpha+\beta}}\ud x=B(\alpha,\beta)=\frac{\Gamma(\alpha)\Gamma(\beta)}{\Gamma(\alpha+\beta)}$$
for $\alpha,\beta \in \mathbb{C}$ with ${\rm Re}(\alpha), {\rm Re}(\beta)>0$. 

We are ready to calculate
\begin{align*}
&\|\nabla u\|_{L^{2}(\mathbb{B}^n)}^{2} \no\\
=& \sum_{i=1}^{N} \int_{\mathcal{A}_{2\delta}\left(x_{i}\right)}|\nabla u|^{2}\ud x+O\left( \ve^{2(3-n)}\log \ve^{-2}\right)\no\\
=& \sum_{i=1}^{N} \nu_i^{\frac{n-2}{n-1}}\int_{\mathcal{A}_{\delta}\left(x_{i}\right)}\left(1+c_1 \nu_i^{-1} \ve^{3-n}\log \ve^{-1}\phi_{\ve, i}^{\frac{2(1-n)}{n-2}}\right)^{-\frac{n}{n-1}}|\nabla \phi_{\ve, i}|^{2} \ud x+O\left(\ve^{2(3-n)}\log \ve^{-2}\right),
\end{align*}
where the last identity follows from the estimate that in $\mathcal{A}_{2\delta}(x_i)$,
\begin{align*}
u^{\frac{2(n-1)}{n-2}}=\nu_i \phi_{\ve, i}^{\frac{2(n-1)}{n-2}}+c_{1} \ve^{3-n}\log \ve^{-1} \quad\Longrightarrow\quad |\nabla u|^2=\nu_i^2 \left(\frac{\phi_{\ve, i}}{u}\right)^{\frac{2n}{n-2}}|\nabla \phi_{\ve, i}|^2.
\end{align*}
Through a direct computation showing that in $\mathcal{A}_\delta(x_i)$, 
\begin{align*}
	|\nabla \phi_{\ve, i}|^2 =(n-2)^2 \frac{(\ve+1-r)^2+r^{-2}\rho^2}{\left((\ve+1-r)^2 +\rho^2\right)^{n}}.
\end{align*}
Hence, the above integral on the right hand side can be estimated by
\begin{align*}
	&\int_{\mathcal{A}_{\delta}\left(x_{i}\right)}\left(1+c_1  \nu_i^{-1} \ve^{3-n}\log \ve^{-1}\phi_{\ve, i}^{\frac{2(1-n)}{n-2}}\right)^{-\frac{n}{n-1}}|\nabla \phi_{\ve, i}|^{2} \ud x\\
	\overset{t=1-r}{=}&(n-2)^2 |\mathbb{S}^{n-2}|\int_{0}^{\delta}\int_{0}^{\delta}\frac{\left((\ve+t)^2+\rho^2\right)^{1-n}(1-t)^{n-1}\sin^{n-2}\rho}{\left(1+c_1  \nu_i^{-1} \ve^{3-n}\log \ve^{-1}\left((\ve+t)^2+\rho^2\right)^{n-1}\right)^{\frac{n}{n-1}}}  \ud\rho \ud t+I_1\\
	\leqslant&(n-2)^2 |\mathbb{S}^{n-2}|\int_{0}^{\delta}\int_{0}^{\delta}\frac{\left((\ve+t)^2+\rho^2\right)^{1-n}\rho^{n-2}}{\left(1+c_1  \nu_i^{-1} \ve^{3-n}\log \ve^{-1}\left((\ve+t)^2+\rho^2\right)^{n-1}\right)^{\frac{n}{n-1}}} \ud\rho \ud t+O\left(\ve^{3-n}\right)\\
	\leqslant &(n-2)^2|\mathbb{S}^{n-2}|\int_{0}^{\delta}\int_{0}^{\delta}{\left((\ve+t)^2+\rho^2\right)^{1-n}}\rho^{n-2} \ud\rho \ud t +O\left(\ve^{3-n}\right)\\
	=&(n-2)^2|\mathbb{S}^{n-2}| \ve^{2-n}\int_{0}^{\frac{\delta}{\ve}}\int_{0}^{\frac{\delta}{\ve}}{\left((1+t)^2+\rho^2\right)^{1-n}}\rho^{n-2} \ud\rho \ud t +O\left(\ve^{3-n}\right)\\
	=&(n-2)^2|\mathbb{S}^{n-2}| \ve^{2-n}\int_{0}^{\infty}\int_{0}^{\infty}{\left((1+t)^2+\rho^2\right)^{1-n}}\rho^{n-2} \ud\rho \ud t +O\left(\ve^{3-n}\right)\\
	=&\frac{n-2}{2}|\mathbb{S}^{n-2}|B(\frac{n-1}{2}, \frac{n-1}{2})\ve^{2-n}+O\left(\ve^{3-n}\right),
\end{align*}
where the first inequality follows from 
\begin{align*}
I_1=&(n-2)^2 |\mathbb{S}^{n-2}|\int_{0}^{\delta}\int_{0}^{\delta}\frac{\left((\ve+t)^2+\rho^2\right)^{1-n}t(2-t)(1-t)^{n-3}}{\left(1+c_1  \nu_i^{-1} \ve^{3-n}\log \ve^{-1}\left((\ve+t)^2+\rho^2\right)^{n-1}\right)^{\frac{n}{n-1}}} \sin^{n-2}\rho \ud\rho \ud t\\
\leqslant& 2(n-2)^2 |\mathbb{S}^{n-2}|\int_{0}^{\delta}\int_{0}^{\delta}\left((\ve+t)^2+\rho^2\right)^{1-n}t\rho^{n-2} \ud\rho \ud t=O(\ve^{3-n}).
\end{align*}

On the other hand, we can further require that 
$$(1-r)^2\leqslant \ve^{1-\frac{2(1-\e_0)}{n-1}}<\delta^2 \quad \mathrm{~~and~~} \quad \rho^2\leqslant \ve^{1-\frac{2(1-\e_0)}{n-1}}<\delta^2$$
 for some $0<\e_0<1/2$,  then
\begin{align*}
	1+c_1  \nu_i^{-1} \ve^{3-n}\log \ve^{-1}\left((\ve+1-r)^2+\rho^2\right)^{n-1}=1+O(\ve^{2\e_0}\log \ve^{-1}).
\end{align*}
With this estimate at hand, it is not hard to show that $I_1=O(\ve^{3-n})$.

Hence, we obtain
\begin{align*}
	&\int_{\mathcal{A}_{\delta}\left(x_{i}\right)}\left(1+c_1  \nu_i^{-1} \ve^{3-n}\log \ve^{-1} \phi_{\ve, i}^{\frac{2(1-n)}{n-2}}\right)^{-\frac{n}{n-1}}|\nabla \phi_{\ve, i}|^{2}  \ud x\\
	=&(n-2)^2 |\mathbb{S}^{n-2}|\int_{0}^{\delta}\int_{0}^{\delta}\frac{\left((\ve+t)^2+\rho^2\right)^{1-n}\rho^{n-2}}{\left(1+c_1 \nu_i^{-1} \ve^{3-n}\log \ve^{-1}\left((\ve+t)^2+\rho^2\right)^{n-1}\right)^{\frac{n}{n-1}}} \ud\rho \ud t\\
	&+I_2 +O\left(\ve^{3-n}\right)\\
	\geqslant &(n-2)^2|\mathbb{S}^{n-2}|(1+O(\ve^{2\e_0}\log \ve^{-1}))\int_{0}^{\ve^{\frac{1}{2}-\frac{1-\e_0}{n-1}}}\int_{0}^{\ve^{\frac{1}{2}-\frac{1-\e_0}{n-1}}}{\left((\ve+t)^2+\rho^2\right)^{1-n}}\rho^{n-2} \ud\rho \ud t \\
	&+O\left(\ve^{3-n}\right)\\
	=&(n-2)^2|\mathbb{S}^{n-2}| \ve^{2-n}\int_{0}^{\ve^{-\frac{1}{2}-\frac{1-\e_0}{n-1}}}\int_{0}^{\ve^{-\frac{1}{2}-\frac{1-\e_0}{n-1}}}{\left((1+t)^2+\rho^2\right)^{1-n}}\rho^{n-2} \ud\rho \ud t\\
	& +O\left(\ve^{2-n+2\e_0}\log \ve^{-1}\right)\\
	=&(n-2)^2|\mathbb{S}^{n-2}| \ve^{2-n}\int_{0}^{\infty}\int_{0}^{\infty}{\left((1+t)^2+\rho^2\right)^{1-n}}\rho^{n-2} \ud\rho \ud t +O\left(\ve^{2-n+2\e_0}\log \ve^{-1}\right)\\
	=&\frac{n-2}{2}|\mathbb{S}^{n-2}|B(\frac{n-1}{2}, \frac{n-1}{2})\ve^{2-n}+O\left(\ve^{2-n+2\e_0}\log \ve^{-1}\right),
\end{align*}
where the above inequality follows from
\begin{align*}
-I_2=&(n-2)^2 |\mathbb{S}^{n-2}|\int_{0}^{\delta}\int_{0}^{\delta}\frac{\left((\ve+t)^2+\rho^2\right)^{1-n}\left(\rho^{n-2} -\sin^{n-2}\rho\right)}{\left(1+c_1  \nu_i^{-1} \ve^{3-n}\log \ve^{-1}\left((\ve+t)^2+\rho^2\right)^{n-1}\right)^{\frac{n}{n-1}}}\ud\rho \ud t\\
&+(n-2)^2 |\mathbb{S}^{n-2}|\int_{0}^{\delta}\int_{0}^{\delta}\frac{\left((\ve+t)^2+\rho^2\right)^{1-n}\left(1-(1-t)^{n-1}\right)\sin^{n-2}\rho}{\left(1+c_1 \nu_i^{-1} \ve^{3-n}\log \ve^{-1}\left((\ve+t)^2+\rho^2\right)^{n-1}\right)^{\frac{n}{n-1}}}  \ud\rho \ud t\\
\leqslant&C  \int_{0}^{\delta}\int_{0}^{\delta}\left((\ve+t)^2+\rho^2\right)^{1-n} \rho^{n-2}(\rho^2+t)\ud\rho \ud t=O(\ve^{3-n}).
\end{align*}

In summary,  for each $1 \leqslant i \leqslant N$ we have
\begin{align*}
	&\int_{\mathcal{A}_{\delta}\left(x_{i}\right)}\left(1+c_1  \nu_i^{-1} \ve^{3-n}\log \ve^{-1}\phi_{\ve, i}^{\frac{2(1-n)}{n-2}}\right)^{-\frac{n}{n-1}}|\nabla \phi_{\ve, i}|^{2}  \ud x\\
=&\ve^{2-n}\frac{n-2}{2}|\mathbb{S}^{n-2}|B(\frac{n-1}{2}, \frac{n-1}{2})+O\left(\ve^{2-n+2\e_0}\log \ve^{-1}\right).
\end{align*}

Consequently, we conclude that
\begin{align}\label{bubble esti3}
	&\|\nabla u\|_{L^{2}(\mathbb{B}^n)}^{2}\no \\
	=& \sum_{i=1}^{N} \nu_i^{\frac{n-2}{n-1}}\int_{\mathcal{A}_{\delta}\left(x_{i}\right)}\left(1+c_1  \nu_i^{-1} \ve^{3-n}\log \ve^{-1}\phi_{\ve, i}^{\frac{2(1-n)}{n-2}}\right)^{-\frac{n}{n-1}}|\nabla \phi_{\ve, i}|^{2} \ud x+O\left(\ve^{2-n+2\e_0}\log \ve^{-1}\right)\no\\
	=& \ve^{2-n}\frac{(n-2)}{2}\big(\sum_{i=1}^{N} \nu_i^{\frac{n-2}{n-1}}\big)|\mathbb{S}^{n-2}|B(\frac{n-1}{2}, \frac{n-1}{2})+O\left(\ve^{2-n+2\e_0}\log \ve^{-1}\right)\no\\
	=&\ve^{2-n} 2^{1-n}(n-2)|\Sp^{n-1}|\Theta(m, \frac{n-2}{n-1}, n-1)(n-2)+O\left(\ve^{2-n+2\e_0}\log \ve^{-1}\right).
\end{align}

Therefore, inserting \eqref{bubble esti1}, \eqref{bubble esti2}, \eqref{bubble esti3} into \eqref{ineq:almost_optimal} and dividing both sides by $\ve^{2-n}$, next letting $\ve\to 0$, we obtain 
\begin{align*}
	2^{2-n}\left|\mathbb{S}^{n-1}\right|^{\frac{n-2}{n-1}} \leqslant a 2^{1-n}\Theta(m, \frac{n-2}{n-1}, n-1)(n-2)|\Sp^{n-1}| ,
\end{align*}
i.e.,
\begin{align*}
	a \geqslant \frac{c(n, \frac{n}{n-2})}{\Theta(m, \frac{n-2}{n-1}, n-1)}.
\end{align*}
This completess our construction.
\end{proof}

\section{Almost sharp  Sobolev trace inequality of order four under constraints in dimension five and higher}\label{Sect4}

We utilize the conic proof to derive Sobolev trace inequality of fourth order for $n\geqslant 5$. With some deep insight on ``\emph{local bubble}", we are able to complete the construction of \emph{almost sharp} example.

\subsection{Fourth order Sobolev trace inequality}\label{subsect:4th order ineq}

Inspired by the proof of the second order Sobolev trace inequality, we know that Ache-Chang's sharp Sobolev trace inequalities of order four in Theorem \ref{Thm:Ache-Chang} is a prerequisite for the following compactness  and concentration lemma. For the \emph{almost sharp} example, we do care about the equality in \eqref{ineq:Ache-Chang_n>4} that if $f=1$, then the extremal metric $u^{4/(n-4)} |\ud x|^2$ is the \emph{adapted metric} in the Poincar\'e model  $(\Bn,\Sp^{n-1},g_{\mathbb{H}})$ of hyperbolic space for $n\geqslant 5$, first introduced by Case-Chang \cite{Case-Chang}; see also \cite[Proposition 2.2]{Ache-Chang}.  

\vskip 8pt

Based on Theorem \ref{Thm:Ache-Chang}, as before we begin with the refinement of the concentration and compactness principle, whose proof presented below is similar in spirit to the one of Lemma \ref{Concentration compactness bdry}.
\begin{lemma}\label{Concentration compactness bdry_4th}
	For $(\mathbb{B}^n, |\ud x|^2)$ and $n\geqslant 5$, let $u_k \in H^2 (\Bn)$ be a sequence of $H^2(\Bn)$ extensions of $f_k \in C^\infty(\Sp^{n-1})$ satisfying the Neumann boundary condition
	$$\frac{\pa u_k}{\pa r} =-\frac{n-4}{2}f_k \qquad \mathrm{~~on~~} \Sp^{n-1}.$$
	Assume that as $k \to \infty$, $u_k \rightharpoonup u$ weakly in $H^2(\mathbb{B}^n)$ and
	\begin{align*}
		(\Delta u_{k})^{2} \ud x \rightharpoonup&~~ (\Delta u)^{2} \ud x+\lambda \qquad~~~\qquad\mathrm{as ~~measure~~in~~} \mathbb{B}^n,\\
		\left|f_{k}\right|^{\frac{2(n-1)}{n-4}} \ud \mu_{\Sp^{n-1}}  \rightharpoonup &~~ |f|^{\frac{2(n-1)}{n-4}} \ud \mu_{\Sp^{n-1}}+\nu \qquad\mathrm{as~~measure~~in~~} \mathbb{S}^{n-1},
			\end{align*}
	where $\lambda$ and $\nu$ are bounded nonnegative measures on $\Bn$ and $\Sp^{n-1}$, respectively.
	Then there exist countably many points $x_i \in \mathbb{S}^{n-1}$ with $\nu_i=\nu(\{x_i\})>0$ such that 
	\begin{align*}
		&\nu=\sum_{i} \nu_{i} \delta_{x_{i}} \qquad \mathrm{~~and~~} \qquad \hat{\lambda} \geqslant \frac{1}{\alpha(n)}\sum_{i}\nu_i^{\frac{n-4}{n-1}} \delta_{x_i}.
		\end{align*}
\end{lemma}
\begin{proof}
	Let $v_k= u_k -u$, then it follows from the assumptions that up to a subsequence, as $k \to \infty$, $v_k \rightharpoonup 0$ weakly in $H^{2}(\mathbb{B}^n)$ and $v_k \to 0$ in $H^1 (\mathbb{B}^n)$ and $f_k \to 0$ in $H^1 (\mathbb{S}^{n-1})$; also  $\nu_{k}:=|u_k|^{\frac{2(n-1)}{n-4}} \ud\mu_{\Sp^{n-1}}-|u|^{\frac{2(n-1)}{n-4}} \ud\mu_{\Sp^{n-1}}\rightharpoonup \nu$ and $\lambda_{k}:=\left(\Delta v_{k}\right)^{2} \ud x \rightharpoonup \lambda$.
	
	 	For any $\varphi \in C^2(\overline{\mathbb{B}^n})$ with $\pa \varphi /\pa r=0$ on $\Sp^{n-1}$, by the Sobolev trace inequality of order four in Theorem \ref{Thm:Ache-Chang} we have 	
		\begin{align*}
\left(\int_{\mathbb{S}^{n-1}}|\varphi|^{\frac{2(n-1)}{n-4}} \ud \nu\right)^{\frac{n-4}{n-1}}=&\left(\lim _{k \rightarrow \infty} \int_{\mathbb{S}^{n-1}}|\varphi|^{\frac{2(n-1)}{n-4}} \ud \nu_{k}\right)^{\frac{n-4}{n-1}}\\
=&\left(\lim _{k \rightarrow \infty} \int_{\mathbb{S}^{n-1}}\left|v_{k} \varphi\right|^{\frac{2(n-1)}{n-4}} \ud \mu_{\Sp^{n-1}} \right)^{\frac{n-4}{n-1}}\\
\leqslant& \alpha(n) \liminf _{k \rightarrow \infty}\left[\int_{\mathbb{B}^{n}}\left|\Delta \left(v_{k} \varphi\right)\right|^{2} \ud x\right]\\
=& \alpha(n) \liminf _{k \rightarrow \infty}\left[\int_{\mathbb{B}^{n}}\varphi^{2}\left(\Delta v_{k}\right)^{2} \ud x\right] \\
=&\alpha(n) \int_{\mathbb{B}^{n}}\varphi^{2} \ud \lambda,
\end{align*}
that is,
	\begin{align}\label{trace esti}
		\left(\int_{\mathbb{S}^{n-1}}|\varphi|^{\frac{2(n-1)}{n-4}} \ud \nu \right)^{\frac{n-4}{n-1}}\leqslant \alpha(n)\int_{\mathbb{B}^{n}}\varphi^{2} \ud \lambda
	\end{align}
	for all $\varphi \in C^2(\overline{\mathbb{B}^n})$ with $\pa \varphi /\pa r=0$ on $\Sp^{n-1}$. 
	
	Now let $J \subset \Sp^{n-1}$ be the set of the atoms of the measure $\nu$. From the assumption and the Sobolev trace inequality that
	 $$\int_{\mathbb{S}^{n-1}}\ud\nu\leqslant \lim_{k\to \infty}\int_{\Sp^{n-1}}|u_k|^{\frac{2(n-1)}{n-4}} \ud \mu_{\Sp^{n-1}} \leqslant C\liminf_{k\to \infty} \|u_k\|_{H^2(\Bn)}^{\frac{2(n-1)}{n-4}}< \infty,$$
	 we know that $J$ is an at most countable set denoted by $J=\{x_i; i \in \N\}$.  We decompose 
	 $$\nu = \nu_0 + \sum_{i}\nu_i \delta_{x_i}$$ 
	where $\nu_0 $ is the singular continuous part of the measure $\nu$ and also a nonnegative measure. For any Borel set $\Omega \subset \Sp^{n-1}$, same as in Lemma \ref{Concentration compactness bdry}, we define 
		$$\hat{\lambda}(\Omega):= \lambda(\mathcal{C}(\Omega)),$$
where $\mathcal{C}(\Omega)$ is the cone defined in \eqref{def:conic_Borel_set}. Since  both $\nu$ and $\lambda$ are bounded Borel measures, so does  $\hat{\lambda}$ by definition. Therefore, given any Borel set $\Omega$ in $\mathbb{S}^{n-1}$ and any $\ve>0$, there exist an open set $\mathcal{O}$ and a closed set $K$ such that $K\subset \Omega \subset \mathcal{O}$ with 
	\begin{align*}
		\nu(\mathcal{O}\setminus \Omega)< \ve, \quad \hat{\lambda}(\mathcal{O}\setminus\Omega)< \ve, 
	\end{align*}
	and \begin{align*}
		\nu(\Omega\setminus K)< \ve, \quad \hat{\lambda}(\Omega\setminus K)< \ve.
	\end{align*}
For all $\varphi \in C^\infty (\overline{\mathbb{B}^n})$ with $\mathrm{supp}~\varphi \subset \mathcal{C}(\mathcal{O})$  such that $\varphi|_{\mathcal{C}(K)}=1$, $0 \leqslant \varphi \leqslant 1$, $\pa \varphi /\pa r=0$ on $\Sp^{n-1}$, we conclude from \eqref{trace esti} that  
 \begin{align*}
	 	\left(\int_{K} \ud\nu\right)^{\frac{n-4}{n-1}} \leqslant \alpha(n)\left(\int_{\mathcal{C}({\mathcal{O}})} \ud \lambda\right).
	 \end{align*}
	 Hence, $$\left(\nu(\Omega)-\ve\right)^{\frac{n-4}{n-1}} \leqslant \alpha(n) \left(\hat{\lambda}(\Omega) +\ve \right).$$
	 Letting $\ve \to 0$, we have 
	 \begin{align*}
	 	\nu(\Omega)^{\frac{n-4}{n-1}} \leqslant \alpha(n) \hat{\lambda}(\Omega)<\infty.
	 \end{align*}
	 This indicates that $\nu$ is  absolutely continuous with respect to $\hat{\lambda}$ and thus from the Radon-Nikodym theorem that there exists $f \in L^{1}\left(\mathbb{B}^{n} ; \hat{\lambda}\right)$ such that $\ud \nu=f \ud \hat{\lambda}$.  Moreover, for $\hat{\lambda}$-almost every $x\in \mathbb{S}^{n-1}$ and  $B_{\rho}(x) \subset \mathbb{S}^{n-1}$, we have  
	 \begin{align*}
	 	f(x)=\lim _{\rho \searrow 0}\frac{\int_{B_{\rho}(x)} \ud \nu}{\int_{B_{\rho}(x)} \ud \hat{\lambda}}.
	 \end{align*}
	 In particular, if  we choose $x \in \Sp^{n-1}$ and the segment $\overline{ox}$ does not carry any atom of $\lambda$, then
	 \begin{align*}
	 	f(x)^{\frac{n-4}{n-1}}=\lim _{\rho \searrow 0}\left(\frac{\int_{B_{\rho}(x)} \ud \nu}{\int_{B_{\rho}(x)} \ud \hat{\lambda}}\right)^{\frac{n-4}{n-1}} \leqslant \alpha(n)\lim _{\rho \searrow 0}\left(\int_{\mathcal{C}(B_{\rho}(x))} \ud \lambda\right)^{\frac{3}{n-1}}=0.
	 \end{align*}
	From this and \eqref{trace esti},  a similar argument in Lemma \ref{Concentration compactness bdry} yields that $\nu_0=0$ and
	$$\hat{\lambda} \geqslant \frac{1}{\alpha(n)}\sum_{i}\nu_i^{\frac{n-4}{n-1}} \delta_{x_i}.$$
	This completes the proof.
	\end{proof}

We are now ready to prove Theorem \ref{Thm1:Ache-Chang-type}.
\begin{proof}[Proof of Theorem \protect\ref{Thm1:Ache-Chang-type}]
Following the same lines in the proof of  \cite[Theorem A]{Ache-Chang}, we only need to consider the case when $u$ is a biharmonic extension of $f$ to $\Bn$. For brevity, we define
$$\beta_\ve =  \frac{\alpha (n)}{\Theta(m, \frac{n-4}{n-1}, n-1)} +\ve.$$

	By contradiction, if \eqref{ineq:n>5} is not true, then we can find some $\beta_\ve$ and a sequence of functions $\{u_k\} \subset H^2(\mathbb{B}^n)$, which are the biharmonic extensions of $f_k$ to $\Bn$ satisfying the Neumann boundary condition
\begin{align*}
	\frac{\pa u_k}{\pa r} =-\frac{n-4}{2}f_k \qquad \mathrm{~~on~~} \Sp^{n-1} ,
\end{align*}
such that
	\begin{align*}
&\left(\int_{\mathbb{S}^{n-1}}|f_k|^{\frac{2(n-1)}{n-4}} \ud \mu_{\Sp^{n-1}}\right)^{\frac{n-4}{n-1}} \\
>&\beta_\ve\int_{\mathbb{B}^{n}}\left(\Delta u_k\right)^{2} \ud x+k  \int_{\mathbb{S}^{n-1}}\left(|\nabla f_k|_{\Sp^{n-1}}^2  +b_n  f_{k}^2\right) \ud\mu_{\Sp^{n-1}}	
\end{align*} 
and
\begin{align*}
		\int_{\mathbb{S}^{n-1}} p |f_k|^{\frac{2(n-1)}{n-4}} \ud \mu_{\Sp^{n-1}}=0, \qquad \forall~ p \in \mathring{\mathcal{P}}_{m}.
	\end{align*}
	
We may normalize $\int_{\mathbb{S}^{n-1}}|f_k|^{2(n-1)/(n-4)} \ud \mu_{\Sp^{n-1}}=1$ such that
	\begin{align*}
		\int_{\mathbb{B}^{n}}|\Delta u_k|^{2} \ud x \leqslant \frac{1}{\beta_\ve}, \qquad\int_{\mathbb{S}^{n-1}}f_k^{2} \ud \mu_{\Sp^{n-1}}\leqslant \frac{1}{b_n k}
	\end{align*}
	and 
	\begin{align*}
		\int_{\mathbb{S}^{n-1}}|\nabla f_k|_{\Sp^{n-1}}^{2} \ud \mu_{\Sp^{n-1}}\leqslant \frac{1}{k}.
	\end{align*}
	It follows from the standard elliptic theory (e.g., \cite[Theorem 2.16]{GGS}) that $u_k$ is uniformly bounded in $H^{2}(\mathbb{B}^n)$.
	Hence, up to a subsequence, as $k \to \infty$, $u_k  \rightharpoonup u$ weakly in $H^{2}(\mathbb{B}^n)$, $f_k  \rightarrow 0$ in $H^{1}(\mathbb{S}^{n-1})$; and there exist two nonnegative Borel measures $\lambda$ and $\nu$ such that
	\begin{align*}
		|\Delta u_k|^{2} \ud x \rightharpoonup&~~ |\Delta u|^{2} \ud x + \lambda \qquad\mathrm{as ~~measures~~in~~} \mathbb{B}^n,\\	
	|f_k|^{\frac{2(n-1)}{n-4}} \ud\mu_{\Sp^{n-1}} \rightharpoonup&~~ \nu \qquad~~~\qquad\qquad\mathrm{as~~measures~~in~~} \mathbb{S}^{n-1}.
	\end{align*}
	
It follows from Lemma \ref{Concentration compactness bdry_4th} that there exist countably many points $x_i \in \mathbb{S}^{n-1}$  with $\nu_i=\nu\left(\left\{x_{i}\right\}\right)>0$ such that 
\begin{align*}
	\nu=\sum_{i} \nu_{i} \delta_{x_{i}} \quad \mathrm{~~and~~} \quad	\hat{\lambda} \geqslant \frac{1}{\alpha(n)}\sum_{i}\nu_i^{\frac{n-4}{n-1}} \delta_{x_i}.
\end{align*}
Notice that
\begin{align*}
	\nu\left(\mathbb{S}^{n-1}\right)=1 \quad \mathrm{~~and~~} \quad \hat{\lambda}\left(\mathbb{S}^{n-1}\right)=\lambda \left(\mathbb{B}^{n}\right) \leqslant \frac{1}{\beta_\ve}.
\end{align*}
By definition of weak convergence for measures, we know that
\begin{align*}
	\int_{\mathbb{S}^{n}} p \ud \nu=0 \qquad \mathrm{~~for~~all~~} p \in \mathring{\mathcal{P}}_{m}.
\end{align*}

Let $\theta=(n-4)/(n-1)$ for simplicity.
By definition of $\Theta(m, \theta, n-1)$ and Lemma \ref{Concentration compactness bdry_4th}, we have 
\begin{align*}
	\Theta(m, \theta, n-1) \leqslant \sum_{j} \nu_{j}^{\theta} \leqslant\alpha (n) \hat \lambda\left(\mathbb{S}^{n-1}\right) \leqslant \frac{\alpha(n)}{\beta_\ve}.
\end{align*}
Hence, 
\begin{align*}
	\beta_\ve \leqslant \frac{\alpha(n)}{\Theta(m, \theta, n-1)}.
\end{align*}
Obviously, this contradicts the choice of $\beta_\ve$. 
\end{proof}

\subsection{Almost sharpness}\label{Sect:n>4}

It is worth pointing out that $\alpha(n)/\Theta(m, \frac{n-4}{n-1}, n-1)$ in \eqref{ineq:n>5} is also \emph{almost optimal}.
Such an example appears as the first one to the fourth order Sobolev trace inequality.
\vskip 8pt

\begin{proof}[Proof of Proposition \protect\ref{example:n>4}]
For each $m$, there exist some natural number $N\geqslant N_m(\Sp^{n-1})$ and   $\nu=\sum_{i=1}^N \nu_i \delta_{x_i} \in  \mathcal{M}_{m}^{c}$ for $1 \leqslant i \leqslant N$, such that
$$\Theta(m, \frac{n-4}{n-1}, n-1)=\sum_{i=1}^N \nu_i^{\frac{n-4}{n-1}}.$$	
	 
	 Let $\stackrel{\frown}{xy}$ denote  the geodesic distance between $x$ and $y$ in $\mathbb{S}^{n-1}$, and  $\mathcal{A}_\delta (x_i)$ be the conic annulus as in \eqref{def:Conic_Annulus} for each $x_i \in \mathbb{S}^{n-1}$. Fix $\delta>0$ small enough such that $\overline{\mathcal{A}_{2 \delta}\left(x_{i}\right)} \cap \overline{\mathcal{A}_{2 \delta}\left(x_{j}\right)}=\emptyset$ for $1 \leqslant i<j \leqslant N$.
		
	 For  each $1 \leqslant i \leqslant N$,  let $\chi_i (x)$ be a smooth cut-off function such that $\chi_i(x)=1$ in $\mathcal{A}_{\delta}(x_i)$ and $\chi_i(x)=0$ outside $\mathcal{A}_{2\delta}(x_i)$. Moreover, $\pa_r \chi_i\geqslant 0$  in  $\overline{\Bn}$ and $\pa_r \chi_i= 0$ on $\Sp^{n-1}$. Indeed, one possible way to choose the cut-off function as $\chi_i(x)=\eta_{1}(r)\widetilde{\chi}_{i}(\rho,\theta)$. Here we can take $\widetilde{\chi}_{i}(\rho,\theta)$
to be a nonnegative smooth cut-off function supported on $B_{2\delta}(x_{i})\subset\mathbb{S}^{n-1}$
such that $\widetilde{\chi}_{i}=1$ on $B_{\delta}(x_{i})$ and $\eta_{1}$
to be a non-decreasing smooth cut-off function on $[1-2\delta,1]$ such
that $0\leqslant \eta_1\leqslant 1$ and $\eta_{1}=1$ on $[1-\delta,1]$, $\eta_{1}=0$ on $[0,1-2\delta]$. 
	 
	 For each $x_i$, we use the polar coordinates $x=(r,\rho,\theta)$ to express the flat metric $|\ud x|^2$ in $\overline{\Bn}$ near $x_i$ as
$$|\ud x|^2=\ud r^2 +r^2 \left(\ud\rho^2 + \sin^2\rho g_{\mathbb{S}^{n-2}}\right)$$
for $x=r\xi \in \Bn, \rho=\ud_{\mathbb{S}^{n-1}}(\xi, x_i)=\stackrel{\frown}{\xi x_i}$ and $\theta \in \Sp^{n-2}$. 
	 
	 Under the above coordinates, we define
	 $$\psi_{\ve,i}(r,\rho)=\left[1+\frac{n-4}{2}\frac{2\ve(1-r)}{(\ve+1-r)^2+\rho^2}\right]\left((\ve+1-r)^2+\rho^2\right)^{\frac{4-n}{2}}$$
	for any $0 < \ve <\delta$ with some sufficiently small $\delta$,  and
	\begin{align*}
		\phi_{\ve, i}(x)= \chi_i(x)\psi(r)\psi_{\ve,i}(r,\rho).
		\end{align*}

		Define 
$$ v(x)=\sum_{i=1}^{N}\nu_i^{\frac{n-4}{2(n-1)}}\phi_{\ve, i} (x),$$
then
\begin{align}\label{est2:4th}
		&\int_{\Sp^{n-1}}v^{\frac{2(n-1)}{n-4}}\ud \mu_{\Sp^{n-1}}\no\\
	=&\sum_{i=1}^{N}\nu_i\int_{B_{2\delta}(x_i)}\phi_{\ve, i}^{\frac{2(n-1)}{n-4}}\ud \mu_{\Sp^{n-1}}\no\\
	=&\big(\sum_{i=1}^{N}\nu_i\big) |\Sp^{n-2}|\int_0^{\delta}(\ve^2+\rho^2)^{1-n}\sin^{n-2}\rho \ud \rho+O(\ve^{1-n})\int_{\frac{\delta}{\ve}}^{\frac{2\delta}{\ve}}(1+t^2)^{1-n}t^{n-2}\ud t\no\\
	=&
	(2\ve)^{1-n}\left|\mathbb{S}^{n-1}\right|+O(\ve^{3-n}).
		\end{align}
		For any $p \in \mathring{\mathcal{P}}_{m}$, we have 
\begin{align*}
& \int_{\mathbb{S}^{n-1}} v^{\frac{2(n-1)}{n-4}} p \ud \mu_{\Sp^{n-1}} \no\\
=& \sum_{i=1}^{N} \nu_i\int_{B_{2 \delta}\left(x_{i}\right)} \phi_{\ve, i}^{\frac{2(n-1)}{n-4}}(\xi) p(\xi) \ud \mu_{\Sp^{n-1}} \no\\
=& \sum_{i=1}^{N}\int_{B_{2 \delta}(x_i)} \left[\phi_{\ve, i}^{\frac{2(n-1)}{n-4}}(\stackrel{\frown}{\xi x_{i}}) \nu_i p\left(x_{i}\right) +\phi_{\ve, i}^{\frac{2(n-1)}{n-2}}(\stackrel{\frown}{\xi x_{i}}) O\left(\stackrel{\frown}{\xi x_{i}}^{2}\right)\right] \ud \mu_{\Sp^{n-1}} \no \\
=& O\left(\ve^{3-n}\right),
\end{align*}
where the last equality follows from
\begin{align*}
	\nu \in \mathcal{M}_{m}^{c}\left(\mathbb{S}^{n-1}\right) \quad \Longrightarrow \quad \sum_{i=1}^{N} \nu_i p\left(x_{i}\right)=0.
\end{align*}

 It is shown in \cite[Theorem 2.1 in Chapter IV]{Stein-Weiss} that there exists a basis 
$\{P_{1}, \cdots, P_{L}\}$ of $\mathring{\mathcal{P}}_{m},$
 such that 
 $$\left.p_{1}=P_1\right|_{\mathbb{S}^{n-1}}, \cdots,\left.p_{L}=P_L\right|_{\mathbb{S}^{n-1}}$$
  are spherical harmonics, where $L=n+\sum_{i=2}^{m}(C_{n+i-1}^{n-1}-C_{n+i-3}^{n-1})$. Then  for each $1\leqslant i \leqslant N$, we claim that there exist $\psi_{1}, \cdots, \psi_{L} \in C_{c}^{\infty}\left(\overline{\mathbb{B}^{n}} \backslash \bigcup_{i=1}^{N} \mathcal{A}_{2 \delta}\left(x_{i}\right)\right)$ such that the determinant 
\begin{align}\label{nonsingular:Gram}
	\det\left[\int_{\mathbb{S}^{n-1}} \psi_{j} p_{k} \ud \mu_{\Sp^{n-1}}\right]_{1 \leqslant j, k \leqslant L} \neq 0.
\end{align}
To this end, we can choose a nonzero smooth function $\eta \in C_{c}^{\infty}\left(\overline{\mathbb{B}^{n}} \backslash \bigcup_{i=1}^{N} \mathcal{A}_{2 \delta}\left(x_{i}\right)\right)$ such that $\eta P_{1}, \cdots, \eta P_{L}$ are linearly independent. It follows that the Gram matrix
\begin{align*}
	\left[\int_{\mathbb{S}^{n-1}} \eta^{2} p_{j} p_{k} \ud \mu_{\Sp^{n-1}}\right]_{1 \leqslant j, k \leqslant L}
\end{align*}
is positive definite, then $\psi_{j}=\eta^{2} P_{j}$ satisfies \eqref{nonsingular:Gram}.

The fact \eqref{nonsingular:Gram} enables us to find $\beta_{1}, \cdots, \beta_{L} \in \mathbb{R}$ such that 
\begin{align}\label{vanishing_moments:4th_1}
	\int_{\mathbb{S}^{n-1}}\left(v^{\frac{2(n-1)}{n-4}}+\sum_{j=1}^{L} \beta_{j} \psi_{j}\right) p_{k} \ud \mu_{\Sp^{n-1}}=0 \qquad \forall~ 1 \leqslant k \leqslant L.
\end{align}
Moreover, it follows from \eqref{est2:4th} that for all $1 \leqslant j \leqslant L$, $\beta_{j}=O\left(\ve^{3-n}\right)$ as $\varepsilon \rightarrow 0$. In the following, we shall use $\pa_r=\sum_{i=1}^n x_i \pa_{x_i}$. As a consequence we can find a constant $c_1 >0$ such that 
\begin{align}\label{lbd1:c_1}
	\sum_{j=1}^{L} \beta_{j} \psi_{j}+c_{1} \ve^{3-n} \geqslant \ve^{3-n}
\end{align} 
and
\begin{equation}\label{lbd2:c_1}
\frac{n-1}{2}\left(\sum_{j=1}^{L}\beta_{j}\psi_{j}+c_{1}\varepsilon^{3-n}\right)+\frac{1}{2}\sum_{j=1}^{L}\beta_{j}\frac{\partial\psi_{j}}{\partial r}>0.
\end{equation}

We need further corrections to satisfy higher order moments constraint and Neumann boundary condition. To this end, we define a test function in the form of
\begin{align}\label{def:test_fcn1}
	u^{\frac{2(n-1)}{n-4}}=&v^{\frac{2(n-1)}{n-4}}+\sum_{j=1}^{L} \beta_{j} \psi_{j}+c_{1} \ve^{3-n}+g(x)(1-r^2)\no\\
	:=&u_1(x)+g(x)(1-r^2).
\end{align}
Here $c_1 \in \R_+$ and $g(x)=g(r,\rho,\theta)$ is a smooth function in $\overline{\Bn}$, which are to be determined later. Our goal is to capture the optimal constant and have a good control of higher order terms at the same time. 

\textbf{Step 1.} We need to find some good candidates of $c_1$ and $g$ such that $u$ satisfies the following conditons:
\begin{itemize}
\item[i)]  Neumann boundary condition:
 $$\frac{\pa u}{\pa r}=-\frac{n-4}{2}u \qquad \mathrm{on~~}\Sp^{n-1}.$$
 \item[ii)] $u> 0$ in $\overline{\Bn}$.
 \item[iii)] Vanishing higher order moments constraint:
 \begin{align*}
	\int_{\mathbb{S}^{n-1}}p_j u^{\frac{2(n-1)}{n-4}}\ud \mu_{\Sp^{n-1}}=0, \qquad 1\leqslant j \leqslant L.
\end{align*}
 \end{itemize}
 
\vskip 8pt
To that end, we shall handle them term by term. Keep in mind that Neumann boundary condition is an additional difficulty to Sobolev inequality on closed manifolds, as well as the second order  example. As we shall see that the property \eqref{Neumann_bdry_local_bubble:n>4} of $\phi_\ve$ plays an important role.
\vskip 8pt

i)~~ We may choose the restriction of $g$ on $\Sp^{n-1}$ as 
\begin{align*}
g(1,\rho,\theta)=&\frac{1}{2}\left(\frac{\pa u_1}{\pa r}+(n-1)u_1\right) \\
=&\frac{n-1}{n-4}v^{\frac{n+2}{n-4}}\frac{\partial v}{\partial r}+\frac{1}{2}\sum_{j=1}^{L}\beta_{j}\frac{\partial\psi_{j}}{\partial r}+\frac{n-1}{2}\left(v^{\frac{2(n-1)}{n-4}}+\sum_{j=1}^{L}\beta_{j}\psi_{j}+c_{1}\varepsilon^{3-n}\right)\\
=&\frac{n-1}{n-4}v^{\frac{n+2}{n-4}}\left(\frac{\partial v}{\partial r}+\frac{n-4}{2}v\right)+\frac{1}{2}\sum_{j=1}^{L}\beta_{j}\frac{\partial\psi_{j}}{\partial r}+\frac{n-1}{2}\left(\sum_{j=1}^{L}\beta_{j}\psi_{j}+c_{1}\varepsilon^{3-n}\right)
\end{align*}
to satisfy the Neumann boundary condition
$$\frac{\pa u}{\pa r}=-\frac{n-4}{2}u \qquad \mathrm{on~~}\Sp^{n-1}.$$

ii)~~We claim that $u>0$ on $\Sp^{n-1}$ if  $c_1$ is sufficiently large.

Thanks to the properties of the \emph{``local bubble"}, we divide our discussion into three distinct domains of $\Sp^{n-1}$ to show
\begin{itemize}

\item On $\mathbb{S}^{n-1}\backslash\cup_{i=1}^{N}B_{2\delta}(x_{i})$,
\[
\frac{\partial v}{\partial r}+\frac{n-4}{2}v=0,
\]
then
\begin{equation*}
g(1,\rho,\theta)=\frac{n-1}{2}\left(\sum_{j=1}^{L}\beta_{j}\psi_{j}+c_{1}\varepsilon^{3-n}\right)+\frac{1}{2}\sum_{j=1}^{L}\beta_{j}\frac{\partial\psi_{j}}{\partial r}>0
\end{equation*}
 by virtue of \eqref{lbd2:c_1}.
 
 \item On $\cup_{i=1}^{N}B_{2\delta}\backslash B_{\delta}(x_{i})$,
\[
\frac{\partial v}{\partial r}+\frac{n-4}{2}v=\sum_{i=1}^{N}\frac{\partial\chi_{i}}{\partial r}\psi\psi_{\varepsilon,i},
\]
then we have
\begin{align*}
g(1,\rho,\theta)=&\frac{n-1}{n-4}v^{\frac{n+2}{n-4}}\left(\sum_{i=1}^{N}\frac{\partial\chi_{i}}{\partial r}\psi\psi_{\varepsilon,i}\right)+\frac{n-1}{2}c_{1}\varepsilon^{3-n}\\
=&\frac{n-1}{2}c_{1}\varepsilon^{3-n}>0,
\end{align*}
since $\pa_r \chi_i= 0$ on $\Sp^{n-1}$. 
 
\item On $\cup_{i=1}^{N}B_{\delta}(x_{i})$,
\[
\frac{\partial v}{\partial r}+\frac{n-4}{2}v=0,
\]
then we have
\begin{align*}
g(1,\rho,\theta)=\frac{n-1}{2}c_{1}\varepsilon^{3-n}>0.
\end{align*}
\end{itemize}

In summary,  we obtain $u>0$ on $\Sp^{n-1}$ and $g(1,\rho,\theta)=O(\ve^{3-n})$. This finishes the proof of the claim.

\vskip 8pt
Next we plan to extend $g(1,\rho,\theta)$ to $g(x)$ in the unit ball $\Bn$. 

To simplify calculations, we may choose $g(x)$ in the form of
$$g(x)=g(r,\rho,\theta)=\eta_2(r) g(1,\rho,\theta)\geqslant 0\qquad \mathrm{in~~} \overline \Bn,$$
where $\eta_2(r) \in C^\infty([0,1])$ such that $\eta_2(r)=1$ for $1-2\delta\leqslant r\leqslant1$,
$\eta_2(r)=0$ for $0\leqslant r\leqslant 1-4\delta$.

For future reference, we define
\begin{align*}
U:=u^{\frac{2(n-1)}{n-4}}=&v^{\frac{2(n-1)}{n-4}}+\sum_{j=1}^{L}\beta_{j}\psi_{j}+c_{1}\varepsilon^{3-n}+g(r,\rho,\theta)\left(1-r^{2}\right)\\
\geqslant& \sum_{j=1}^{L}\beta_{j}\psi_{j}+c_{1}\varepsilon^{3-n}\geqslant \ve^{3-n}>0
\end{align*}
by virtue of \eqref{lbd1:c_1}.
This directly implies that $U>0$ in $\overline \Bn$, so does $u$.

\vskip 8pt

\vskip 8pt

iii)~~ Clearly, it follows from \eqref{vanishing_moments:4th_1} and \eqref{def:test_fcn1} that 
\begin{align*}
	\int_{\mathbb{S}^{n-1}}p_j u^{\frac{2(n-1)}{n-4}}\ud \mu_{\Sp^{n-1}}=0, \qquad 1\leqslant j \leqslant L.
\end{align*}

		 \vskip 8pt
		 \textbf{Step 2.} Sharp constant and a good control of higher order terms.
		 \vskip 8pt
		 
		 We are now in a position to estimate the involved terms.

For the first main term, it follows from \eqref{est2:4th} and definition of $u$ in \eqref{def:test_fcn1} that
\begin{equation}\label{mainterm_1}
\|u\|_{L^{\frac{2(n-1)}{n-4}}(\Sp^{n-1})}^2= 2^{4-n} \left|\mathbb{S}^{n-1}\right|^{\frac{n-4}{n-1}}\ve^{4-n}+O(\ve^{6-n}).
\end{equation}

In the following, we use $A\lesssim B$ to denote that there exists a positive constant $C$ independent of $\ve$ such that $|A|\leqslant CB$.

We turn to estimate $\|u\|_{L^2(\Sp^{n-1})}$. On $\Sp^{n-1}$ we have
\begin{align*}
		u^{\frac{2(n-1)}{n-4}} \leqslant v^{\frac{2(n-1)}{n-4}} +  C\ve^{3-n},
\end{align*}
which directly yields
\begin{align*}
	u^2 \leqslant \left(v^{\frac{2(n-1)}{n-4}} + C \ve^{3-n}\right)^{\frac{n-4}{n-1}} \leqslant v^2 + C  \ve^{\frac{(3-n)(n-4)}{n-1}}. 
\end{align*}
Hence,  we obtain
\begin{align}\label{bubble:est1}
	\int_{\mathbb{S}^{n-1}}u^2 \ud \mu_{\Sp^{n-1}}\leqslant &\int_{\mathbb{S}^{n-1}}\left(v^2 + C\ve^{\frac{(3-n)(n-4)}{n-1}}\right) \ud \mu_{\Sp^{n-1}}\no\\
	 \lesssim & \sum_{i=1}^{N}\nu_i^{\frac{n-4}{n-1}}\int_{B_{2\delta}(x_i)} \phi_{\ve, i}^2 \ud \mu_{\Sp^{n-1}} + C\ve^{\frac{(3-n)(n-4)}{n-1}}\no\\
	 \lesssim&\varepsilon^{(3-n)\frac{n-4}{n-1}}+\sum_{i=1}^{N}\int_{0}^{2\delta}\left(\varepsilon^{2}+\rho^{2}\right)^{4-n}\rho^{n-2}\ud\rho\no\\
  \lesssim&\varepsilon^{\frac{(3-n)(n-4)}{n-1}}.
\end{align}

To estimate $\|\nabla u\|_{L^2(\Sp^{n-1})}^2$. It follows from \eqref{def:test_fcn1} that
\[
\pa_\rho u|_{\Sp^{n-1}}=U^{\frac{n-4}{2(n-1)}-1}v^{\frac{2(n-1)}{n-4}-1}\pa_\rho v+\frac{n-4}{2(n-1)}U^{\frac{n-4}{2(n-1)}-1}\sum_{j=1}^{L}\beta_{j}\pa_\rho\psi_j
\]
and
$$\nabla_\theta u|_{\Sp^{n-1}}=U^{\frac{n-4}{2(n-1)}-1}v^{\frac{2(n-1)}{n-4}-1}\nabla_\theta v+\frac{n-4}{2(n-1)}U^{\frac{n-4}{2(n-1)}-1}\sum_{j=1}^{L}\beta_{j}\nabla_\theta \psi_j.$$

Notice that
\[
|\nabla u|_{\mathbb{S}^{n-1}}^{2}=|\pa_\rho u|^2+\frac{1}{\sin^{2}\rho}\left|\nabla u\right|_{\mathbb{S}^{n-2}}^{2}.
\]

Since
\begin{align*}
\int_{B_{2\delta}\left(x_{i}\right)}|\pa_\rho u|^2\ud \mu_{\Sp^{n-1}}  =&\int_{B_{2\delta}(x_{i})}\left(U^{\frac{n-4}{2(n-1)}-1}v^{\frac{2(n-1)}{n-4}-1}\pa _{\rho} v\right)^{2}\ud\mu_{\Sp^{n-1}}\\
  \lesssim& \int_{B_{2\delta}(x_{i})}|\pa_{\rho}\chi_i  \psi_{\ve, i}+ \chi_i \pa_\rho \psi_{\ve, i}|^{2}\rho^{n-2}\ud\rho\\
 \lesssim&\int_{B_{2\delta}(x_{i})}\left[\left(\varepsilon^{2}+\rho^{2}\right)^{4-n}+\rho^2\left(\varepsilon^{2}+\rho^{2}\right)^{2(\frac{4-n}{2}-1)}\right]\cdot\rho^{n-2}\ud\rho\\
 \lesssim&\varepsilon^{2(2-n)+n+1}\cdot\int_{0}^{\frac{\delta}{\varepsilon}}\left(1+t^{2}\right)^{2-n}t^{n}\ud t\\
  \lesssim&  \ve^{5-n}\begin{cases}
O(1),& n\geqslant 6\\
O(\log\frac{1}{\varepsilon}), & n=5
\end{cases}  \lesssim  \ve^{5-n}\log \frac{1}{\varepsilon}
\end{align*}
and
\begin{align*}
&\int_{\Sp^{n-1}}\frac{1}{\sin^{2}\rho}\left|\nabla u\right|_{\mathbb{S}^{n-2}}^{2} \ud \mu_{\Sp^{n-1}}\\
 \lesssim& \sum_{i=1}^{n+1}\int_{B_{2\delta}(x_{i})}\frac{|\nabla \chi_i|_{\Sp^{n-2}}^2 |\psi_{\ve,i}|^2}{\sin^{2}\rho}\ud \mu_{\Sp^{n-1}} +\int_{\Sp^{n-1}}\frac{1}{\sin^{2}\rho} U^{\frac{n-4}{n-1}-2}\sum_{j=1}^{L}|\beta_{j}|^2|\nabla \psi_j|_{\Sp^{n-2}}^2 \ud \mu_{\Sp^{n-1}}\\
 \lesssim& \sum_{i=1}^{n+1}\int_0^{2\delta}(\ve^2+\rho^2)^{4-n}\rho^{n-4}\ud \rho +\varepsilon^{\frac{(3-n)(n-4)}{n-1}}\int_0^\pi \sin^{n-4} \rho \ud \rho\\
 \lesssim& \varepsilon^{\frac{(3-n)(n-4)}{n-1}}
 +\ve^{5-n}\log \frac{1}{\varepsilon},
\end{align*}
we obtain
\begin{equation}\label{bubble:est2}
\int_{\mathbb{S}^{n-1}}|\nabla u|_{\mathbb{S}^{n-1}}^{2}\ud \mu_{\Sp^{n-1}}=O(\varepsilon^{\frac{(3-n)(n-4)}{n-1}})+O(\ve^{5-n}\log \frac{1}{\varepsilon})=o(\ve^{4-n}).
\end{equation}

		 It is left to calculate $\|\Delta u\|_{L^2(\Bn)}^2$ containing the second main term.
		 
		 Notice that
		 \begin{align*}
\Delta u & =\operatorname{div}\left(\frac{n-4}{2(n-1)}(U)^{\frac{n-4}{2(n-1)}-1}\nabla U\right)\\
 & =\frac{n-4}{2(n-1)}U^{\frac{n-4}{2(n-1)}-1}\Delta U+\frac{n-4}{2(n-1)}\left(\frac{n-4}{2(n-1)}-1\right)U^{\frac{n-4}{2(n-1)}-2}|\nabla U|^{2}.
\end{align*}

		 We divide $\Bn$ into three distinct domains accordingly.

\begin{itemize}

\item On $\overline{\Bn}\backslash\mathcal{A}_{2\delta}(x_{i})$, 
\begin{align*}
U=\sum_{j=1}^{L}\beta_{j}\psi_{j}+c_{1}\varepsilon^{3-n}+\eta(r)\left[\frac{n-1}{2}\left(\sum_{j=1}^{L}\beta_{j}\psi_{j}+c_{1}\varepsilon^{3-n}\right)+\frac{1}{2}\sum_{j=1}^{L}\beta_{j}\frac{\partial\psi_{j}}{\partial r}\right]\left(1-r^{2}\right).
\end{align*}
A direct computation shows $|\nabla U|\lesssim\varepsilon^{3-n}$ and $|\Delta U|\lesssim\varepsilon^{3-n}$. Thus, we obtain
\[
|\Delta u|\lesssim\varepsilon^{(3-n)\frac{(n-4)}{2(n-1)}}
\]

and 
\[
\int_{\Bn\backslash\overline{\mathcal{A}_{2\delta}(x_{i})}}|\Delta u|^{2}\ud x\lesssim\varepsilon^{(3-n)\frac{(n-4)}{(n-1)}}.
\]

\item On $\overline{\mathcal{A}_{2\delta}}\backslash\mathcal{A}_{\delta}(x_{i})$, according to the selection of cut-off functions $\chi_i$ we have
\begin{align*}
U=\nu_i \phi_{\varepsilon,i}^{\frac{2(n-1)}{n-4}}+c_{1}\varepsilon^{3-n}+\frac{n-1}{2}c_{1}\varepsilon^{3-n}\left(1-r^{2}\right).
\end{align*}
A direct computation yields
\begin{align}\label{psi_r}
\partial_{r}\psi_{\varepsilon,i} =\frac{(n-4)(1-r)}{((\varepsilon+1-r)^{2}+\rho^{2})^{\frac{n}{2}}}\left[(\ve+1-r)^2+\rho^2+(n-2)\ve(\ve+1-r)\right]
\end{align}
and 
\begin{align} 
 \partial_{\rho}\psi_{\varepsilon,i}=-\frac{(n-4)\rho}{((\varepsilon+1-r)^{2}+\rho^{2})^{\frac{n}{2}}}\left[(\ve+1-r)^2+\rho^2+(n-2)\ve(1-r)\right], \label{psi_rho}
 \end{align}
\begin{align}
 \partial_\rho^{2}\psi_{\varepsilon,i} =&-\frac{n-4}{((\varepsilon+1-r)^{2}+\rho^{2})^{\frac{n}{2}+1}}\left[(n-1)(\varepsilon+1-r)^{2}+\rho^{2})^2-n(n-2) \rho^2\ve(1-r)\right.\no\\
 &\left.\hskip 128pt -2(n-2)\rho^2((\varepsilon+1-r)^{2}+\rho^{2})\right]
 .\label{psi_rho^2}
\end{align}
Then
\begin{align}\label{est:|D_bubble|^2}
|\nabla\psi_{\varepsilon,i}|^{2} & =\left(\partial_r\psi_{\varepsilon,i}\right)^{2}+r^{-2}\left(\partial_\rho\psi_{\varepsilon,i}\right)^{2}\no\\
=&\left\{(1-r)^2\left[(n-1)\ve^2+\rho^2+(1-r)^2+n\ve(1-r)\right]^2\right.\no\\
&\left.~~+\frac{\rho^2}{r^2}\left[\ve^2+(1-r)^2+n\ve(1-r)+\rho^2\right]^2\right\}\frac{(n-4)^2}{(\varepsilon+1-r)^{2}+\rho^{2})^n}\\
=&O(1).\no
\end{align}
and
\begin{align*}
\Delta\psi_{\varepsilon,i}=&r^{1-n}\pa_r\left(r^{n-1}\pa_r \psi_{\varepsilon,i}\right)+r^{-2}\sin^{2-n}\rho \pa_\rho(\sin^{n-2}\rho \pa_\rho \psi_{\varepsilon,i}) \\
=&\partial_r^{2}\psi_{\varepsilon,i}+\frac{n-1}{r}\partial_r\psi_{\varepsilon,i}+r^{-2}\left(\partial_\rho^{2}\psi_{\varepsilon,i}+(n-2)\cot\rho\partial_\rho\psi_{\varepsilon,i}\right).
\end{align*}
Thus, we have 
 
\[
|\nabla U|\lesssim\varepsilon^{3-n} \qquad \mathrm{and~~}\quad |\Delta U|\lesssim\varepsilon^{3-n}.
\]
Recall that $\phi_{\ve, i}= \chi_i\psi\psi_{\ve,i}$. Hence, putting these facts together we obtain
\[
\int_{\mathcal{A}_{2\delta}\backslash\overline{\mathcal{A}_{\delta}(x_{i})}}|\Delta u|^{2}\ud x\lesssim\varepsilon^{(3-n)\frac{(n-4)}{(n-1)}}.
\]

\item On $\mathcal{A}_{\delta}(x_{i})$,
$$U=\nu_i(\psi\psi_{\varepsilon,i})^{\frac{2(n-1)}{n-4}}+c_{1}\varepsilon^{3-n}+\frac{n-1}{2}c_{1}\varepsilon^{3-n}\left(1-r^{2}\right).$$
Note that
\[
\nabla\left(\psi\psi_{\varepsilon,i}\right)^{\frac{2(n-1)}{n-4}}=\frac{2(n-1)}{n-4}\left(\psi\psi_{\varepsilon,i}\right)^{\frac{2(n-1)}{n-4}-1}\nabla\left(\psi\psi_{\varepsilon,i}\right).
\]

All  terms   $|\nabla\psi_{\varepsilon,i}|^{2}$  involved in the expression of $\Delta u$ are
\begin{align*}
 & \left(\frac{2(n-1)}{n-4}-1\right)U^{\frac{n-4}{2(n-1)}-2}\left(\psi\psi_{\varepsilon,i}\right)^{\frac{2(n-1)}{n-4}-2}\left|\nabla\psi_{\varepsilon,i}\right|^{2}\psi^{2}\left(U-\nu_i\left(\psi\psi_{\varepsilon,i}\right)^{\frac{2(n-1)}{n-4}}\right)\\
= & \left(\frac{2(n-1)}{n-4}-1\right)U^{\frac{n-4}{2(n-1)}-2}\left(\psi\psi_{\varepsilon,i}\right)^{\frac{2(n-1)}{n-4}-2}\left|\nabla\psi_{\varepsilon,i}\right|^{2}\psi^{2}c_{1}\varepsilon^{3-n}\left(1+\frac{n-1}{2}\left(1-r^{2}\right)\right).
\end{align*}

Thus, we obtain
\begin{align*}
|\Delta u|\leqslant& \frac{n-4}{2(n-1)}U^{\frac{n-4}{2(n-1)}-1}\left|\Delta\left(\frac{n-1}{2}c_{1}\varepsilon^{3-n}\left(1-r^{2}\right)\right)\right|\\
  & +\frac{n-4}{2(n-1)}\left(1-\frac{n-4}{2(n-1)}\right)U^{\frac{n-4}{2(n-1)}-2}\left|\nabla\left(\frac{n-1}{2}c_{1}\varepsilon^{3-n}\left(1-r^{2}\right)\right)\right|^{2}\\
  &+\frac{n-1}{2}\nu_i c_1\ve^{3-n}\left(1-\frac{n-4}{2(n-1)}\right)U^{\frac{n-4}{2(n-1)}-2}(\psi\psi_{\varepsilon,i})^{2\frac{n-1}{n-4}-1}\left|\langle\nabla(\psi\psi_{\varepsilon,i}),\nabla r^2\rangle\right|\\
   & +\frac{n-4}{2(n-1)}\left(1-\frac{n-4}{2(n-1)}\right)\nu_i^2 U^{\frac{n-4}{2(n-1)}-2}|\nabla\psi|^{2}\psi_{\varepsilon,i}^{2}\left(\frac{2(n-1)}{n-4}\left(\psi\psi_{\varepsilon,i}\right)^{\frac{2(n-1)}{n-4}-1}\right)^{2}\\
   & +\left(\frac{2(n-1)}{n-4}-1\right)\nu_i U^{\frac{n-4}{2(n-1)}-1}\left(\psi\psi_{\varepsilon,i}\right)^{\frac{2(n-1)}{n-4}-2}|\nabla\psi|^{2}\psi_{\varepsilon,i}^{2}\\
   &+C U^{\frac{n-4}{2(n-1)}-1}\left(\psi\psi_{\varepsilon,i}\right)^{\frac{2(n-1)}{n-4}-1}|\langle\nabla\psi,\nabla\psi_{\varepsilon,i}\rangle|\\
    & +C U^{\frac{n-4}{2(n-1)}-2}\left(\psi\psi_{\varepsilon,i}\right)^{\frac{4(n-1)}{n-4}-2}\psi\psi_{\varepsilon,i}|\langle\nabla\psi,\nabla\psi_{\varepsilon,i}\rangle|\\
    &+ \frac{n+2}{n-4}\nu_i U^{\frac{n-4}{2(n-1)}-2}\left(\psi\psi_{\varepsilon,i}\right)^{\frac{2(n-1)}{n-4}-2}\left|\nabla\psi_{\varepsilon,i}\right|^{2}\psi^{2}c_{1}\varepsilon^{3-n}\left(1+\frac{n-1}{2}\left(1-r^{2}\right)\right)\\
    & +U^{\frac{n-4}{2(n-1)}-1}\nu_i \left(\psi\psi_{\varepsilon,i}\right)^{\frac{2(n-1)}{n-4}-1}|\Delta\left(\psi\psi_{\varepsilon,i}\right)|.
\end{align*}

We first deal with two easier terms:
\begin{align*}
 & \int_{\mathcal{A}_{\delta}(x_{i})}\left[U^{\frac{n-4}{2(n-1)}-1}\Delta\left(\frac{n-1}{2}c_{1}\varepsilon^{3-n}\left(1-r^{2}\right)\right)\right]^{2}\ud x\\
\lesssim & \varepsilon^{2(3-n)}\int_{\mathcal{A}_{\delta}(x_{i})}U^{\frac{n-4}{n-1}-2}\ud x\\
\lesssim & \varepsilon^{2(3-n)}\varepsilon^{(3-n)(\frac{n-4}{n-1}-2)}|\mathcal{A}_{\delta}(x_{i})|\\
\lesssim & \varepsilon^{(3-n)\frac{n-4}{n-1}}
\end{align*}
and
\begin{align*}
 & \int_{\mathcal{A}_{\delta}(x_{i})}\left[U^{\frac{n-4}{2(n-1)}-2}\left|\nabla\left(\frac{n-1}{2}c_{1}\varepsilon^{3-n}\left(1-r^{2}\right)\right)\right|^{2}\right]^{2}\ud x\\
\lesssim & \varepsilon^{4(3-n)}\int_{\mathcal{A}_{\delta}(x_{i})}U^{\frac{n-4}{(n-1)}-4}\ud x\\
\lesssim & \varepsilon^{(3-n)\frac{n-4}{n-1}}.
\end{align*}

In the following, we shall use the change of variables: $s=(1-r)/\ve, t=\rho/\ve$ and $\tau=t/(1+s)$. 

By \eqref{psi_r} we have
\begin{align*}
&\int_{\mathcal{A}_{\delta}(x_{i})}\left[c_1\ve^{3-n}U^{\frac{n-4}{2(n-1)}-2}(\psi\psi_{\varepsilon,i})^{2\frac{n-1}{n-4}-1}\left|\langle\nabla(\psi\psi_{\varepsilon,i}),\nabla r^2\rangle\right|\right]^2\ud x\\
\lesssim & \int_{\mathcal{A}_{\delta}(x_{i})}\left[ \left|\langle\nabla(\psi\psi_{\varepsilon,i}),\nabla r^2\rangle\right|\right]^2\ud x\\
\lesssim &  \int_{\mathcal{A}_{\delta}(x_{i})}|\partial_r \psi_{\varepsilon,i}|^2\ud x\\
\lesssim & \int_{1-\delta}^1\int_0^\delta \left((\varepsilon+1-r)^{2}+\rho^{2}\right)^{-n}(1-r)^2\\
&\qquad\quad \cdot \left[(\ve+1-r)^2+\rho^2+(n-2)\ve(\ve+1-r)\right]^2 \rho^{n-2}\ud \rho \ud r\\
\lesssim & \varepsilon^{6-n}\int_{0}^{\delta/\varepsilon}\int_{0}^{\delta/\varepsilon}((1+s)^{2}+t^{2})^{-n}s^2\left[(1+s)^2+t^2+(n-2)(1+s)\right]^2t^{n-2}\ud s \ud t\\
\lesssim & \varepsilon^{6-n}.
\end{align*}

For a real number $\alpha\leqslant 0$, we have
\begin{align*}
 & \int_{\mathcal{A}_{\delta}(x_{i})}\psi_{\varepsilon,i}^{2}\ud x\\
\lesssim&\int_{\mathcal{A}_{\delta}(x_{i})}\left(1+\frac{n-4}{2}\frac{2\varepsilon(1-r)}{(\varepsilon+1-r)^{2}+\rho^{2}}\right)^2\left((\varepsilon+1-r)^{2}+\rho^{2}\right)^{4-n}\ud x\\
\lesssim&\delta^{-2\alpha} \int_{1-\delta}^\delta \int_0^\delta \left((\varepsilon+1-r)^{2}+\rho^{2}\right)^{4-n+\alpha} \rho^{n-2}\ud \rho\ud r\\
\lesssim& \varepsilon^{8-n+2\alpha}\int_{0}^{\delta/\varepsilon}\int_{0}^{\delta/\varepsilon}\left((1+s)^{2}+t^{2}\right)^{4-n+\alpha}t^{n-2}\ud s \ud t\\
 \lesssim&\varepsilon^{8-n+2\alpha}\int_{0}^{\delta/\varepsilon}(1+s)^{7-n+2\alpha}\ud s\int_{0}^{\delta/\varepsilon}(1+\tau^2)^{4-n+\alpha}\tau^{n-2} \ud \tau.
\end{align*}
We now choose $-\alpha=\frac{8-n}{2}+\e_0$ for some $\e_0\in (0,1/2)$ if $5 \leqslant n \leqslant 8$ and $\alpha=0$ if $n\geqslant 9$ to obtain
\begin{align*}
  \int_{\mathcal{A}_{\delta}(x_{i})}\psi_{\varepsilon,i}^{2}\ud x
 \lesssim\begin{cases}
 \ve^{-2\e_0} \quad & \mathrm{if~~}5 \leqslant n \leqslant 8;\\
 \ve^{8-n} \quad & \mathrm{if ~}n \geqslant 9.
 \end{cases}
 \end{align*}
 
 The above estimate will be used to control the next three terms, which are thus of order $o(\ve^{4-n})$. The first term is
 \begin{align*}
 & \int_{\mathcal{A}_{\delta}(x_{i})}\left[U^{\frac{n-4}{2(n-1)}-2}|\nabla\psi|^{2}\psi_{\varepsilon,i}^{2}\left(\frac{2(n-1)}{n-4}\left(\psi\psi_{\varepsilon,i}\right)^{\frac{2(n-1)}{n-4}-1}\right)^{2}\right]^{2}\ud x\\
\lesssim & \int_{\mathcal{A}_{\delta}(x_{i})}U^{\frac{n-4}{n-1}-4}\psi_{\varepsilon,i}^{\frac{8(n-1)}{n-4}}\ud x\\
\lesssim &  \int_{\mathcal{A}_{\delta}(x_{i})}\psi_{\varepsilon,i}^{(\frac{n-4}{n-1}-4)\cdot\frac{2(n-1)}{n-4}+\frac{8(n-1)}{n-4}}\ud x\\
\lesssim & \int_{\mathcal{A}_{\delta}(x_{i})}\psi_{\varepsilon,i}^{2}\ud x=o(\ve^{4-n})
\end{align*}
and the left two terms are
\begin{align*}
 & \int_{\mathcal{A}_{\delta}(x_{i})}\left[U^{\frac{n-4}{2(n-1)}-1}\left(\frac{2(n-1)}{n-4}-1\right)\left(\psi\psi_{\varepsilon,i}\right)^{\frac{2(n-1)}{n-4}-2}|\nabla\psi|^{2}\psi_{\varepsilon,i}^{2}\right]^{2}\ud x\\
\lesssim & \int_{\mathcal{A}_{\delta}(x_{i})}\psi_{\varepsilon,i}^{\left(\frac{2(n-1)}{n-4}(\frac{n-4}{2(n-1)}-1)+\frac{2(n-1)}{n-4}-2+2\right)2}\ud x\\
\lesssim & \int_{\mathcal{A}_{\delta}(x_{i})}\psi_{\varepsilon,i}^{2}\ud x=o(\ve^{4-n})
\end{align*}
and
\begin{align*}
 & \int_{\mathcal{A}_{\delta}(x_{i})}\left[U^{\frac{n-4}{2(n-1)}-1}\left(\psi\psi_{\varepsilon,i}\right)^{\frac{2(n-1)}{n-4}-1}\Delta\psi\psi_{\varepsilon,i}\right]^{2}\ud x\\
\lesssim & \int_{\mathcal{A}_{\delta}(x_{i})}\psi_{\varepsilon,i}^{\frac{4(n-1)}{n-4}(\frac{n-4}{2(n-1)}-1)+\frac{4(n-1)}{n-4}}\ud x\\
\lesssim & \int_{\mathcal{A}_{\delta}(x_{i})}\psi_{\varepsilon,i}^{2}\ud x=o( \varepsilon^{4-n}).
\end{align*}

Similarly, for a real number $\alpha\leqslant 0$, by \eqref{est:|D_bubble|^2} we have
\begin{align*}
 & \varepsilon^{2(3-n)}\int_{\mathcal{A}_{\delta}(x_{i})}\left[U^{\frac{n-4}{2(n-1)}-2}\left(\psi\psi_{\varepsilon,i}\right)^{\frac{2(n-1)}{n-4}-2}\left|\nabla\psi_{\varepsilon,i}\right|^{2}\psi^{2}\right]^{2}\ud x\\
\lesssim & \varepsilon^{2(3-n)}\int_{\mathcal{A}_{\delta}(x_{i})}U^{\frac{n-4}{(n-1)}-4}\left(\psi_{\varepsilon,i}\right)^{\frac{4(n-1)}{n-4}-4}\left|\nabla\psi_{\varepsilon,i}\right|^{4}\ud x\\
\lesssim & \varepsilon^{2(3-n)}\int_{\mathcal{A}_{\delta}(x_{i})}U^{\alpha}U^{\frac{n-4}{(n-1)}-4-\alpha}\left(\psi_{\varepsilon,i}\right)^{\frac{4(n-1)}{n-4}-4}\left|\nabla\psi_{\varepsilon,i}\right|^{4}\ud x\\
\lesssim & \varepsilon^{(2+\alpha)(3-n)}\int_{\mathcal{A}_{\delta}(x_{i})}U^{\frac{n-4}{(n-1)}-4-\alpha}\left(\psi_{\varepsilon,i}\right)^{\frac{4(n-1)}{n-4}-4}\left|\nabla\psi_{\varepsilon,i}\right|^{4}\ud x\\
\lesssim & \varepsilon^{(2+\alpha)(3-n)}\int_{1-\delta}^{1}\ud r\int_{0}^{\delta}\left((\varepsilon+1-r)^{2}+\rho^{2}\right)^{-(n-1)(\frac{n-4}{n-1}-4-\alpha)+\frac{4-n}{2}(\frac{4(n-1)}{n-4}-4)}\\
 & \cdot\left((\varepsilon+1-r)^{2}+\rho^{2}\right)^{2(2-n)}\left[(1-r)^2+\rho^{2}\right]^2\rho^{n-2}\ud\rho \\
\lesssim & \varepsilon^{(2+\alpha)(3-n)}\int_{1-\delta}^{1}\ud r\int_{0}^{\delta}\left((\varepsilon+1-r)^{2}+\rho^{2}\right)^{n+\alpha(n-1)}\rho^{n-2}\ud\rho\\
\lesssim &\varepsilon^{\alpha(n+1)+n+6}\int_{0}^{\delta/\varepsilon}\int_{0}^{\delta/\varepsilon}\left((1+s)^{2}+t^{2}\right)^{n+\alpha(n-1)}t^{n-2}\ud s \ud t\\
\lesssim& \varepsilon^{\alpha(n+1)+n+6}\int_{0}^{\delta/\varepsilon}(1+s)^{2(n+\alpha(n-1))+n-1}\ud s\int_{0}^{\delta/\varepsilon}(1+\tau^2)^{n+\alpha(n-1)}\tau^{n-2}\ud \tau.
\end{align*}
We emphasize that the condition $n\geqslant 5$ has been used to choose 
\begin{align*}
&-2<\alpha<-\frac{3n}{2(n-1)}\\
\Longrightarrow&~~ \alpha(n+1)+n+6>4-n \quad \mathrm{and}\quad 3n+2\alpha(n-1)-1<-1
\end{align*}
such that
\[\int_{0}^{\delta/\varepsilon}(1+s)^{2(n+\alpha(n-1))+n-1}\ud s<\infty\]
and 
$$\int_{0}^{\delta/\varepsilon}(1+\tau^2)^{n+\alpha(n-1)}\tau^{n-2}\ud \tau<\infty.$$
Thus, we obtain
\begin{align*}
 \varepsilon^{2(3-n)}\int_{\mathcal{A}_{\delta}(x_{i})}\left[U^{\frac{n-4}{2(n-1)}-2}\left(\psi\psi_{\varepsilon,i}\right)^{\frac{2(n-1)}{n-4}-2}\left|\nabla\psi_{\varepsilon,i}\right|^{2}\psi^{2}\right]^{2}\ud x=o(\ve^{4-n}).
 \end{align*}

By \eqref{psi_r} we have
\begin{align*}
 & \int_{\mathcal{A}_{\delta}(x_{i})}\left|U^{\frac{n-4}{2(n-1)}-1}\left(\psi\psi_{\varepsilon,i}\right)^{\frac{2(n-1)}{n-4}-1}\langle\nabla\psi,\nabla\psi_{\varepsilon,i}\rangle\right|^{2}\ud x\\
\lesssim & \int_{\mathcal{A}_{\delta}(x_{i})}U^{\frac{n-4}{n-1}-2}\left(\psi\psi_{\varepsilon,i}\right)^{\frac{4(n-1)}{n-4}-2}\left(r\partial_{r}\psi_{\varepsilon,i}\right)^{2}r^{n-1}\sin^{n-2}\rho \ud\rho \ud r\\
\lesssim & \int_{\mathcal{A}_{\delta}(x_{i})}U^{\frac{n-4}{n-1}-2}\left(\psi\psi_{\varepsilon,{i}}\right)^{\frac{4(n-1)}{n-4}-2}\left((\varepsilon+1-r)^{2}+\rho^{2}\right)^{-n}\\
 & \qquad \cdot(1-r)^2\left[(\ve+1-r)^2+\rho^2+(n-2)\ve(\ve+1-r)\right]^2 \rho^{n-2}\ud \rho \ud r
\\
\lesssim & \int_{1-\delta}^1\int_0^\delta \left((\varepsilon+1-r)^{2}+\rho^{2}\right)^{-n}(1-r)^2\\
&\qquad\quad \cdot \left[(\ve+1-r)^2+\rho^2+(n-2)\ve(\ve+1-r)\right]^2 \rho^{n-2}\ud \rho \ud r\\
\lesssim & \varepsilon^{6-n}\int_{0}^{\delta/\varepsilon}\int_{0}^{\delta/\varepsilon}((1+s)^{2}+t^{2})^{-n}s^2\left[(1+s)^2+t^2+(n-2)(1+s)\right]^2t^{n-2}\ud s \ud t\\
\lesssim & \varepsilon^{5-n}.
\end{align*}
Similarly, 
\begin{align*}
 & \int_{\mathcal{A}_{\delta}(x_{i})}\left[U^{\frac{n-4}{2(n-1)}-2}\left(\psi\psi_{\varepsilon,i}\right)^{\frac{4(n-1)}{n-4}-2}\psi\psi_{\varepsilon,i}\right]^{2}|\langle\nabla\psi,\nabla\psi_{\varepsilon,i}\rangle|^{2}\ud x\\
\lesssim & \int_{\mathcal{A}_{\delta}(x_{i})}\left[(\psi\psi_{\varepsilon,i})^{\frac{2(n-1)}{n-4}(\frac{n-4}{2(n-1)}-2)}\left(\psi\psi_{\varepsilon,i}\right)^{\frac{4(n-1)}{n-4}-2}\psi\psi_{\varepsilon,i}\right]^{2}\left((\varepsilon+1-r)^{2}+\rho^{2}\right)^{-n}\\
 & \qquad\quad \cdot(1-r)^2\left[(\ve+1-r)^2+\rho^2+(n-2)\ve(\ve+1-r)\right]^2 \rho^{n-2}\ud \rho \ud r\\
\lesssim & \int_{1-\delta}^1\int_0^\delta \left((\varepsilon+1-r)^{2}+\rho^{2}\right)^{-n}(1-r)^2\\
&\qquad\quad \cdot \left[(\ve+1-r)^2+\rho^2+(n-2)\ve(\ve+1-r)\right]^2 \rho^{n-2}\ud \rho \ud r\\
\lesssim & \varepsilon^{5-n}.
\end{align*}

The remaining term is exactly the second main term
\[
\nu_i^2 \int_{\mathcal{A}_{\delta}(x_{i})}\left[ U^{\frac{n-4}{2(n-1)}-1}\left(\psi\psi_{\varepsilon,i}\right)^{\frac{2(n-1)}{n-4}-1}\psi\Delta\psi_{\varepsilon,i}\right]^{2}\ud x.
\]

Keep in mind that
\begin{align*}
U & =\nu_i (\psi\psi_{\varepsilon,i})^{\frac{2(n-1)}{n-4}}+c_{1}\varepsilon^{3-n}+\frac{n-1}{2}c_{1}\varepsilon^{3-n}\cdot\left(1-r^{2}\right).
\end{align*}

Observe that
\begin{align*}
 & \nu_i^2 \int_{\mathcal{A}_{\delta}(x_{i})}\left[U^{\frac{n-4}{2(n-1)}-1}\left(\psi\psi_{\varepsilon,i}\right)^{\frac{2(n-1)}{n-4}-1}\psi\Delta\psi_{\varepsilon,i}\right]^{2}\ud x\\
\leqslant & \nu_i^{\frac{n-4}{n-1}}\int_{\mathcal{A}_{\delta}(x_{i})}(\psi\psi_{\varepsilon,i})^{\frac{2(n-1)}{n-4}(\frac{n-4}{2(n-1)}-1)+\frac{2(n-1)}{n-4}-1}(\psi\Delta\psi_{\varepsilon,i})^{2}\ud x\\
\leqslant & \nu_i^{\frac{n-4}{n-1}}\int_{\mathcal{A}_{\delta}(x_{i})}(\psi\Delta\psi_{\varepsilon,i})^{2}\ud x\\
= & \nu_i^{\frac{n-4}{n-1}} \int_{\mathcal{A}_{\delta}(x_{i})}(\Delta\psi_{\varepsilon,i})^{2}\ud x+O(1)\int_{\mathcal{A}_{\delta}(x_{i})}(1-r^{2})(\Delta\psi_{\varepsilon,i})^{2}\ud x.
\end{align*}

Notice that
\begin{align*}
\psi_{\varepsilon,i} (r,\rho)=&\varepsilon^{4-n}\left[1+\frac{n-4}{2}\frac{\frac{2(1-r)}{\varepsilon}}{(1+\frac{1-r}{\varepsilon})^{2}+(\frac{\rho}{\varepsilon})^{2}}\right]\left((1+\frac{1-r}{\varepsilon})^{2}+(\frac{\rho}{\varepsilon})^{2}\right)^{\frac{4-n}{2}}\\
 =&\varepsilon^{4-n}\hat{\psi}(\frac{1-r}{\varepsilon},\frac{\rho}{\varepsilon})
\end{align*}
by recalling that $s=(1-r)/\ve, t=\rho/\ve$ and
\[
\hat{\psi}(s,t)=\left[1+\frac{n-4}{2}\frac{2s}{(1+s)^{2}+t{}^{2}}\right]((1+s)^{2}+t^{2})^{\frac{4-n}{2}}.
\]

Let us deal with 
\begin{align*}
 \Delta\psi_{\varepsilon,i}=&r^{1-n}\pa_r\left(r^{n-1}\pa_r \psi_{\varepsilon,i}\right)+r^{-2}\sin^{2-n}\rho \pa_\rho(\sin^{n-2}\rho \pa_\rho \psi_{\varepsilon,i}) \\
=&\partial_r^{2}\psi_{\varepsilon,i}+\frac{n-1}{r}\partial_r\psi_{\varepsilon,i}+r^{-2}\left(\partial_\rho^{2}\psi_{\varepsilon,i}+(n-2)\cot\rho\partial_\rho\psi_{\varepsilon,i}\right)\\
= & \partial_r^{2}\psi_{\varepsilon,i}+\partial_\rho^{2}\psi_{\varepsilon,i}+\frac{n-2}{\rho}\partial_\rho\psi_{\varepsilon,i}+\frac{n-1}{r}\partial_r\psi_{\varepsilon,i}+(r^{-2}-1)\partial_\rho^{2}\psi_{\varepsilon,i}\\
 & +(n-2)(r^{-2} \cot\rho -\rho^{-1})\partial_\rho\psi_{\varepsilon,i}\\
= & \varepsilon^{2-n}\left(\partial_s^{2}\hat{\psi}+\partial_t^{2}\hat{\psi}+\frac{n-2}{t}\partial_t\hat{\psi}\right)\\
 & +\frac{n-1}{r}\partial_r\psi_{\varepsilon,i}+(r^{-2}-1)\partial_\rho^{2}\psi_{\varepsilon,i} +(n-2)(r^{-2} \cot\rho -\rho^{-1})\partial_\rho\psi_{\varepsilon,i}\\
 =&\varepsilon^{2-n}\Delta \hat \psi+\frac{n-1}{r}\partial_r\psi_{\varepsilon,i}+(r^{-2}-1)\partial_\rho^{2}\psi_{\varepsilon,i} +(n-2)(r^{-2} \cot\rho -\rho^{-1})\partial_\rho\psi_{\varepsilon,i}.
\end{align*}
Here we regard $\Delta \hat \psi$ as the Laplacian of a function $\hat \psi(s,t)=\hat \psi(z)$ with $s=z_n,t=|z'|$ for $z=(z',z_n) \in \Rn_+$.

By \eqref{psi_r} we have
\begin{align*}
 & \int_{\mathcal{A}_{\delta}(x_{i})}\left(\frac{n-1}{r}\partial_r\psi_{\varepsilon,i}\right)^{2}\ud x\\
\lesssim & \int_{1-\delta}^1\int_0^\delta \frac{(1-r)^2 r^{n-1}\rho^{n-2}}{(\varepsilon+1-r)^{2}+\rho^{2})^{n}}\left[(\ve+1-r)^2+\rho^2+(n-2)\ve(\ve+1-r)\right]^{2} \ud r \ud \rho \\
\lesssim &\varepsilon^{6-n}\int_0^\delta\int_0^\delta \frac{s^2 t^{n-2}}{((1+s)^2+t^2)^n}\left[(1+s)^2+t^2+(n-2)(1+s)\right]^2 \ud s \ud t\\
\lesssim & \varepsilon^{6-n}\left[\int_0^{\delta/\ve}s^2(1+s)^{3-n}\ud s\int_0^{\delta/\ve}(1+\tau^2)^{2-n}\tau^{n-2} \ud \tau\right.\\
&\qquad~~+2(n-2)\int_0^{\delta/\ve}s^2(1+s)^{2-n}\ud s\int_0^{\delta/\ve}(1+\tau^2)^{1-n}\tau^{n-2} \ud \tau\\
&\qquad~~\left.+(n-2)^2\int_0^{\delta/\ve}s^2(1+s)^{1-n}\ud s\int_0^{\delta/\ve}(1+\tau^2)^{-n}\tau^{n-2} \ud \tau\right]\\
\lesssim&\varepsilon^{5-n}.
\end{align*}

By \eqref{psi_rho^2} and \eqref{psi_rho} we obtain
\begin{align*}
 & \int_{\mathcal{A}_{\delta}(x_{i})}(r^{-2}-1)^2\left(\partial_\rho^{2}\psi_{\varepsilon,i}\right)^{2}\ud x\\
\lesssim &\int_{1-\delta}^1\int_0^\delta\frac{(1-r)^2}{((\ve+1-r)^2+\rho^2)^{n-2}}r^{n-1}\rho^{n-2} \ud r \ud \rho\\
\lesssim&\varepsilon^{6-n}\int_0^{\delta/\ve}(1+s)^{5-n}\ud s\int_0^{\delta/\ve}(1+\tau^2)^{2-n}\tau^{n-2} \ud \tau\\
\lesssim & \varepsilon^{5-n}
\end{align*}
and
\begin{align*}
 & \int_{\mathcal{A}_{\delta}(x_{i})}(n-2)^2(r^{-2} \cot\rho -\rho^{-1})^2\left(\partial_\rho\psi_{\varepsilon,i}\right)^{2}\ud x\\
 \lesssim &\int_{\mathcal{A}_{\delta}(x_{i})}\frac{(1-r)^2}{((\ve+1-r)^2+\rho^2)^{n-2}}r^{n-1}\rho^{n-2} \ud r \ud \rho\\
\lesssim&\varepsilon^{5-n}.
\end{align*}

A direct computation yields
\begin{align*}
 \Delta \hat \psi(s,t)
= -2(n-4)\left((s+1)^2+t^2+(n-2)(1+s)\right)((1+s)^{2}+t^{2})^{-\frac{n}{2}}.
\end{align*}
Then we have 
\begin{align*}
 & \int_{\mathcal{A}_{\delta}(x_{i})}(\Delta\psi_{\varepsilon,i})^{2}\ud x\\
= &|\mathbb{S}^{n-2}| \int_{\mathcal{A}_{\delta}(x_{i})}(\Delta\psi_{\varepsilon,i})^{2}r^{n-1}(\sin\rho)^{n-2}\ud r\ud\rho\\
= & |\mathbb{S}^{n-2}|\int_{\mathcal{A}_{\delta}(x_{i})}(\Delta\psi_{\varepsilon,i})^{2}\rho^{n-2}\ud r\ud\rho+|\mathbb{S}^{n-2}|\int_{\mathcal{A}_{\delta}(x_{i})}(\Delta\psi_{\varepsilon,i})^{2}[r^{n-1}(\sin\rho)^{n-2}-\rho^{n-2}]\ud r\ud\rho\\
= &4(n-4)^2 \varepsilon^{4-n}|\mathbb{S}^{n-2}|\int_{0}^{\delta/\varepsilon}\int_{0}^{\delta/\varepsilon}\left((s+1)^2+t^2+(n-2)(1+s)\right)^2\frac{t^{n-2}\ud s\ud t}{((1+s)^{2}+t^{2})^n}\\
 & +O(\varepsilon^{5-n}).
\end{align*}

Now we deal with the opposite direction of the above main term. Fix $\delta>\ve^{1/2}$, notice that on $\mathcal{A}_{\varepsilon^{1/2}},$
\[
\psi_{\varepsilon,i}\gtrsim\varepsilon^{\frac{4-n}{2}},
\]
whence
\begin{align*}
U =& (\psi\psi_{\varepsilon,i})^{\frac{2(n-1)}{n-4}}+c_{1}\varepsilon^{3-n}+\frac{n-1}{2}c_{1}\varepsilon^{3-n}\cdot\left(1-r^{2}\right)\\
 \leqslant& (\psi\psi_{\varepsilon,i})^{\frac{2(n-1)}{n-4}}(1+O(\varepsilon^{3-n+n-1}))\\
 =&(\psi\psi_{\varepsilon,i})^{\frac{2(n-1)}{n-4}}(1+O(\varepsilon^{2})).
\end{align*}
Then we obtain
\begin{align*}
 & \int_{\mathcal{A}_{\delta}(x_{i})}\left[U^{\frac{n-4}{2(n-1)}-1}\left(\psi\psi_{\varepsilon,i}\right)^{\frac{2(n-1)}{n-4}-1}\psi\Delta\psi_{\varepsilon,i}\right]^{2}\ud x\\
\geqslant & \int_{\mathcal{A}_{\varepsilon^{1/2}}}\left[U^{\frac{n-4}{2(n-1)}-1}\left(\psi\psi_{\varepsilon,i}\right)^{\frac{2(n-1)}{n-4}-1}\psi\Delta\psi_{\varepsilon,i}\right]^{2}\ud x\\
\geqslant& \int_{\mathcal{A}_{\varepsilon^{1/2}}}(1+O(\varepsilon^{2}))^{\frac{n-4}{2(n-1)}-1}[\psi\Delta\psi_{\varepsilon,i}]^{2}\ud x\\
\geqslant & \int_{\mathcal{A}_{\varepsilon^{1/2}}}[\psi\Delta\psi_{\varepsilon,i}]^{2}\ud x+O(\varepsilon^{2})\int_{\mathcal{A}_{\varepsilon^{1/2}}}[\psi\Delta\psi_{\varepsilon,i}]^{2}\ud x.
\end{align*}

The left proof is similar to the above by noting 
\[
\varepsilon^{1/2}/\varepsilon\rightarrow\infty.
\]

Based on the above estimates, it is not hard to see
$$\int_{\mathcal{A}_{\delta}(x_{i})}(1-r^{2})(\Delta\psi_{\varepsilon,i})^{2}\ud x\lesssim \ve^{5-n}.$$

Hence, with the above choice of $\delta$, the main term is 
\begin{align*}
4(n-4)^2 \varepsilon^{4-n}|\mathbb{S}^{n-2}|\int_{0}^{\frac{\delta}{\varepsilon}}\int_{0}^{\frac{\delta}{\varepsilon}}\left((s+1)^2+t^2+(n-2)(1+s)\right)^2((1+s)^{2}+t^{2})^{-n}t^{n-2}\ud s\ud t.
\end{align*}

A direct calculation shows
\begin{align*}
 &4(n-4)^2 \int_{0}^{\infty}\int_{0}^{\infty}\left[(s+1)^2+t^2+(n-2)(1+s)\right]^2((1+s)^{2}+t^{2})^{-n}t^{n-2}\ud s\ud t\\
 =&n(n-2)(n-4)B(\frac{n-1}{2},\frac{n-1}{2}).
 \end{align*}

\end{itemize}

Therefore, putting these facts together we conclude that
\begin{align}\label{mainterm_2}
\int_{\mathbb{B}^{n}}(\Delta u)^{2}\ud x=&|\mathbb{S}^{n-2}|\big(\sum_{i=1}^N \nu_i^{\frac{n-4}{n-1}}\big)n(n-2)(n-4)B(\frac{n-1}{2},\frac{n-1}{2})\varepsilon^{4-n}+o(\varepsilon^{4-n})\no\\
=&2^{2-n}|\mathbb{S}^{n-1}|\Theta(m, \frac{n-4}{n-1}, n-1)n(n-2)(n-4)\varepsilon^{4-n}+o(\varepsilon^{4-n})
\end{align}
by virtue of $2^{n-2}B\left(\frac{n-1}{2},\frac{n-1}{2}\right)|\mathbb{S}^{n-2}|=\left|\mathbb{S}^{n-1}\right|$.

Finally, we combine \eqref{mainterm_1}, \eqref{bubble:est1}, \eqref{bubble:est2} and \eqref{mainterm_2}, as well as other higher order terms, to show
\begin{align*}
&2^{4-n}\left|\mathbb{S}^{n-1}\right|^{\frac{n-4}{n-1}}\ve^{4-n}\leqslant a2^{2-n}|\mathbb{S}^{n-1}|\Theta(m, \frac{n-4}{n-1}, n-1)n(n-2)(n-4)\varepsilon^{4-n}+o(\varepsilon^{4-n}),
\end{align*}
which implies
\begin{align*}
a & \geqslant \frac{4}{(n-2)(n-4)n} \left|\mathbb{S}^{n-1}\right|^{-\frac{3}{n-1}}\cdot\frac{1}{\Theta(m, \frac{n-4}{n-1}, n-1)}= \frac{\alpha(n)}{\Theta(m, \frac{n-4}{n-1}, n-1)}.
\end{align*}
This completes our construction.
\end{proof}

\section{Almost sharp  Sobolev trace inequality of order four under constraints in dimension four}\label{Sect5}

For the Sobolev trace inequality, our conic proof in dimension four can be compared to the one in the pioneering work of Chang-Hang \cite{Chang-Hang}. For the four dimensional example,  as mentioned before, we need to settle the existing obstructions and make an improvement of Chang-Hang type estimate. These enable us to complete the construction.

\subsection{Fourth order Sobolev trace inequality}\label{Subsect:4th_Sobolev_trace_ineq}

Due to the same reason as in Section \ref{subsect:4th order ineq}, it is important to understand Ache-Chang's sharp Sobolev trace inequality in Theorem \ref{Thm:Ache-Chang_n=4} well first. Moreover, notice in the equality case in \eqref{ineq:Ache-Chang_n=4} with $f=0$ that $u(x)=(1-|x|^2)/2$ and thus the extremal metric $e^{2u}|\ud x|^2$ is the \emph{Fefferman-Graham metric} defined in \cite{Fefferman-Graham}.

\vskip 8pt
	
	We employ the regularity theory for bi-Laplace boundary value problem in Gazzola-Grunau-Sweers \cite{GGS} to prove a preliminary result.
	\begin{lemma}\label{lem:H^2-est}
Given $f_i \in C^\infty(\Sp^3)$, let $u_i$ be the biharmonic extension of $ f_i$ to the unit ball $\mathbb{B}^4$ satisfying zero Neumann boundary condition. Assume that as $i \to \infty$, $u_i \rightharpoonup 0$ in $H^2(\mathbb{B}^4)$ and $u_i \to 0$ in $H^1(\Bn)$, $f_i \to 0$ in $H^{3/2}(\mathbb{S}^3).$ Denote by $U_i$ the biharmonic extension of $\varphi f_i$ to $\mathbb{B}^4$, where $\varphi\in C^\infty(\overline{\mathbb{B}^4})$ satisfies $\partial \varphi/\partial r=0$ on $\mathbb{S}^3$. Then there holds
\begin{align*}
    \lim_{i \to \infty}\|\varphi u_i- U_i \|_{H^2(\mathbb{B}^4)}= 0.
\end{align*}
\end{lemma}
\begin{proof}
    By definition of $U_i$, we know that $U_i$ satisfies
    $$
    \begin{cases}
        \Delta^{2} U_i=0 & \mathrm{~~in~~} \mathbb{B}^4,\\ 
        U_i=\varphi f_i &\mathrm{~~on~~}  \mathbb{S}^3,\\ 
        \frac{\pa U_i}{\pa r}=0 & \mathrm {~~on~~}  \mathbb{S}^3.
    \end{cases}
    $$
By assumption, a direct computation yields that $\varphi u_i$ satisfies
    \begin{align*}
        \begin{cases}
        \Delta^{2} \left(\varphi u_i\right)=\Delta^2\varphi u_i + 4\langle \nabla u_i, \nabla \Delta \varphi \rangle + 2 \Delta \varphi   \Delta u_i + 4\langle \nabla^2 \varphi , \nabla^2 u_i \rangle +4 \langle \nabla \varphi, \nabla \Delta u_i\rangle & \mathrm {in~~} \mathbb{B}^4, \\ 
        \varphi u_i =\varphi f_i &\mathrm{on~~}  \mathbb{S}^3, \\
        \frac{\pa }{\pa r}(\varphi u_i)=0 & \mathrm{on~~}  \mathbb{S}^3.
       \end{cases}
    \end{align*}
    Thus, if we let $v_i \in H^2(\Bn)\cap H_0^1(\Bn)$ be a weak solution\footnote{See \cite[(2.42) on p.41 and Theorem 2.31 on p.52]{GGS} for the precise definition of weak solution.}  of 
    \begin{align*}
         \begin{cases}
         \Delta^{2} v_i=\Delta^2 \varphi u_i + 4\langle \nabla u_i, \nabla \Delta \varphi \rangle  & \mathrm {~~in~~} \mathbb{B}^4, \\ 
         v_i=0 &\mathrm {~~on~~}  \mathbb{S}^3, \\ 
         \frac{\pa v_{i}}{\pa r}=0 & \mathrm {~~on~~}  \mathbb{S}^3, 
         \end{cases}
    \end{align*}
    then $\varphi u_i - U_i - v_i$ weakly satisfies 
    \begin{align*}
         \begin{cases}
         \Delta^{2} \left(\varphi u_i - U_i - v_i\right)= 2 \Delta \varphi \Delta u_i + 4\langle \nabla^2 \varphi , \nabla^2 u_i \rangle +4 \langle \nabla \varphi, \nabla \Delta u_i\rangle & \mathrm{~~in~~} \mathbb{B}^4, \\ 
         \varphi u_i - U_i - v_i =0 &\mathrm{~~on~~}  \mathbb{S}^3, \\ 
         \frac{\pa}{\pa r}\left(\varphi u_i - U_i - v_i\right)=0 & \mathrm {~~on~~}  \mathbb{S}^3.
         \end{cases}
    \end{align*}
   Hence, it follows from \cite[Theorem 2.22]{GGS} that 
    \begin{align*}
        \|\varphi u_i - U_i - v_i \|_{H^2(\mathbb{B}^4)}\leqslant C \|u_i \|_{H^1(\mathbb{B}^4)}.
    \end{align*}
    Moreover, by   \cite[Theorem 2.16]{GGS} we have
    \begin{align*}
        \|v_i\|_{H^2(\mathbb{B}^4)}\leqslant C \|u_i\|_{H^1(\mathbb{B}^4)}.
    \end{align*}
    
    Consequently, we conclude that
    \begin{align*}
        & \|\varphi u_i - U_i \|_{H^2(\mathbb{B}^4)}\\
         \leqslant& \|v_i\|_{H^2(\mathbb{B}^4)}+ \|\varphi u_i - U_i - v_i \|_{H^2(\mathbb{B}^4)}\leqslant C \|u_i \|_{H^1(\mathbb{B}^4)}\to 0 \qquad \mathrm{as~~} i \to \infty.
    \end{align*}
    This completes the proof.
\end{proof}

\begin{lemma}\label{lem8.2} 
For $f_{i}\in  C^{\infty}\left( \mathbb{S}^3\right) $ with $\overline{%
f_{i}}=0$,  let $u_i$ be the biharmonic extension of $f_i$ to the unit ball $\mathbb{B}^4$ satisfying zero Neumann boundary condition. We also
assume that as $i\to \infty$, $u_{i}\rightharpoonup u$ weakly in $H^{2}\left( \mathbb{B}^4\right) $, $%
u_{i}\rightarrow u$ a.e. in $\mathbb{B}^4$ and%
\begin{equation}
\left(\Delta u_i\right)^{2} \ud x \rightharpoonup \left(\Delta u\right)^{2} \ud x + \sigma \qquad \mathrm{~~as~~measures.}
\label{eq8.3}
\end{equation}%
If $K\subset \mathbb{S}^3$ is a compact subset with $\sigma \left( \mathcal{C}(K)\right)
<1$, then for any $1\leqslant p<\frac{1}{\sigma \left( \mathcal{C}(K)\right) }$, $%
e^{12\pi^2 f_{i}^{2}}$ is bounded in $L^{p}\left( K,g_{\Sp^3}\right) $, i.e.,%
\begin{equation}
\sup_{i}\int_{K}e^{12\pi^2 p f_{i}^{2}}\ud\mu_{\Sp^3} <\infty .  \label{eq8.4}
\end{equation}
\end{lemma}
\begin{proof}
		Let $%
v_{i}=u_{i}-u$ and $g_i =f_i -f$, then as $i \to \infty$, $v_{i}\rightharpoonup 0$ weakly in $H^{2}\left(
\mathbb{B}^4\right) $, $ v_i\rightarrow 0$ in $H^1\left(\mathbb{B}^4\right)$, $g_i \rightarrow 0$ in $H^{1}\left(\mathbb{S}^3\right)$. Thus, for any $%
\varphi \in C^{\infty }( \overline{\mathbb{B}^4}) $,  as $i\rightarrow \infty$ we obtain
\begin{align*}
&\| \Delta \left( \varphi v_{i}\right) \| _{L^2(\mathbb{B}^4)}^{2} \\
=&\int_{\mathbb{B}^4}\left(\Delta \varphi v_i +2 \langle \nabla \varphi,\nabla v_{i}\rangle+\varphi
 \Delta v_{i}\right)^2 \ud x \\
=&\int_{\mathbb{B}^4}\left\vert \Delta \varphi \right\vert ^{2}v_{i}^{2}\ud x
+4\int_{\mathbb{B}^4}\left(\langle\nabla \varphi,\nabla v_{i}\rangle\right)^2 \ud x
+ \int_{\mathbb{B}^4}\left(\varphi \Delta v_{i}\right)^2 \ud x \\
&+\int_{\mathbb{B}^4}\left( 4 v_i \Delta\varphi\langle\nabla \varphi,\nabla v_{i}\rangle -2\varphi \Delta v_i \langle\nabla u,\nabla v_{i}\rangle+2\varphi\Delta\varphi v_i \Delta v_i\right) \ud x \\
\rightarrow &\int_{\mathbb{B}^4}\varphi ^{2}\ud\sigma.
\end{align*}%

If $1\leqslant p_{1}<\frac{1}{\sigma\left(
\mathcal{C}(K)\right) }$, then $\sigma \left( \mathcal{C}(K)\right) <\frac{1}{p_{1}}$. Hence, there
exists $\varphi \in C^{\infty }(\overline{\mathbb{B}^4}) $ such that $\pa \varphi/\pa r=0$ on $\Sp^3$, $\left. \varphi
\right\vert _{\mathcal{C}(K)}=1$ and $\int_{\mathbb{B}^4}\varphi ^{2}\ud\sigma <\frac{1}{p_{1}}$. 

The fractional Graham, Jenne, Mason, Sparling (GJMS) operator of order three on $\Sp^3$ defined in \cite{Beckner} is
$$P_3^{\Sp^3}=-\Delta_{\Sp^3}(-\Delta_{\Sp^3}+1)^{1/2},$$
which coincides with the one in a Poincar\'e-Einstein manifold introduced by C. Graham and M. Zworski \cite{Graham-Zworski} via scattering theory, and enjoys the conformal covariance property that 
$$P_3^{e^{2\psi}g_{\Sp^3}}=e^{-3\psi}P_3^{\Sp^3} \qquad \mathrm{~~for~~} \psi \in C^\infty(\Sp^3).$$
Denote by $U_i$ the biharmonic function in $\mathbb{B}^4$ of $\varphi g_i$ with zero Neumann boundary condition, then it follows from \cite[(5.5)]{Ache-Chang} that
\begin{equation*}
2\int_{\Sp^3}\varphi g_{i} P_3^{\mathbb{S}^3}\left( \varphi g_{i}\right) \ud \mu_{\Sp^3} =\int_{\mathbb{B}^4}|\Delta U_i|^2 \ud x + 2\int_{\mathbb{S}^3} |\nabla g_{i}|_{\Sp^3}^2 \ud \mu_{\Sp^3}.
\end{equation*}
On the other hand, a direct consequence of Lemma \ref{lem:H^2-est} is
$$\lim_{i \to \infty} \int_{\mathbb{B}^4}|\Delta(U_i-\varphi v_i)|^2 \ud x=0.$$

Therefore, putting these facts together, we obtain that for $i$ sufficiently large, 
$$2\int_{\Sp^3}\varphi g_{i} P_3^{\mathbb{S}^3}\left( \varphi g_{i}\right) \ud \mu_{\Sp^3} <\frac{1}{p_{1}}.$$

We are now ready to estimate
\begin{align*}
\int_{K}e^{12\pi^2 p_{1}\left( g_{i}-\overline{\varphi g_{i}}\right) ^{2}}\ud\mu_{\Sp^3}
\leqslant &\int_{\mathbb{S}^3}e^{12\pi^2 p_{1}\left( \varphi g_{i}-\overline{\varphi g_{i}}%
\right) ^{2}}\ud\mu_{\Sp^3} \\
\leqslant &\int_{\mathbb{S}^3}e^{6\pi^2 \frac{\left( \varphi g_{i}-\overline{\varphi g_{i}}%
\right)^{2}}{\int_{\Sp^3}\varphi g_{i} P_3^{\mathbb{S}^3}\left( \varphi g_{i}\right) \ud \mu_{\Sp^3} }} \ud\mu_{\Sp^3} \\
\leqslant &C,
\end{align*}%
where the second inequality follows from the Moser-Trudinger type inequality for $P_3^{\mathbb{S}^3}$; see Chang-Yang \cite[Proposition 4.4]{Chang-Yang3}.

With the above estimates at hand, we obtain that for any $\varepsilon >0$,%
\begin{align*}
f_{i}^{2} =&\left( \left( g_{i}-\overline{\varphi g_{i}}\right) +f+%
\overline{\varphi g_{i}}\right) ^{2} \\
=&\left( g_{i}-\overline{\varphi g_{i}}\right) ^{2}+2\left( g_{i}-\overline{%
\varphi g_{i}}\right) \left( f+\overline{\varphi g_{i}}\right) +\left( f+%
\overline{\varphi g_{i}}\right) ^{2} \\
\leqslant &\left( 1+\varepsilon \right) \left( g_{i}-\overline{\varphi g_{i}}%
\right) ^{2}+\left( 1+\varepsilon ^{-1}\right) \left( f+\overline{\varphi
g_{i}}\right) ^{2} \\
\leqslant &\left( 1+\varepsilon \right) \left( g_{i}-\overline{\varphi g_{i}}%
\right) ^{2}+2\left( 1+\varepsilon ^{-1}\right) f^{2}+2\left( 1+\varepsilon
^{-1}\right) \overline{\varphi g_{i}}^{2}.
\end{align*}%
Hence,%
\begin{equation*}
e^{12\pi^2 f_{i}^{2}}\leqslant e^{12\pi^2 \left( 1+\varepsilon \right) \left( g_{i}-%
\overline{\varphi g_{i}}\right) ^{2}}e^{24\pi^2 \left( 1+\varepsilon
^{-1}\right) f^{2}}e^{24\pi^2 \left( 1+\varepsilon ^{-1}\right) \overline{%
\varphi g_{i}}^{2}}.
\end{equation*}%

Given $1\leqslant p<\frac{1}{\sigma \left( \mathcal{C}(K)\right) }$, we can choose some $p_{1}\in
\left( p,\frac{1}{\sigma \left( \mathcal{C}(K)\right) }\right) $ and small enough $%
\varepsilon >0$ such that $\frac{p_{1}}{1+\varepsilon }>p$. Notice that $%
e^{12\pi^2 \left( 1+\varepsilon \right) \left( g_{i}-\overline{\varphi g_{i}}%
\right) ^{2}}$ is bounded in $L^{\frac{p_{1}}{1+\varepsilon }}\left(
K\right) $, $e^{24\pi^2 \left( 1+\varepsilon ^{-1}\right) f^{2}}\in L^{q}\left(
K,g_{\Sp^3}\right) $ for any $q>0$ (e.g., see \cite[Lemma 2.1]{Chang-Hang}) and $e^{24\pi^2 \left(
1+\varepsilon ^{-1}\right) \overline{\varphi g_{i}}^{2}}\rightarrow 1$ as $%
i\rightarrow \infty $. Therefore,  by H\"older's inequality we conclude that $e^{12\pi^2
f_{i}^{2}}$ is bounded in $L^{p}\left( K,g_{\Sp^3}\right) $.
\end{proof}

\begin{corollary}
\label{cor8.1}
With the same assumption as in Lemma \ref{lem8.2}, let%
\begin{equation*}
\kappa =\max_{x\in \mathbb{S}^3}\sigma \left( \left\{x\right\} \right) \leqslant 1,
\end{equation*}
then
\begin{enumerate}
\item[i)] if $\kappa <1$, then for any $1\leqslant p<\frac{1}{\kappa }$, $e^{12\pi^2
f_{i}^{2}}$ is bounded in $L^{p}\left( \mathbb{S}^3\right) $. In particular, as $i \to \infty$, $e^{12\pi^2
f_{i}^{2}}\rightarrow e^{12\pi^2 f^{2}}$ in $L^{1}\left( \mathbb{S}^3\right) $;

\item[ii)] if $\kappa =1$, then $\sigma =\delta _{x_0}$ for some $x_{0}\in \mathbb{S}^3$, $%
u=f=0$ and after passing to a subsequence, as $i \to \infty$,%
\begin{equation*}
e^{12\pi^2 f_{i}^{2}} \rightharpoonup 1+c_{0}\delta _{x_{0}} \qquad \mathrm{~~as~~measures,}
\end{equation*}%
for some constant $c_{0}\geqslant 0$.
\end{enumerate}
\end{corollary}
\begin{proof}
Since the proof is similar in spirit to the one of \cite[Corollary 2.1]{Chang-Hang}, we omit it here. 
\end{proof}

\begin{proposition}
\label{prop8.1}
Assume $\alpha >0$, $m_{i}>0$, $m_{i}\rightarrow \infty $. For $
f_{i}\in C^{\infty}\left(\mathbb{S}^3\right)$ with $\overline{f_{i}}=0$ and $\int_{\Sp^3} f_{i} P_3^{\mathbb{S}^3}f_{i} \ud \mu_{\Sp^3} =\left(\int_{\mathbb{B}^4}\left( \Delta  u_{i}\right)^2 \ud x + 2\int_{\Sp^3} |\nabla f_{i}|_{\Sp^3}^2 \ud \mu_{\Sp^3}\right)/2 =1$, where $u_i$ be the biharmonic extension of $f_i$ to the unit ball $\mathbb{B}^4$ satisfying zero Neumann boundary condition, and%
\begin{equation*}
\log \int_{\Sp^3}e^{3m_{i}f_{i}}\ud\mu_{\Sp^3} \geqslant \alpha m_{i}^{2}.
\end{equation*}%
We also assume as $i \to \infty$, $u_{i}\rightharpoonup u$ weakly in $H^{2}\left( \mathbb{B}^4\right) $, $%
\int_{\Sp^3} f_{i} P_3^{\mathbb{S}^3}f_{i} \ud \mu_{\Sp^3} \rightharpoonup \int_{\Sp^3} f P_3^{\mathbb{S}^3}f \ud \mu_{\Sp^3} +\sigma $ as measures and%
\begin{equation*}
\frac{e^{3m_{i}f_{i}}}{\int_{\mathbb{S}^3}e^{3m_{i}f_{i}}\ud\mu_{\Sp^3} }\rightharpoonup \nu \qquad \mathrm{as~~measures}.
\end{equation*}%
Let%
\begin{equation*}
\left\{ x\in \mathbb{S}^3; \sigma \left(\{x\}\right) \geqslant \frac{16}{3}\pi^2 \alpha \right\} =\left\{
x_{1},\cdots ,x_{N}\right\} ,
\end{equation*}%
then%
\begin{equation*}
\nu =\sum_{i=1}^{N}\nu _{i}\delta _{x_{i}},
\end{equation*}%
where $\nu _{i}\geqslant 0$ and $\sum_{i=1}^{N}\nu _{i}=1$.
\end{proposition}

\begin{proof}
First we claim that if $K$ is a compact subset of $\mathbb{S}^3$ with $\sigma \left(
\mathcal{C}(K)\right) <\frac{16}{3}\pi^2 \alpha $, then $\nu \left( K\right) =0$. To this end, we can find
another compact set $K_{1}$ such that $K\subset \mathring{K_1}$, the interior of $K_1$, and $%
\sigma \left(\mathcal{C}(K_{1})\right) <\frac{16}{3}\pi^2 \alpha $. Fix a number $p$ such that%
\begin{equation*}
\frac{3}{16\pi^2 \alpha }<p<\frac{1}{\sigma \left(\mathcal{C}(K_{1})\right) },
\end{equation*}%
then it follows from Lemma \ref{lem8.2} that with a positive constant $C$ independent of $i$, there holds
\begin{equation*}
\int_{K_{1}}e^{12\pi^2 pu_{i}^{2}}\ud\mu_{\Sp^3} \leqslant C.
\end{equation*}%
Notice that%
\begin{equation*}
3m_{i}u_{i}\leqslant 12\pi^2 pu_{i}^{2}+\frac{3m_{i}^{2}}{16\pi^2 p},
\end{equation*}%
this together with Lemma \ref{lem8.2} yields
\begin{equation*}
\int_{K_{1}}e^{3m_{i}u_{i}}\ud\mu_{\Sp^3} \leqslant C e^{\frac{3m_{i}^{2}}{16\pi^2 p}}.
\end{equation*}%
It follows that%
\begin{equation*}
\frac{\int_{K_{1}}e^{4m_{i}u_{i}}\ud\mu_{\Sp^3} }{\int_{\Sp^3}e^{3m_{i}u_{i}}\ud\mu_{\Sp^3} }\leqslant
Ce^{\left( \frac{3}{16\pi^2 p}-\alpha \right) m_{i}^{2}}.
\end{equation*}%
Hence%
\begin{equation*}
\nu \left( K\right) \leqslant \nu \left( \mathring{K_1}\right) \leqslant 
\liminf_{i\rightarrow \infty }\frac{\int_{K_{1}}e^{4m_{i}u_{i}}\ud\mu_{\Sp^3} }{%
\int_{M}e^{4m_{i}u_{i}}\ud\mu_{\Sp^3} }=0,
\end{equation*}%
whence $\nu \left( K\right) =0$.

If $\sigma \left(\{x\}\right) <\frac{16}{3}\pi^2 \alpha $, then we choose $r_{x}>0$ small enough such that $\sigma \left( \mathcal{C}(\overline{B_{r_{x}}\left( x\right) })\right) <\frac{16}{3}\pi^2 \alpha 
$. It follows from the above claim that $\nu \left( \overline{B_{r_{x}}\left(
x\right) }\right) =0$. Hence,%
\begin{equation*}
\nu \left( \mathbb{S}^3\backslash \left\{ x_{1},\cdots ,x_{N}\right\} \right) =0.
\end{equation*}%
In other words, $\nu =\sum_{i=1}^{N}\nu _{i}\delta _{x_{i}}$ with $\nu
_{i}\geqslant 0$ and $\sum_{i=1}^{N}\nu _{i}=1$.
\end{proof}

\begin{proof}[Proof of Theorem \protect\ref{Thm:Sobolev trace n=4}]
Let $\alpha_\ve =\frac{3}{16\pi^2 N_{m}(\Sp^3)}+\varepsilon $. By the proof of  \cite[Theorem A]{Ache-Chang}, we only need to prove the case when $u$ is a biharmonic extension of $f$. 

 By contradiction, if \eqref{ineq:n=4} is not
true, then
there exist some $\ve>0$ and $v_{i}\in H^{2}\left( \mathbb{B}^4\right) $ to be the biharmonic extension of $f_i \in C^\infty(\Sp^3)$,  where $\overline{f_{i}}=0$, $%
 \int_{\mathbb{S}^3}p e^{3f_{i}}\ud\mu_{\Sp^3} =0$ for all $p \in \mathring{\mathcal{P}}_{m}$, such that%
\begin{equation*}
\log \fint_{\mathbb{S}^3}e^{3 f_{i}}\ud\mu_{\Sp^3} -2\alpha_\ve \int_{\Sp^3}  f_{i} P_3^{\mathbb{S}^3}f_{i}\ud \mu_{\Sp^3}\rightarrow \infty  \qquad \mathrm{~~as~~} i\rightarrow \infty.\label{eq8.2}
\end{equation*}%
Then $\log \int_{\mathbb{S}^3}e^{3 f_{i}}\ud\mu_{\Sp^3} \rightarrow \infty $. It follows from \cite[Theorem 1]{Beckner} that
\begin{equation*}
\log \fint_{\mathbb{S}^3}e^{3 f_{i}}\ud\mu_{\Sp^3} \leqslant \frac{3}{8\pi^2 }\int_{\Sp^3}  f_{i} P_3^{\mathbb{S}^3}f_{i}\ud \mu_{\Sp^3},  
\end{equation*}%
whence $\int_{\Sp^3}  f_{i} P_3^{\mathbb{S}^3}f_{i}\ud \mu_{\Sp^3} \rightarrow \infty $. Let 
$$
m_{i}=\left(\int_{\Sp^3}  f_{i} P_3^{\mathbb{S}^3}f_{i}\ud \mu_{\Sp^3} \right)^{\frac{1}{2}}, \quad u_{i}=\frac{v_{i}}{m_{i}} \quad \mathrm{~~and~~} \quad g_{i}=\frac{f_{i}}{m_{i}},$$ 
then  $m_{i}\rightarrow \infty $ and $u_i$ is the biharmonic extension of $g_i$ such that $\int_{\Sp^3} g_{i} P_3^{\mathbb{S}^3} g_{i}\ud \mu_{\Sp^3}=1$, $\overline{g_{i}}=0$. Up to a
subsequence, as $i \to \infty$ we have%
\begin{align*}
u_{i} \rightharpoonup &~u \qquad \mathrm{~~weakly~~in~~}H^{2}\left( \mathbb{B}^4\right), \\
g_{i} \rightarrow &~g \qquad \mathrm{~~in~~}H^{1}\left( \mathbb{S}^3\right),\\
\log \fint_{\mathbb{S}^3}e^{3m_{i} g_{i}}\ud\mu_{\Sp^3} -\alpha_\ve m_{i}^{2} \rightarrow &~\infty , \\
\int_{\mathbb{B}^4} \left(\Delta u_{i}\right)^2 \ud x \rightharpoonup &~\int_{\mathbb{B}^4} \left(\Delta u\right)^2 \ud x +\sigma \quad\mathrm{~~as~~measures}, \\
\frac{e^{3m_{i} g_{i}}}{\int_{\mathbb{S}^3}e^{3m_{i} g_{i}}\ud\mu_{\Sp^3} } \rightharpoonup &~\nu \quad \mathrm{
~~as~~measures}.
\end{align*}

Let%
\begin{equation*}
\left\{ x\in \mathbb{S}^3; \sigma \left( \{x\}\right) \geqslant \frac{16}{3} \pi^2 \alpha_\ve \right\} =\left\{
x_{1},\cdots ,x_{N}\right\} , 
\end{equation*}%
then it follows from Proposition \ref{prop8.1} that%
\begin{equation*}
\nu =\sum_{i=1}^{N}\nu _{i}\delta _{x_{i}}, 
\end{equation*}%
here $\nu _{i}\geqslant 0$ and $\sum_{i=1}^{N}\nu _{i}=1$. On the other hand, we
have%
\begin{equation*}
\int_{\mathbb{B}^4}p \ud\nu =0
\end{equation*}%
for all $p \in \mathring{\mathcal{P}}_{m}$. In other words, we conclude that%
\begin{align*}
\frac{16}{3}\pi^2 \alpha_\ve N \leqslant &1 \qquad \mathrm{~~and} \\
\sum_{i=1}^{N}\nu _{i} p \left( x_{i}\right) =&0 \qquad \mathrm{~~for~~all~~} p \in \mathring{\mathcal{P}}_{m}.
\end{align*}%
This indicates that $N\in \mathcal{N}_{m}(\Sp^3)$ and thus $N\geqslant
N_{m} $. Moreover,%
\begin{equation*}
\alpha_\ve \leqslant \frac{3}{16\pi^2 N}\leqslant \frac{3}{16\pi^2 N_{m}(\Sp^3)}.
\end{equation*}%
However, this contradicts the choice of $\alpha_\ve$.
\end{proof} 

\subsection{Almost sharpness}\label{Sect:n=4}

Provided that we adopt Chang-Hang type test function in \cite{Chang-Hang} as the restriction of our test function to the boundary, we shall face a dilemma: How to extend it to the interior of the unit ball, since Chang-Hang's test function  is only piecewise Lipschitz? To demonstrate our idea,  as a good warm-up we construct an example to show the sharpness of Widom inequality in Theorem \ref{Thm:Widom}. For the completeness, it is left to Appendix \ref{Append}. Thanks to this stimulating example together `\emph{nice}' properties of ``\emph{geometric local bubble}", it extricates us from the above dilemma. 

\vskip 8pt
\subsubsection{An improvement of the Chang-Hang type estimate}\label{Subsect:improved_estimate}

With the help of the lower bounds \eqref{lbd:DGS} of $N_m(\Sp^{n-1})$, we give an elementary proof of the exact value of $N_3(\Sp^{n-1})$, which is originally due to Mysovskih \cite{Mysovskih}.
\begin{lemma}\label{lem:N_3}
For $n \geqslant 2$, there holds $N_3(\Sp^{n-1})=2n$.
\end{lemma}
\begin{proof}
It suffices to consider $n\geqslant 3$.  On one hand, we have $N_3(\Sp^{n-1})\geqslant 2n$ by virtue of \eqref{lbd:DGS}. On the other hand, a basis of $\mathring{\mathcal{P}}_3$ can be chosen as 
\begin{align*}
\bigg\{& x_i,\quad 1 \leqslant i \leqslant n; ~ x_i^2-\frac{|x|^2}{n}, \quad 1\leqslant i\neq j \leqslant n-1;~x_i x_j, \quad 1\leqslant i\neq j \leqslant n;\\
&~x_i(|x|^2 -(n+2)x_j^2),\quad \quad 1 \leqslant i\neq j\leqslant n;~x_ix_jx_k, \quad 1 \leqslant i\neq j\neq k \leqslant n~\bigg\}.
\end{align*}
We can choose $\nu_i=1/(2n)$ for $1 \leqslant i \leqslant n$ and vertices as $\{\pm e_i; 1 \leqslant i \leqslant n\}$. Then it is straightforward to check that for any element $p_j$ belonging to the above basis, 
$$\sum_{i=1}^n \nu_i\left( p_j(e_i)+p_j(-e_i)\right)=0.$$

Putting these facts together, the desired assertion follows.
\end{proof}
A key observation on $\mathring{\mathcal{P}}_1$ is very important to our improved estimate.
\begin{proposition}\label{prop:3-dim_m=1}
Let
$$\phi_{\ve,1}(x)=-\log (\ve^2+\mathrm{dist}_{\Sp^3}(x,N)^2)\quad \mathrm{~~and~~} \quad \phi_{\ve,2}(x)=-\log (\ve^2+\mathrm{dist}_{\Sp^3}(x,S)^2),$$
where $N$ and $S$ are the north and south poles on $\Sp^3$, respectively, then for small $\delta>0$,
$$\int_{B_\delta(N)}e^{3\phi_{\ve,1}(x)}x_i \ud \mu_{\Sp^3}+\int_{B_\delta(S)}e^{3\phi_{\ve,2}(x)} x_i \ud \mu_{\Sp^3}=0$$
for all $1 \leqslant i \leqslant 4$.
\end{proposition}
\begin{proof}
For brevity, we set $\rho=\mathrm{dist}_{\Sp^3}(x,N)$ and $x=(\sin \rho~ \xi ,\cos \rho), \xi \in \Sp^2$. Then we distinguish our discussion into two cases:
\begin{itemize}
\item For $1 \leqslant i \leqslant 3$, it follows from symmetry that
\begin{align*}
&\int_{B_\delta(N)}e^{3\phi_{\ve,1}(x)}x_i \ud \mu_{\Sp^3}+\int_{B_\delta(S)}e^{3\phi_{\ve,2}(x)} x_i \ud \mu_{\Sp^3}\\
=&\left[\int_0^\delta \frac{1}{(\ve^2+\rho^2)^3}  \sin^3 \rho\ud \rho+\int_{\pi-\delta}^\pi \frac{1}{(\ve^2+(\pi-\rho)^2)^3} \sin^3 \rho\ud \rho \right]\cdot \left(\int_{\Sp^2}\xi_i \ud \mu_{\Sp^2}\right)\\
=&0.
\end{align*}
\item For $i=4$,
\begin{align*}
&\int_{B_\delta(N)}e^{3\phi_{\ve,1}(x)}x_i \ud \mu_{\Sp^3}+\int_{B_\delta(S)}e^{3\phi_{\ve,2}(x)} x_i \ud \mu_{\Sp^3}\\
=&|\Sp^2|\left[\int_0^\delta \frac{1}{(\ve^2+\rho^2)^3} \cos \rho \sin^2 \rho\ud \rho+\int_{\pi-\delta}^\pi  \frac{1}{(\ve^2+(\pi-\rho)^2)^3} \cos \rho \sin^2 \rho\ud \rho \right]\\
=&0.
\end{align*}
\end{itemize}

Putting these facts together, the desired assertion follows.
\end{proof}

Though $n=4$ is enough to our later use, Proposition \ref{prop:3-dim_m=1} above motives us to prove a generic result.

\begin{proposition}\label{prop:n-dim_N_2}
 For every $p\in \mathring{\mathcal P}_m$ with $m\geqslant 1$, there exist $\nu_1,\cdots,\nu_{n+1} \in (0,1)$ and  pairwise distinct points $\{\bar x_k; 1\leqslant k \leqslant n+1\}\subset \Sp^{n-1}$ with $n\geqslant 3$ such that 
$$\sum_{k=1}^{n+1}\nu_k=1 \qquad \mathrm{~~and~~} \qquad \sum_{k=1}^{n+1}\nu_k p(\bar x_k)=0.$$
Choose $\delta>0$ sufficiently small such that geodesic balls $\{B_\delta(\bar x_k); 1 \leqslant k\leqslant n+1\} \subset \Sp^{n-1}$ satisfy
$$B_\delta(\bar x_i)\cap B_\delta(\bar x_j)=\emptyset, \qquad \forall~ 1 \leqslant i\neq j \leqslant n+1.$$
Then
$$\sum_{k=1}^{n+1}\nu_k\int_{B_\delta(\bar x_k)}e^{(n-1)\phi_{\ve,k}(x)}p(x)\ud \mu_{\Sp^{n-1}}=\sum_{k=1}^{n+1}\nu_k\int_{B_\delta(\bar x_k)}e^{(n-1)\phi_{\ve,k}}O(\rho^3)\ud \mu_{\Sp^{n-1}},$$
where 
$$\phi_{\ve,k}(x)=-\log\bigg(\ve^2+\mathrm{dist}_{\Sp^{n-1}}(x,\bar x_k)^2\bigg), \qquad 1 \leqslant k \leqslant n+1.$$
\end{proposition}
\begin{proof}

Since 
$$\mathring{\mathcal P}_m=\bigoplus_{l=1}^m \mathcal{H}_l \qquad \mathrm{on~~} \Sp^{n-1}$$
by virtue of \cite[Theorem 2.1]{Stein-Weiss},
where $\mathcal{H}_l$ is the set of all homogeneous harmonic polynomials of degree $l$ in $\Rn$, it suffices to consider each $\mathcal{H}_l$ instead of $\mathring{\mathcal P}_m$. Without loss of generality, we assume $p \in \mathcal{H}_{l_0}$ for some $1 \leqslant l_0 \leqslant m$.

Fix each $\bar x_k$, rotate $\bar x_k$ properly to the north pole $N$. In other words, we can find some  $Q \in SO(n)$ such that $Q \bar x_k=N$. Under the change of variables: $y=Qx$, we let $\tilde p(y)=p(Q^\top y)$ and now choose local coordinates around $N$ such that
$$y=(\sin \rho~ \xi, \cos \rho)\in \Sp^{n-1}, \quad \xi \in \Sp^{n-2}.$$ 
Notice that
\begin{align*}
\frac{\pa \tilde p}{\pa \rho}=&\sum_{i=1}^{n-1}\frac{\pa  \tilde p}{\pa y_i} \frac{\pa y_i}{\pa \rho}+\frac{\pa  \tilde  p}{\pa y_n} \frac{\pa y_n}{\pa \rho},\\
\frac{\pa^2 \tilde p}{\pa \rho^2}=&\sum_{i,j=1}^{n-1}\frac{\pa^2 \tilde p}{\pa y_i \pa y_j} \frac{\pa y_i}{\pa \rho} \frac{\pa y_j}{\pa \rho}+\sum_{i=1}^{n-1}\frac{\pa \tilde p}{\pa y_i} \frac{\pa^2 y_i}{\pa \rho^2}+2\sum_{i=1}^{n-1}\frac{\pa^2 \tilde p}{\pa y_i \pa y_n} \frac{\pa y_i}{\pa \rho}\frac{\pa y_n}{\pa \rho}+\frac{\pa \tilde p}{\pa y_n} \frac{\pa^2 y_n}{\pa \rho^2}.
\end{align*}
In particular, at $\rho=0$ we have
\begin{align*}
\frac{\pa \tilde p}{\pa \rho}(N)=&\sum_{i=1}^{n-1}\frac{\pa \tilde p}{\pa y_i}(N) \xi_i,\\
\frac{\pa^2 \tilde p}{\pa \rho^2}(N)=&\sum_{i,j=1}^{n-1}\frac{\pa^2 \tilde p}{\pa y_i \pa y_j}(N) \xi_i\xi_j-\left(y_n\frac{\pa  \tilde  p}{\pa y_n}\right)(N).
\end{align*}

For brevity, we define
$$a_{ij}=\frac{1}{2}\frac{\pa^2 \tilde p}{\pa y_i \pa y_j}(N), \qquad b_i=\frac{\pa \tilde p}{\pa y_i}(N), \qquad c_n=\frac{1}{2}\left(y_n\frac{\pa \tilde p}{\pa y_n}\right)(N)$$
for $1\leqslant i,j \leqslant n-1$. Then the Taylor's expansion of $\tilde p(y)$ in a neighborhood around $N$ is
$$\tilde p(y)=\tilde p(N)+\left(\sum_{i=1}^{n-1}b_i \xi_i\right)\rho+\left(\sum_{i,j=1}^{n-1}a_{ij} \xi_i\xi_j-c_n\right)\rho^2+O(\rho^3).$$

Observe that
$$\int_{B_\delta(\bar x_k)}e^{(n-1)\phi_{\ve,k}}p(x)\ud \mu_{\Sp^{n-1}}=\int_{B_\delta(N)}e^{(n-1)\phi_{\ve,k}(Q^\top y)}\tilde p(y)\ud \mu_{\Sp^{n-1}}.$$
Thus by symmetry we obtain 
\begin{align*}
&\int_{B_\delta(N)}e^{(n-1)\phi_{\ve,k}(Q^\top y)}\tilde p(y)\ud \mu_{\Sp^{n-1}}+c_n\int_{B_\delta(N)}e^{(n-1)\phi_{\ve,k}(Q^\top y)}\rho^2 \ud \mu_{\Sp^{n-1}}\\
=&\tilde p(N) |\Sp^{n-2}|\int_0^\delta  \frac{1}{(\ve^2+\rho^2)^{n-1}} \sin^{n-2} \rho \ud \rho\\
&+\sum_{i=1}^{n-1}b_i \int_0^\delta \frac{\rho}{(\ve^2+\rho^2)^{n-1}} \sin^{n-2} \rho \ud \rho \int_{\Sp^{n-2}} \xi_i\ud \mu_{\Sp^{n-2}}\\
&+\sum_{i,j=1}^{n-1} \left(a_{ij}\int_{\Sp^{n-2}} \xi_i\xi_j \ud \mu_{\Sp^{n-2}}\right)\int_0^\delta \frac{\rho^2}{(\ve^2+\rho^2)^{n-1}} \sin^{n-2} \rho \ud \rho \\
&+\int_{B_\delta(N)}e^{(n-1)\phi_{\ve,k}(Q^\top y)}O(\rho^3) \ud \mu_{\Sp^{n-1}}\\
=&\tilde p(N) |\Sp^{n-2}|\int_0^\delta  \frac{1}{(\ve^2+\rho^2)^{n-1}} \sin^{n-2} \rho \ud \rho\\
&+ \frac{ |\Sp^{n-2}|}{n-1}\left(\sum_{i=1}^{n-1}a_{ii}\right)\int_0^\delta \frac{\rho^2}{(\ve^2+\rho^2)^{n-1}} \sin^{n-2} \rho \ud \rho\\
&+\int_{B_\delta(N)}e^{(n-1)\phi_{\ve,k}(Q^\top y)}O(\rho^3) \ud \mu_{\Sp^{n-1}}.
\end{align*}

Next,  we notice that $p(x) \in \mathcal{H}_{l_0}$ implies $\tilde p(y)\in \mathcal{H}_{l_0}$.
Observe that
\begin{align*}
\sum_{i=1}^{n-1}a_{ii}=\sum_{i=1}^{n-1}\frac{\pa^2 \tilde p}{\pa y_i^2}=\Delta_{\Sp^{n-1}}\tilde p=-l_0(n+l_0-2)\tilde p \quad \mathrm{~~at~~}N.
\end{align*}
Moreover, we claim that
\begin{align*}
2c_n=\left[y_n\pa_{y_n}\tilde p(y)\right]\big|_{y=N}=l_0 p(\bar x_k).
\end{align*}

To that end, denote by $Q=(q_{ij})_{n \times n}$ and then $x_i=q_{ij} y_j$ by definition, hereafter we shall use Einstein summation notation. Notice that
$$Q \bar x_k=N \quad  \Longleftrightarrow \quad \bar x_k=(q_{1n},\cdots,q_{nn})^\top.$$
Using
$$\pa_{y_n}\tilde p= \frac{\pa x_l}{\pa y_n}\pa_{x_l} p=q_{l n}\pa_{x_l} p,$$
we have
\begin{align*}
\left(y_n\pa_{y_n}\tilde p\right)\big|_{y=N}=&\left(q_{ln}x_l q_{sn}\pa_{x_s} p\right)\big|_{x=\bar x_k}\\
=&(q_{ln}q_{ln})(\bar x_k)_s\pa_{x_s} p(\bar x_k)=\left(x_s\pa_{x_s} p\right)\big|_{x=\bar x_k}=l_0 p(\bar x_k).
\end{align*}
Hence, the desired claim follows.

\vskip 8pt

Going back to the original variable $x$, we combine all facts together to conclude that
\begin{align*}
&\int_{B_\delta(\bar x_k)}e^{(n-1)\phi_{\ve,k}(x)}p(x)\ud \mu_{\Sp^{n-1}}\\
=&\mathscr{C}_{\ve,\delta,n} p(\bar x_k)+\int_{B_\delta(\bar x_k)}e^{(n-1)\phi_{\ve,k}}O(\rho^3)\ud \mu_{\Sp^{n-1}},
\end{align*}
where
\begin{align*}
\mathscr{C}_{\ve,\delta,n}=&|\Sp^{n-2}|\int_0^\delta  \frac{1}{(\ve^2+\rho^2)^{n-1}} \sin^{n-2} \rho \ud \rho\\
&-\frac{ |\Sp^{n-2}|}{n-1}l_0(n+l_0-2)\int_0^\delta \frac{\rho^2}{(\ve^2+\rho^2)^{n-1}} \sin^{n-2} \rho \ud \rho\\
&-\frac{l_0}{2} |\Sp^{n-2}|\int_0^\delta  \frac{\rho^2}{(\ve^2+\rho^2)^{n-1}} \sin^{n-2} \rho \ud \rho
\end{align*}
is a constant independent of $k$. Consequently, we arrive at
\begin{align*}
&\sum_{k=1}^{n+1}\nu_k\int_{B_\delta(\bar x_k)}e^{(n-1)\phi_{\ve,k}(x)}p(x)\ud \mu_{\Sp^{n-1}}\\
=&\mathscr{C}_{\ve,\delta,n}\sum_{k=1}^{n+1}\nu_k p(\bar x_k)+\sum_{k=1}^{n+1}\nu_k\int_{B_\delta(\bar x_k)}e^{(n-1)\phi_{\ve,k}}O(\rho^3) \ud \mu_{\Sp^{n-1}}\\
=&\sum_{k=1}^{n+1}\nu_k\int_{B_\delta(\bar x_k)}e^{(n-1)\phi_{\ve,k}}O(\rho^3) \ud \mu_{\Sp^{n-1}}.
\end{align*}
This finishes the proof.
\end{proof}

\subsubsection{Construction of test functions}\label{Subsect:n=4}

With preparations above, we are now in a position to construct our example in dimension four.
\vskip 8pt
\begin{proof}[Proof of Proposition \protect{\ref{example:n=4}}]
For each $m$, we can find a  $\nu=\sum_{i=1}^{N_m(\Sp^3)} \nu_i \delta_{x_i} \in \mathcal{M}_{m}^{c}\left(\mathbb{S}^3\right)$. Due to the precise estimate of $m=1$  in  Proposition \ref{prop:3-dim_m=1}, we shall point out modifications if necessary.

	 \begin{figure}[h]
\centering
\includegraphics[width=0.5\textwidth]{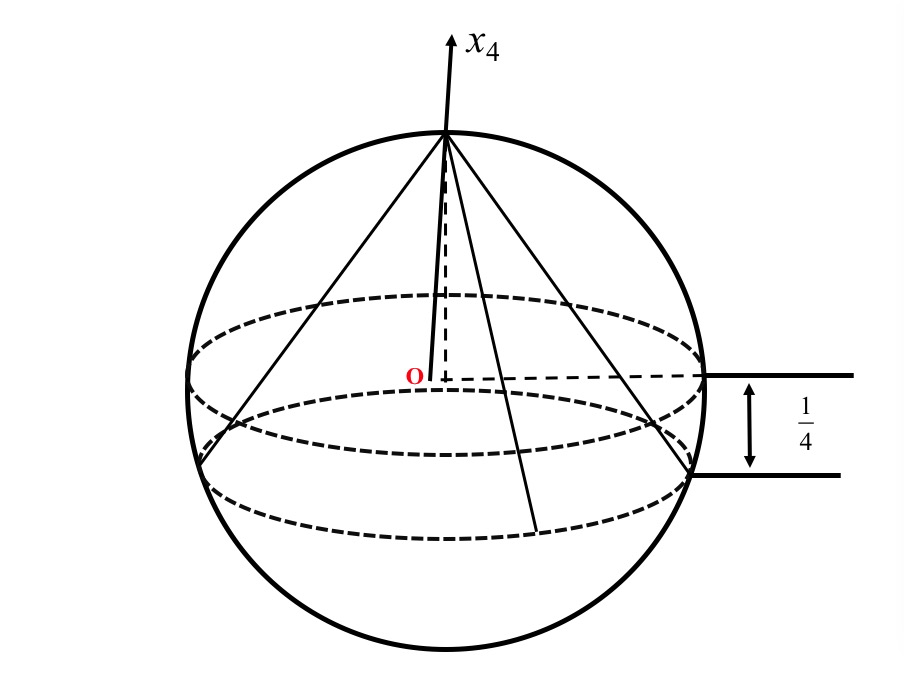}
\vskip -4pt
\caption{A regular $4$-simplex $\mathscr{R}$}\label{fig:4-simplex}
\end{figure}

\begin{figure}
    \centering
    \subfigure[]{\includegraphics[width=0.4\textwidth]{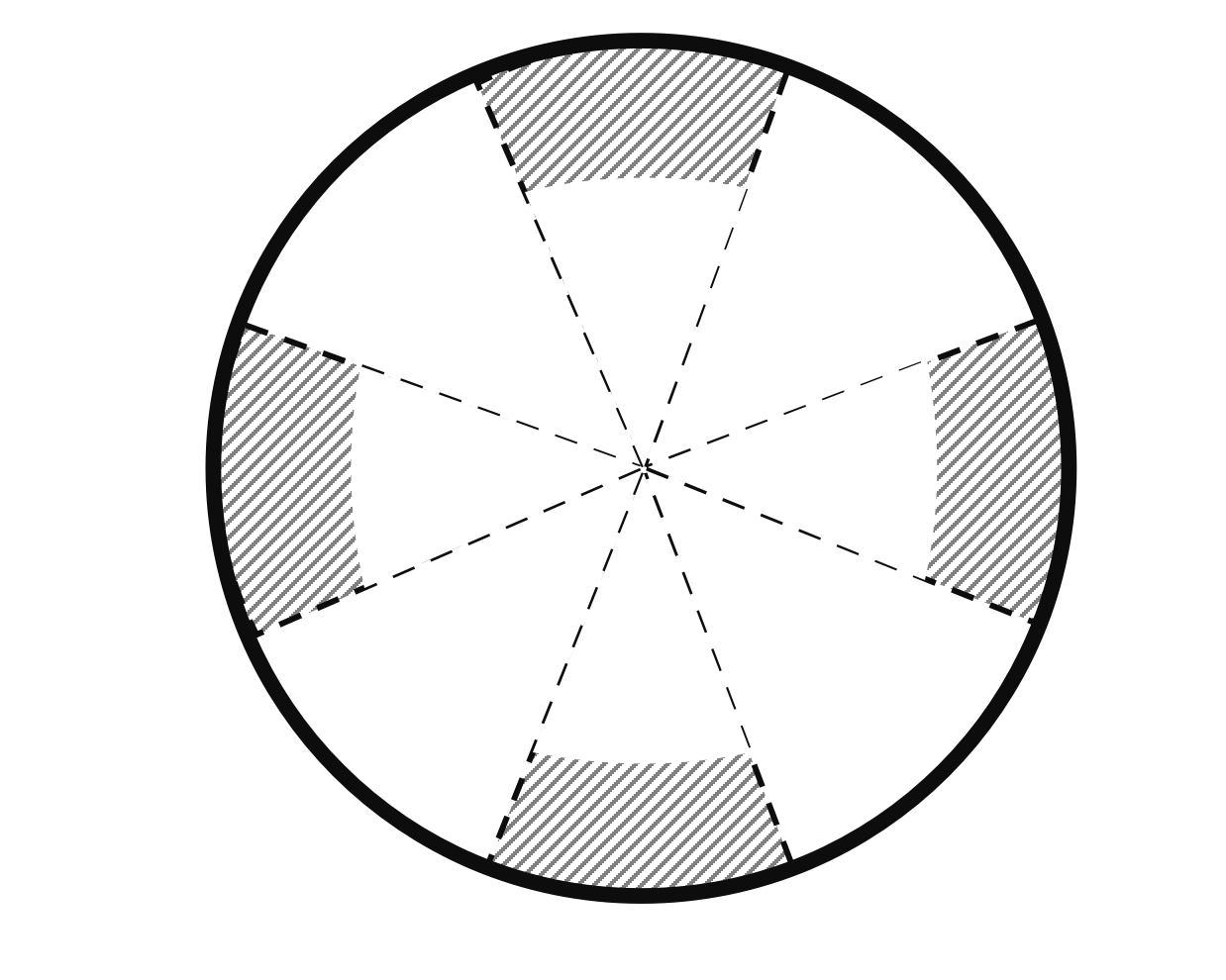}} 
    \subfigure[]{\includegraphics[width=0.4\textwidth]{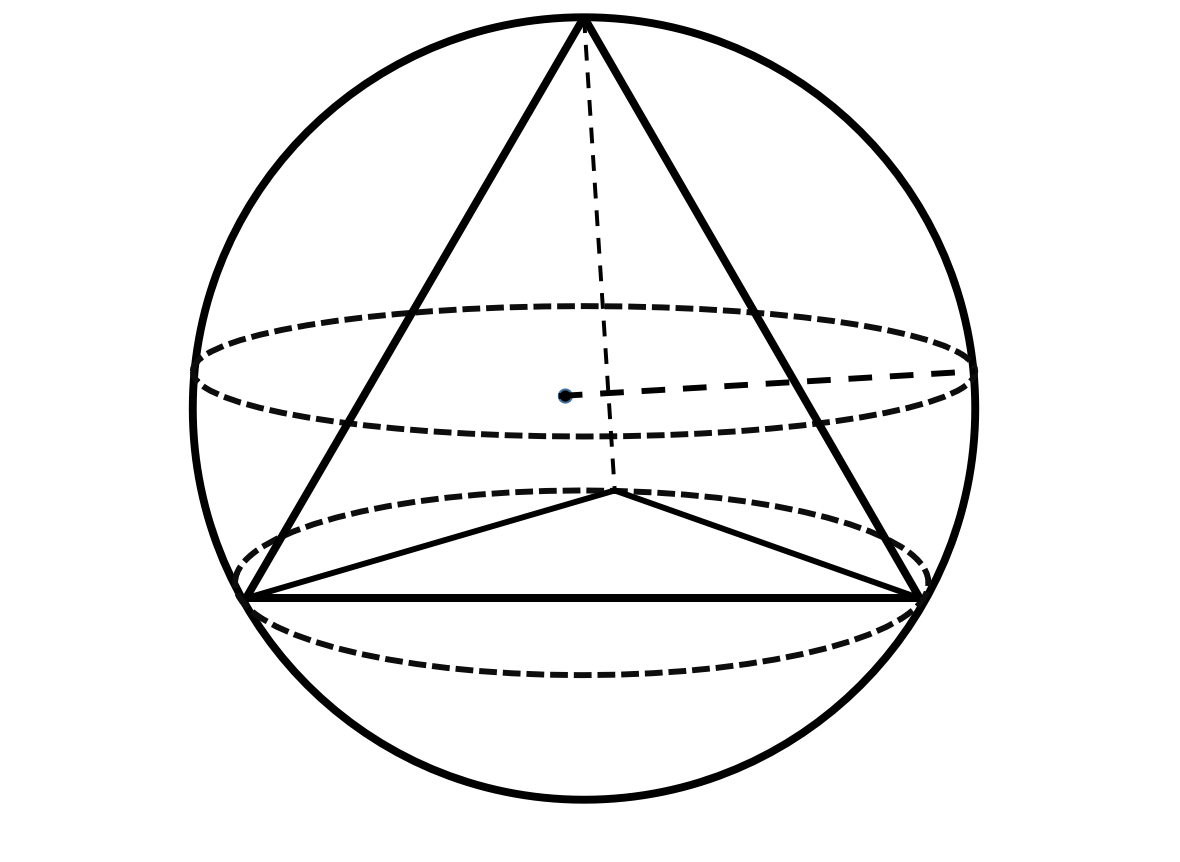}} 
    \vskip -4pt
        \caption{On  $\left\{x_4=-\frac{1}{4}\right\}\cap \mathscr{R}$: (a) Conic Annuli  \qquad (b) A regular $3$-simplex}
    \label{fig:Conic Annuli}
\end{figure}
To illustrate the number $N_m(\Sp^3)$, we take $m=2$ for example. It follows from \cite[Lemma 3.4]{Hang} that $N_2(\Sp^3)=5$ can be attained by some $\nu=\sum_{i=1}^5 \nu_i \delta_{x_i} \in \mathcal{M}_{2}^{c}\left(\mathbb{S}^3\right)$, where $\nu_i=1/5$  and distinct points $x_i\in \Sp^3$ for $1\leqslant i \leqslant 5$ are vertices of a regular $4$-simplex embedded in $\mathbb{B}^4$. One possible choice of vertices could be
\begin{align*}
x_1=(0,0,0,1), \qquad&\\
x_2=(0,0,\frac{\sqrt{15}}{4},-\frac{1}{4}), \qquad&
x_3=(0,\frac{\sqrt{30}}{6},-\frac{\sqrt{15}}{12},-\frac{1}{4}), \\
x_4=(\frac{\sqrt{10}}{4},-\frac{\sqrt{30}}{12},-\frac{\sqrt{15}}{12},-\frac{1}{4}),\qquad&
x_5=(-\frac{\sqrt{10}}{4},-\frac{\sqrt{30}}{12},-\frac{\sqrt{15}}{12},-\frac{1}{4}). 
\end{align*}

In the following, we set $N_m=N_m(\Sp^3)$ for brevity.

For  each $1 \leqslant i \leqslant N_m$,  choose conic annuli $\mathcal{A}_{\delta}(x_i)$ and local coordinates around $x_i$ as before, and let $\chi_i (x)$ be a smooth cut-off function such that $\chi_i(x)=1$ in $\mathcal{A}_{\delta}(x_i)$ and $\chi_i(x)=0$ outside $\mathcal{A}_{2\delta}(x_i)$. Under the above coordinates, for any $0 < \ve <\delta$ we define
$$\phi_\ve(r,\rho)=-\log\left((\ve+1-r)^2+\rho^2\right)+\frac{2\ve(1-r)}{(\ve+1-r)^{2}+\rho^{2}}$$
and
	\begin{align*}
		\phi_{\ve, i}(r,\rho)=\phi_\ve(r,\rho)+\frac{1}{3}\log \nu_i,
		\end{align*}
		where $x=r\xi$ with $r=|x|$ and $\rho=\stackrel{\frown}{\xi x_i}$.
	
We estimate
\begin{align}\label{est2:n=4}
		&\sum_{i=1}^{N_m}\int_{\Sp^3}\chi_i(x)e^{3\phi_{\ve, i}(x)}\ud \mu_{\Sp^3}=\sum_{i=1}^{N_m}\int_{B_{2\delta}(x_i)}\chi_ie^{3\phi_{\ve, i}}\ud \mu_{\Sp^3}\no\\
		=&\sum_{i=1}^{N_m}\int_{B_{2\delta}\left(x_{i}\right)}\frac{\chi_i(x) \nu_{i}}{(\varepsilon^{2}+\rho^{2})^{3}}\ud \mu_{\Sp^3}\no\\
	=&|\Sp^2| \int_0^{\delta}(\ve^2+\rho^2)^{-3}\sin^2\rho \ud \rho+O(\ve^{-3})\int_{\frac{\delta}{\ve}}^{\frac{2\delta}{\ve}}(1+t^2)^{-3}t^{2}\ud t\no\\
	=&\frac{\pi^2}{4}\ve^{-3}+O(\frac{1}{\ve}).
		\end{align}
and then
\begin{align*}
\int_{\mathbb{S}^{3}}e^{3u}\ud \mu_{\Sp^3} & =\sum_{i=1}^{N_m}\int_{B_{2\delta}(x_i)}\chi_ie^{3\phi_{\ve, i}}\ud \mu_{\Sp^3}+\int_{\mathbb{S}^{3}}\bigg(c_1\log\frac{1}{\varepsilon}+\sum_{j=1}^{L}\beta_{j}\psi_{j}\bigg) \ud \mu_{\Sp^3}\no\\
 & =\frac{\pi^2}{4}\ve^{-3}+O(\log\frac{1}{\varepsilon}).
\end{align*}
This implies 
\begin{equation}\label{area-ln}
\log\left(\frac{1}{2\pi^2} \int_{\mathbb{S}^{3}}e^{3u}\ud \mu_{\Sp^3} \right)=3\log\frac{1}{\varepsilon}+O(\log\log\frac{1}{\varepsilon}).
\end{equation}
		
		For every $p \in \mathring{\mathcal{P}}_{m}$, it follows from Proposition \ref{prop:n-dim_N_2} that \footnote{If $m=1$, then $N_1(\Sp^3)=2$ can be achieved by a Dirac probability measure $\nu=(\delta_N+\delta_S)/2$, and  by Proposition \ref{prop:3-dim_m=1} we have $\int_{\Sp^3}\left(\sum_{i=1}^2 \chi_i e^{3\phi_{\ve,i}}\right)p\ud \mu_{\Sp^3}=O(1)$ for every $p\in\mathring{\mathcal{P}}_1$, where $\phi_{\ve,1}(x)=-\log(\ve^2+\mathrm{dist}_{\Sp^3}(x,N)^2)+\frac{1}{3}\ln \frac{1}{2}$ and $\phi_{\ve,2}(x)=-\log(\ve^2+\mathrm{dist}_{\Sp^3}(x,S)^2)+\frac{1}{3}\ln \frac{1}{2}$.}
\begin{align}\label{est1:4th_n=4}
 \int_{\Sp^3} \left(\sum_{i=1}^{N_m}\chi_i(x)e^{3\phi_{\ve, i}(x)}\right)p \ud \mu_{\Sp^3} =& \sum_{i=1}^{N_m} \nu_i \int_{B_\delta(x_i)}e^{3\phi_\ve}O(\rho^3) \ud\mu_{\Sp^3}\no\\
 =& |\Sp^2| \sum_{i=1}^{N_m} \nu_i \int_0^\delta \frac{O(\rho^3)}{(\ve^2+\rho^2)^3} \sin^2 \rho \ud\rho+O(1)\no\\
 =&O(\log\frac{1}{\varepsilon}).
      \end{align}

Similar to the proof of Proposition \ref{example:n>4}, we can find a basis of spherical harmonics
$\{p_{1}, \cdots, p_{L}\}$ of $\mathring{\mathcal{P}}_{m}$ with $L=m(2m^2+9m+13)/6$,
 and functions $\psi_{1}, \cdots, \psi_{L} \in C_{c}^{\infty}\left(\mathbb{B}^{4} \backslash \bigcup_{i=1}^{N_m} \overline{\mathcal{A}_{2 \delta}\left(x_{i}\right)}\right)$ such that the determinant 
\begin{equation}\label{nonsingular:Gram_n=4}
	\det\left[\int_{\mathbb{S}^{3}} \psi_{j} p_{k} \ud \mu_{\Sp^{3}}\right]_{1 \leqslant j, k \leqslant L} \neq 0.
\end{equation}
To this end, we can choose a nonzero smooth function $\eta \in C_{c}^{\infty}\left(\overline{\mathbb{B}^{4} }\backslash \bigcup_{i=1}^{N_m} \mathcal{A}_{2 \delta}\left(x_{i}\right)\right)$ such that $\eta p_{1}, \cdots, \eta p_{L}$ are linearly independent. Actually, one possible way of $\eta$ is to choose $\eta(x)=\eta_1(r)\hat \chi(\rho,\theta)$, where $\hat \chi(\rho,\theta)$ is a smooth cut-off function supported in $\Sp^3\setminus \bigcup_{i=1}^{N_m} B_{2\delta}(x_i)$ and $\eta_1(r)$ is the smooth cut-off function used in Section \ref{Sect:n>4}.  Then  it follows that the Gram matrix
\begin{align*}
	\left[\int_{\mathbb{S}^{3}} \eta^{2} p_{j} p_{k} \ud \mu_{\Sp^{3}}\right]_{1 \leqslant j, k \leqslant L}
\end{align*}
is positive definite, then $\psi_{j}=\eta^{2}p_j$ satisfies \eqref{nonsingular:Gram_n=4} and $\pa_r \psi_j=0$ on $\Sp^3$.

The fact \eqref{nonsingular:Gram_n=4} enables us to find $\beta_{1}, \cdots, \beta_{L} \in \mathbb{R}$ such that 
\begin{align}\label{vanishing_moments:n=4_1}
	\int_{\mathbb{S}^3}\left(\sum_{i=1}^{N_m}\chi_i(x)e^{3\phi_{\ve, i}(x)}+\sum_{j=1}^{L} \beta_{j} \psi_{j}\right) p_{k} \ud \mu_{\Sp^{3}}=0 \qquad \forall~ 1 \leqslant k \leqslant L.
\end{align}
Moreover, it follows from \eqref{est1:4th_n=4} that for all $1 \leqslant j \leqslant L$, $\beta_{j}=O\left(\log\frac{1}{\varepsilon}\right)$. As a consequence we can find a constant $c_1 >0$ such that 
\begin{align*}
	\sum_{j=1}^{L} \beta_{j} \psi_{j}+c_{1} \log\frac{1}{\varepsilon} \geqslant \log\frac{1}{\varepsilon}.
\end{align*}

In dimension four, we may choose our test function as
\begin{equation}\label{test fcn:n=4}
e^{3u}=\sum_{i=1}^{N_m}\chi_ie^{3\phi_{\ve, i}}+\sum_{j=1}^{L} \beta_{j} \psi_{j}+c_{1}\log\frac{1}{\varepsilon}.
\end{equation}
Here we would like to take $\chi_{i}$ to be a variable separation cut-off
function in order to ensure zero Neumann boundary condition without additional correction. To that end, as in Section \ref{Sect:n>4} we may choose $\chi_i(x)=\eta_{1}(r) \widetilde{\chi}_{i}(\rho,\theta)$ such that $\pa_r \chi_i=0$ on $\Sp^3$. Since $\pa_r \phi_{\ve,i}=0$ on $\Sp^3$ by virtue of \eqref{bubble:Neumann_bdry_cond},  it follows follows from \eqref{test fcn:n=4} that 
\begin{align*}
&3 e^{3u} \pa_r u =\sum_{i=1}^{N_m}\pa_r \chi_ie^{3\phi_{\ve, i}}+3\sum_{i=1}^{N_m}\chi_i e^{3\phi_{\ve, i}}\pa_r \phi_{\ve,i}+\sum_{j=1}^{L} \beta_{j} \pa_r \psi_{j}=0 \qquad \mathrm{on~~}\Sp^3\\
\Longrightarrow&~~ \pa_r u=0 \qquad  \mathrm{on~~}\Sp^3.
\end{align*}

\vskip 8pt

By \eqref{vanishing_moments:n=4_1}  and \eqref{test fcn:n=4} we have 
$$\int_{\Sp^3}e^{3u}p_k \ud\mu_{\Sp^3}=0 \qquad \mathrm{for~~} 1 \leqslant k \leqslant L. $$

On one hand, we have
\begin{align*}
&\bar{u}	|\mathbb{S}^{3}|\\
=&\int_{\mathbb{S}^{3}}u \ud\mu_{\Sp^3}=\int_{\mathbb{S}^{3}}\frac{1}{3}\log(e^{3u})\ud\mu_{\Sp^3}\\
	=&\frac{1}{3}\int_{\mathbb{S}^{3}}\log\bigg(\sum_{i=1}^{N_m}\chi_{i}e^{3\phi_{\varepsilon,i}}+c_1\log\frac{1}{\varepsilon}+\sum_{j=1}^{L}\beta_{j}\psi_{j}\bigg)\ud\mu_{\Sp^3}\\
	\leqslant&\frac{1}{3} \sum_{i=1}^{N_m}\int_{B_{\delta}(x_{i})}\left[\log\left(2 e^{3\phi_{\varepsilon,i}}\right)+\log\left(2c_1\log\frac{1}{\varepsilon}\right)\right]\ud\mu_{\Sp^3}\\
	&+\frac{1}{3} \sum_{i=1}^{N_m}\int_{\mathbb{S}^{3}\backslash B_{\delta}(x_{i})}\log \left(O(\log\frac{1}{\varepsilon})\right)\ud\mu_{\Sp^3}\\
	\lesssim&\frac{1}{3}\bigg(\varepsilon^{3}\int_{0}^{\delta}(\frac{\rho}{\varepsilon})^{2} \log\frac{2\nu_{i}}{\varepsilon^{2}(1+(\frac{\rho}{\varepsilon})^{2})}\ud(\frac{\rho}{\varepsilon})\bigg)+O(\log\log\frac{1}{\varepsilon})\\
	\lesssim&\frac{1}{3}\varepsilon^{3}\int_{0}^{\delta/\varepsilon}\bigg(t^{2}\log\frac{1}{\varepsilon^{2}}+t^{2}\log\frac{1}{1+t^{2}}\bigg)\ud t+O(\log\log\frac{1}{\varepsilon})\\
	\lesssim&\frac{1}{3}\varepsilon^{3}\bigg[\int_{0}^{\delta/\varepsilon}t^{2}\log\frac{1}{\varepsilon^{2}}\ud t-\left(\frac{2}{3}\frac{\delta}{\varepsilon}-\frac{2}{9}(\frac{\delta}{\varepsilon})^{3}+\frac{1}{3}(\frac{\delta}{\varepsilon})^{3}\log[1+(\frac{\delta}{\varepsilon})^{2}]\right)\bigg]+O(\log\log\frac{1}{\varepsilon})\\
	\lesssim &O(\log\log\frac{1}{\varepsilon}),
\end{align*}
where 
\begin{align*} 
& \int_{0}^{\delta/\varepsilon}t^{2}\log\left(1+t^{2}\right)\ud t\\
=&\frac{1}{9}\left[6t-2t^{3}-6\arctan t+3t^{3}\log(1+t^{2})\right]\big|_{t=\frac{\delta}{\varepsilon}}\\
= & \frac{2}{3}\frac{\delta}{\varepsilon}-\frac{2}{9}(\frac{\delta}{\varepsilon})^{3}+\frac{1}{3}(\frac{\delta}{\varepsilon})^{3}\log[1+(\frac{\delta}{\varepsilon})^{2}]+O(1).
\end{align*}
On the other hand,
\begin{align*}
\bar{u}|\mathbb{S}^{3}|	=&\int_{\mathbb{S}^{3}}\frac{1}{3}\log\bigg(\sum_{i=1}^{N_m}\chi_{i}e^{3\phi_{\varepsilon,i}}+c_1\log\frac{1}{\varepsilon}+\sum_{j=1}^{L}\beta_{j}\psi_{j}\bigg)\ud\mu_{\Sp^3}\\
	\geqslant&\int_{\mathbb{S}^{3}}\frac{1}{3}\log\left(c_1\log\frac{1}{\varepsilon}\right)\ud\mu_{\Sp^3}\\
	\geqslant& O(\log\log\frac{1}{\varepsilon}).
\end{align*}
Thus, we put these together to show
\begin{equation}\label{average-u}
\bar{u}=O(\log\log\frac{1}{\varepsilon}).
\end{equation}

Under local coordinates of the flat metric  near each $x_i$
$$|\ud x|^2=\ud r^{2}+r^{2}(\ud\rho^{2}+\sin^{2}\rho g_{\mathbb{S}^{2}}),$$ 
there holds 
\begin{align*}
\Delta u	=&\partial_r^{2} u+\frac{3}{r}\partial_r u+\frac{1}{r^{2}}\left(\partial_\rho^{2} u+2\cot \rho \partial_\rho u+\sin^{-2}\rho\Delta_{\mathbb{S}^{2}}u\right)\\
	=&\partial_r^{2} u+\partial_\rho^{2} u+\frac{2}{\rho}\partial_\rho u+\left(\frac{2\cot\rho}{r^{2}}-\frac{2}{\rho}\right) \partial_\rho u+(\frac{1}{r^{2}}-1)\partial_\rho^{2} u+\frac{1}{r^{2}}\frac{1}{\sin^{2}\rho}\Delta_{\mathbb{S}^{2}}u.
	\end{align*}
	
In $\overline{\mathbb{B}^4}\backslash\mathscr{\mathcal{A}}_{\delta}(x_{i})$, by \eqref{test fcn:n=4} we have
$$u=\frac{1}{3}\log \bigg(\chi_i e^{3\phi_{\varepsilon,i}}+c_1\log\frac{1}{\varepsilon}\bigg) \qquad \mathrm{in~~} \mathcal{A}_{2\delta}\setminus \mathcal{A}_\delta(x_i)$$
and
$$u=\frac{1}{3}\log\bigg(c_1\log\frac{1}{\varepsilon}+\sum_{j=1}^{L} \beta_{j} \psi_{j}\bigg) \qquad \mathrm{in~~} \overline{\mathbb{B}^4}\setminus \mathcal{A}_{2\delta}(x_i).$$
Then it is not hard to verify that
 $$|\nabla u|_{\Sp^3}+|\Delta u|=O(1)\quad \overline{\mathbb{B}^4}\backslash\mathscr{\mathcal{A}}_{\delta}(x_{i})$$
	
In each $\mathscr{\mathcal{A}}_{\delta}(x_{i})$, again by \eqref{test fcn:n=4} we have 
\begin{equation*}
	u=\frac{1}{3}\log\left(e^{3\phi_{\varepsilon,i}}+c_1\log\frac{1}{\varepsilon}\right).
	\end{equation*}
Recall that $\phi_{\ve, i}=\phi_\ve+\frac{1}{3}\log \nu_i$. For convenience, we decompose $\phi_\ve=\phi_1+\phi_2$, where
$$\phi_1(r,\rho)=-\log\left((\ve+1-r)^2+\rho^2\right)$$ 
and 
$$\phi_2(r,\rho)=\frac{2\ve(1-r)}{(\ve+1-r)^{2}+\rho^{2}}.$$	

A direct computation yields
\begin{align*}
\partial_{r}\phi_{1}=&\frac{2(\varepsilon+1-r)}{(\varepsilon+1-r)^{2}+\rho^{2}},\\
\partial_{\rho}\phi_{1}=&\frac{-2\rho}{(\varepsilon+1-r)^{2}+\rho^{2}}
\end{align*}
and
\begin{align*}
\partial_r^{2}\phi_{1}=&-\frac{2}{(\varepsilon+1-r)^{2}+\rho^{2}}+\frac{4(\varepsilon+1-r)^{2}}{\left((\varepsilon+1-r)^{2}+\rho^{2}\right)^{2}},\\
\partial_\rho^{2}\phi_{1}=&-\frac{2}{(\varepsilon+1-r)^{2}+\rho^{2}}+\frac{4\rho^{2}}{\left((\varepsilon+1-r)^{2}+\rho^{2}\right)^{2}}.
\end{align*}
Also, 
\begin{align*}
\partial_{r}\phi_{2}=&\frac{-2\varepsilon}{(\varepsilon+1-r)^{2}+\rho^{2}}+\frac{4\varepsilon(1-r)(\varepsilon+1-r)}{((\varepsilon+1-r)^{2}+\rho^{2})^{2}},\\
\partial_{\rho}\phi_{2}=&-\frac{4\varepsilon(1-r)\rho}{((\varepsilon+1-r)^{2}+\rho^{2})^{2}}
\end{align*}
and
\begin{align*}
\partial_{r}^2\phi_{2}=&\frac{-8\varepsilon^{2}-12\varepsilon(1-r)}{((\varepsilon+1-r)^{2}+\rho^{2})^{2}}+\frac{16\varepsilon(1-r)(\varepsilon+1-r)^{2}}{((\varepsilon+1-r)^{2}+\rho^{2})^{3}},\\
\partial_{\rho}^2\phi_{2}=&\frac{-4\varepsilon(1-r)}{((\varepsilon+1-r)^{2}+\rho^{2})^{2}}+\frac{16\varepsilon(1-r)\rho^{2}}{((\varepsilon+1-r)^{2}+\rho^{2})^{3}}.
\end{align*}

Hence, it is straightforward to show 
\[
\partial_\rho u=\frac{e^{3\phi_{\varepsilon,i}} \pa_\rho \phi_{\varepsilon,i} }{e^{3\phi_{\varepsilon,i}}+c_1\log\frac{1}{\varepsilon}}=\frac{\pa_{\rho}\phi_\ve}{1+c_1\nu_{i}^{-1}\log\frac{1}{\varepsilon}\left((\varepsilon+1-r)^{2}+\rho^{2}\right)^{3}}
\]
and 
\begin{align*}
\partial_\rho^2 u =&\frac{\pa_{\rho}^2 \phi_\ve e^{3\phi_{\varepsilon,i}}}{e^{3\phi_{\varepsilon,i}}+c_1\log\frac{1}{\varepsilon}}+\frac{3c_1\nu_{i}^{-1}\log\frac{1}{\varepsilon}(\pa_{\rho}\phi_\ve)^{2}e^{3\phi_{\varepsilon,i}}}{\big(e^{3\phi_{\varepsilon,i}}+c_1\log\frac{1}{\varepsilon}\big)^{2}}\\
 =&\frac{\pa_{\rho}^2 \phi_\ve}{1+c_1\nu_{i}^{-1}\log\frac{1}{\varepsilon}\left((\varepsilon+1-r)^{2}+\rho^{2}\right)^{3}}+\frac{3c_1\nu_i^{-1}\log\frac{1}{\varepsilon}e^{-3\phi_{\varepsilon}}(\pa_{\rho}\phi_\ve)^{2}}{(1+c_1\nu_{i}^{-1}\log\frac{1}{\varepsilon}\left((\varepsilon+1-r)^{2}+\rho^{2}\right)^{3})^{2}},
\end{align*}
as well as similar formulae hold for $\pa_r u$ and $\pa_r^2 u$.

In the following estimates, as before we shall use the change of variables: $s=(1-r)/\varepsilon, t=\rho/\varepsilon$ and $\tau=t/(1+s)$.

We first handle higher order terms. By \eqref{test fcn:n=4} we estimate
\begin{align}\label{est:|Du|_L^2}
&\int_{\mathbb{S}^{3}}|\nabla u|_{\mathbb{S}^{3}}^{2}\ud\mu_{\Sp^3}\no\\
=&\int_{\mathbb{S}^{3}}\bigg[(\pa_\rho u)^2+\sin^{-2} \rho |\nabla u|_{\Sp^2}^2\bigg]\ud\mu_{\Sp^3}\no\\
 =&O(1)+O(1)\int_0^\delta \frac{\rho^{2}}{(\varepsilon^{2}+\rho^{2})^{2}}\sin^2\rho\ud \rho\no\\
=&O(1).
\end{align}

Based on the above calculations, we have
\begin{align}\label{extra term-1}
 & \int_{\mathcal{A}_{\delta}(x_{i})}\left|\bigg(\frac{2\cot\rho}{r^{2}}-\frac{2}{\rho}\bigg) \partial_\rho u\right|^{2} \ud x\no\\
\lesssim & \int_{\mathcal{A}_{\delta}(x_{i})}\big(O(\rho^{2})+O((1-r)^{2})\big)(\partial_\rho u)^{2}\ud x\no\\
\lesssim & \int_{\mathcal{A}_{\delta}(x_{i})}\big(O(\rho^{2})+O((1-r)^{2})\big)(\pa_\rho \phi_\ve)^{2}\ud x\no\\
\lesssim & \int_{\mathcal{A}_{\delta}(x_{i})}\big(O(\rho^{2})+O((1-r)^{2})\big)\frac{\rho^{2}}{((\varepsilon+1-r)^{2}+\rho^{2})^{2}}\rho^{2}\ud\rho \ud r\no\\
\lesssim & \varepsilon^{4}\int_{0}^{\frac{\delta}{\varepsilon}}\int_{0}^{\frac{\delta}{\varepsilon}}(O(t^{2})+O(s^{2}))\frac{t^{4}}{((1+s)^{2}+t^{2})^{2}}\ud t\ud s=O(1)
\end{align}
and
\begin{align}\label{extra term-2}
 & \int_{\mathcal{A}_{\delta}(x_{i})}\bigg((\frac{1}{r^{2}}-1)\partial_\rho^{2} u\bigg)^{2}\ud x\no\\
 \lesssim &  \int_{\mathcal{A}_{\delta}(x_{i})}O((1-r)^{2}) ((\pa_\rho \phi_\ve)^{4}+|\partial_\rho^{2} \phi_\ve|^2)\rho^{2}\ud\rho \ud r\no\\
\lesssim & \int_{\mathcal{A}_{\delta}(x_{i})}O((1-r)^{2})\big(\frac{1}{(\varepsilon+1-r)^{2}+\rho^{2}}\big)^{2}\rho^{2}\ud\rho \ud r\no\\
\lesssim & \int_{\mathcal{A}_{\delta}(x_{i})}O((1-r)^{2})\big(\frac{1}{(\varepsilon+1-r)^{2}+\rho^{2}}\big)^{2}\rho^{2}\ud\rho \ud r\no\\
\lesssim & \varepsilon^{2}\int_{0}^{\frac{\delta}{\varepsilon}}\int_{0}^{\frac{\delta}{\varepsilon}}\frac{s^{2}t^{2}\ud t\ud s}{((1+s)^{2}+t^{2})^{2}}= O(1).
\end{align}

Next, we focus on the remaining main term. Collecting the above calculations together we obtain
\begin{align*}
\Delta_{1}u:= & \partial_r^{2} u+\partial_\rho^{2} u+\frac{2}{\rho}\partial_\rho u\\
= & \frac{\partial_r^{2}\phi_{\ve}+\partial_\rho^{2}\phi_{\ve}+\frac{2}{\rho}\partial_\rho \phi_\ve}{1+c_1\nu_{i}^{-1}\log\frac{1}{\varepsilon}\left((\varepsilon+1-r)^{2}+\rho^{2}\right)^{3}}\\
 & +3(\pa_r \phi_\ve)^{2}\frac{c_1\nu_{i}^{-1}\log\frac{1}{\varepsilon}\left((\varepsilon+1-r)^{2}+\rho^{2}\right)^{3}}{\left[1+c_1\nu_{i}^{-1}\log\frac{1}{\varepsilon}\left((\varepsilon+1-r)^{2}+\rho^{2}\right)^{3}\right]^{2}}\\
 & +3(\pa_\rho \phi_\ve)^{2}\frac{c_1\nu_{i}^{-1}\log\frac{1}{\varepsilon}\left((\varepsilon+1-r)^{2}+\rho^{2}\right)^{3}}{\left[1+c_1\nu_{i}^{-1}\log\frac{1}{\varepsilon}\left((\varepsilon+1-r)^{2}+\rho^{2}\right)^{3}\right]^{2}}\\
=: & \frac{\partial_r^{2}\phi_{1}+\partial_\rho^{2}\phi_{1}+\frac{2}{\rho}\partial_\rho \phi_1+\partial_r^{2}\phi_{2}+\partial_\rho^{2}\phi_{2}+\frac{2}{\rho}\partial_\rho \phi_2}{1+c_1\nu_i^{-1}\log\frac{1}{\varepsilon}\left((\varepsilon+1-r)^{2}+\rho^{2}\right)^{3}}+B_{1}+B_{2}\\
 =:& B_{0}+B_{1}+B_{2}+B_{3},
\end{align*}
where 
\begin{align*}
B_{0}=&-\frac{4}{(\varepsilon+1-r)^{2}+\rho^{2}}\frac{1}{1+c_1\nu_i^{-1}\log\frac{1}{\varepsilon}\left((\varepsilon+1-r)^{2}+\rho^{2}\right)^{3}},\\
B_{1} =&3(\pa_r \phi_\ve)^{2}\frac{c_1\nu_{i}^{-1}\log\frac{1}{\varepsilon}\left((\varepsilon+1-r)^{2}+\rho^{2}\right)^{3}}{\left[1+c_1\nu_{i}^{-1}\log\frac{1}{\varepsilon}\left((\varepsilon+1-r)^{2}+\rho^{2}\right)^{3}\right]^{2}},\\
B_{2} 
  =&3(\pa_\rho \phi_\ve)^{2}\frac{c_1\nu_{i}^{-1}\log\frac{1}{\varepsilon}\left((\varepsilon+1-r)^{2}+\rho^{2}\right)^{3}}{\left[1+c_1\nu_{i}^{-1}\log\frac{1}{\varepsilon}\left((\varepsilon+1-r)^{2}+\rho^{2}\right)^{3}\right]^{2}},\\
B_{3}=&-\frac{8\varepsilon(\varepsilon+1-r)}{((\varepsilon+1-r)^{2}+\rho^{2})^{2}}\frac{1}{1+c_1\nu_i^{-1}\log\frac{1}{\varepsilon}\left((\varepsilon+1-r)^{2}+\rho^{2}\right)^{3}}.
\end{align*}

On one hand, we obtain an upper bound of
\begin{align*}
 &\int_{\mathcal{A}_{\delta}(x_i)}B_{0}^{2}\ud x\\
= & 16\int_{\mathcal{A}_{\delta}(x_{i})}\frac{1}{((\varepsilon+1-r)^{2}+\rho^{2})^{2}}\frac{1}{\left(1+c_1\nu_{i}^{-1}\log\frac{1}{\varepsilon}\left((\varepsilon+1-r)^{2}+\rho^{2}\right)^{3}\right)^{2}}\ud x\\
\leqslant & 16\left|\Sp^2\right| \int_{1-\delta}^{1}\int_{0}^{\delta}\frac{\rho^{2}}{\left(\left(\ve+1-r\right)^{2}+\rho^{2}\right)^{2}}\ud\rho \ud r\\
= & 16\left|\Sp^2\right| \int_{0}^{\frac{\delta}{\varepsilon}}\ud s\int_{0}^{\frac{\delta}{\varepsilon}}\frac{t^{2}}{\left((1+s)^{2}+t^{2}\right)^{2}}\ud t\\
= & 16\left|\Sp^2\right|\frac{\pi}{4}\log\frac{1}{\varepsilon}+O(1).
\end{align*}
On the other hand, we have
\begin{align*}
 & \int_{\mathcal{A}_{\delta}(x_i)}\left(\frac{1}{(\varepsilon+1-r)^{2}+\rho^{2}}\right)^{2}\frac{1}{\left(1+c_1\nu_{i}^{-1}\log\frac{1}{\varepsilon}\left((\varepsilon+1-r)^{2}+\rho^{2}\right)^{3}\right)^{2}}\ud x\\
= &\left|\Sp^2\right| \int_{0}^{\frac{\delta}{\varepsilon}}\ud s\int_{0}^{\frac{\delta}{\varepsilon}}\frac{t^2}{((1+s)^{2}+t^{2})^{2}}\frac{1}{\left(1+c_1\nu_{i}^{-1}\log\frac{1}{\varepsilon}\varepsilon^{6}\left((1+s)^{2}+t^{2}\right)^{3}\right)^{2}}\ud t+O(1)\\
\geqslant & \left|\Sp^2\right| \int_{0}^{\varepsilon^{-1}\frac{1}{(\log\frac{1}{\varepsilon})^{\frac{1}{3}}}}\frac{1}{1+s}\ud s\int_{0}^{\varepsilon^{-1}\frac{1}{(\log\frac{1}{\varepsilon})^{\frac{1}{3}}}\frac{1}{1+s}}\frac{\tau^{2}}{(1+\tau^{2})^{2}}\bigg(1+O(\frac{1}{\log\frac{1}{\varepsilon}})\bigg)\ud\tau+O(1)\\
= & \left|\Sp^2\right|\bigg(1+O(\frac{1}{\log\frac{1}{\varepsilon}})\bigg)\int_{0}^{\varepsilon^{-1}\frac{1}{(\log\frac{1}{\varepsilon})^{\frac{1}{3}}}}\frac{1}{1+s}\bigg(\int_{0}^{\infty}\frac{\tau^{2}}{(1+\tau^{2})^{2}}\ud\tau+O(\varepsilon(\log\frac{1}{\varepsilon})^{\frac{1}{3}}(1+s))\bigg)\ud s\\
&+O(1)\\
= & \frac{\pi}{4}\left|\Sp^2\right|\log\frac{1}{\varepsilon}+O(\log\log\frac{1}{\varepsilon}).
\end{align*}
Hence, we know 
\[
\int_{\mathcal{A}_{\delta}(x_i)}B_{0}^{2}\ud x=4 \pi \left|\Sp^2\right|\log\frac{1}{\varepsilon}+O(\log\log\frac{1}{\varepsilon}).
\]

Notice that
\begin{align*}
&B_{1}+B_{2}\\
=&\frac{3c_1\nu_{i}^{-1}\log\frac{1}{\varepsilon}}{\left(1+c_1\nu_{i}^{-1}\log\frac{1}{\varepsilon}\left(\varepsilon^{2}+\gamma(1-r)^{2}+\rho^{2}\right)^{3}\right)^{2}}\\
&\cdot \bigg[(\varepsilon+1-r)^{2}+\rho^{2})(4(1-r)^{2}+4\rho^{2})+16\varepsilon^{2}(1-r)^{2}\\
&\qquad+16\varepsilon(1-r)^2(\varepsilon+1-r)+16\varepsilon\rho^{2}(1-r)\bigg].
\end{align*}
Now let us consider
\begin{align*}
 & \int_{\mathcal{A}_{\delta}(x_{i})}(B_{1}+B_{2})^{2}\ud x\\
\lesssim & O((\log\frac{1}{\varepsilon})^{2})\bigg[\int_{\mathcal{A}_{\delta}(x_{i})}\frac{(\varepsilon+1-r)^{2}+\rho^{2})^{2}(4(1-r)^{2}+4\rho^{2})^{2}}{\left(1+c_1\nu_{i}^{-1}\log\frac{1}{\varepsilon}\left((\varepsilon+1-r)^{2}+\rho^{2}\right)^{3}\right)^{4}}\ud x\\
 & +\int_{\mathcal{A}_{\delta}(x_{i})}\frac{\varepsilon^{4}(1-r)^{4}}{\left(1+c_1\nu_{i}^{-1}\log\frac{1}{\varepsilon}\left((\varepsilon+1-r)^{2}+\rho^{2}\right)^{3}\right)^{4}}\ud x\\
 & +\int_{\mathcal{A}_{\delta}(x_{i})}\frac{\varepsilon^{2}(1-r)^{4}(\varepsilon+1-r)^{2}}{\left(1+c_1\nu_{i}^{-1}\log\frac{1}{\varepsilon}\left((\varepsilon+1-r)^{2}+\rho^{2}\right)^{3}\right)^{4}}\ud x\\
 & +\int_{\mathcal{A}_{\delta}(x_{i})}\frac{\varepsilon^{2}\rho^{4}(1-r)^{2}}{\left(1+c_1\nu_{i}^{-1}\log\frac{1}{\varepsilon}\left((\varepsilon+1-r)^{2}+\rho^{2}\right)^{3}\right)^{4}}\ud x\bigg]\\
:=& O\big(\log\frac{1}{\varepsilon}\big)^{2}|\Sp^2|\bigg(I_{1}+I_{2}+I_{3}+I_{4}\bigg),
\end{align*}
where 
\begin{align*}
I_{1}=&\int_{1-\delta}^1 \int_0^\delta \frac{((\varepsilon+1-r)^{2}+\rho^{2})^{2}(4(1-r)^{2}+4\rho^{2})^{2}}{\left(1+c_1\nu_{i}^{-1}\log\frac{1}{\varepsilon}\left((\varepsilon+1-r)^{2}+\rho^{2}\right)^{3}\right)^{4}}\rho^{2}\ud\rho \ud r,\\
I_{2}=&\int_{1-\delta}^1 \int_0^\delta\frac{\varepsilon^{4}(1-r)^{4}}{\left(1+c_1\nu_{i}^{-1}\log\frac{1}{\varepsilon}\left((\varepsilon+1-r)^{2}+\rho^{2}\right)^{3}\right)^{4}}\rho^{2}\ud\rho \ud r,\\
I_{3}=&\int_{1-\delta}^1 \int_0^\delta\frac{\varepsilon^{2}(1-r)^{4}(\varepsilon+1-r)^{2}}{\left(1+c_1\nu_{i}^{-1}\log\frac{1}{\varepsilon}\left((\varepsilon+1-r)^{2}+\rho^{2}\right)^{3}\right)^{4}}\rho^{2}\ud\rho \ud r,\\
I_{4}=&\int_{1-\delta}^1 \int_0^\delta\frac{\varepsilon^{2}\rho^{4}(1-r)^{2}}{\left(1+c_1\nu_{i}^{-1}\log\frac{1}{\varepsilon}\left((\varepsilon+1-r)^{2}+\rho^{2}\right)^{3}\right)^{4}}\rho^{2}\ud\rho \ud r.
\end{align*}

We shall estimate $I_i$ term by term. Let 
$$\gamma=\frac{1}{\varepsilon(\log\frac{1}{\varepsilon})^{b}},$$
where $b \in \R_+$ is to be determined later, then
\begin{align*}
I_{1}= &\int_{1-\delta}^1 \int_0^\delta\frac{((\varepsilon+1-r)^{2}+\rho^{2})^{2}(4(1-r)^{2}+4\rho^{2})^{2}}{\left(1+c_1\nu_{i}^{-1}\log\frac{1}{\varepsilon}\left((\varepsilon+1-r)^{2}+\rho^{2}\right)^{3}\right)^{4}}\rho^{2}\ud\rho \ud r\\
\lesssim & \varepsilon^{12}\bigg[\int_{0}^{\gamma}\ud s\int_{0}^{\gamma}((1+s)^{2}+t^{2})^{2}(s^{2}+t^{2})^{2}t^{2}\ud t\\
 & \qquad+\int_{0}^{\gamma}\ud s\int_{\gamma}^{\frac{\delta}{\varepsilon}}\frac{((1+s)^{2}+t^{2})^{2}(s^{2}+t^{2})^{2}t^{2}}{(\log\frac{1}{\varepsilon})^{4}\varepsilon^{24}\left((1+s)^{2}+t^{2}\right)^{12}}\ud t\\
 &\qquad +\int_{\gamma}^{\frac{\delta}{\varepsilon}}\ud s\int_{0}^{\frac{\delta}{\varepsilon}}\frac{((1+s)^{2}+t^{2})^{2}(s^{2}+t^{2})^{2}t^{2}}{(\log\frac{1}{\varepsilon})^{4}\varepsilon^{24}\left((1+s)^{2}+t^{2}\right)^{12}}\ud t\bigg]\\
\lesssim & \varepsilon^{12}\big(\gamma^{12}+\frac{(\log\frac{1}{\varepsilon})^{12b-4}}{\varepsilon^{12}}\big)\\
\lesssim & (\log\frac{1}{\varepsilon})^{12b-4},
\end{align*}
where we have used the following estimates:
\begin{align*}
 & \int_{0}^{\gamma}\ud s\int_{\gamma}^{\frac{\delta}{\varepsilon}}\frac{((1+s)^{2}+t^{2})^{2}(s^{2}+t^{2})^{2}t^{2}}{(\log\frac{1}{\varepsilon})^{4}\varepsilon^{24}\left((1+s)^{2}+t^{2}\right)^{12}}\ud t\\
\lesssim & \frac{1}{(\log\frac{1}{\varepsilon})^{4}\varepsilon^{24}}\int_{0}^{\gamma}\frac{\ud s}{(1+s)^{13}}\int_{\frac{\gamma}{1+s}}^{\frac{\delta}{\varepsilon(1+s)}}\frac{\tau^{2}}{(1+\tau^{2})^{8}}\ud\tau\\
\lesssim & \frac{\gamma^{-12}}{(\log\frac{1}{\varepsilon})^{4}\varepsilon^{24}}\lesssim \frac{(\log\frac{1}{\varepsilon})^{12b-4}}{\varepsilon^{12}}
\end{align*}
and 
\begin{align*}
 & \int_{\gamma}^{\frac{\delta}{\varepsilon}}\ud s\int_{0}^{\frac{\delta}{\varepsilon}}\frac{((1+s)^{2}+t^{2})^{2}(s^{2}+t^{2})^{2}t^{2}}{(\log\frac{1}{\varepsilon})^{4}\varepsilon^{24}\left((1+s)^{2}+t^{2}\right)^{12}}\ud t\\
\lesssim & \frac{1}{(\log\frac{1}{\varepsilon})^{4}\varepsilon^{24}}\int_{\gamma}^{\frac{\delta}{\varepsilon}}\frac{\ud s}{(1+s)^{13}}\int_{0}^{\frac{\delta}{\varepsilon(1+s)}}\frac{\tau^{2}}{(1+\tau^{2})^{8}}\ud\tau\\
\lesssim & \frac{1}{(\log\frac{1}{\varepsilon})^{4}\varepsilon^{24}}\int_{\gamma}^{\frac{\delta}{\varepsilon}}(1+s)^{-13}\ud s\\
\lesssim & \frac{\gamma^{-12}}{(\log\frac{1}{\varepsilon})^{4}\varepsilon^{24}}\lesssim \frac{(\log\frac{1}{\varepsilon})^{12b-4}}{\varepsilon^{12}}.
\end{align*}
Hence, we may take $12b-2<1$, i.e., $b \in (0,1/4)$ such that 
\[
(\log\frac{1}{\varepsilon}\big)^{2}I_{1}=o(\log\frac{1}{\varepsilon}).
\]

Let us deal with the other three terms:
\begin{align*}
 I_{2}=&\int_{1-\delta}^1 \int_0^\delta\frac{\varepsilon^{4}(1-r)^{4}}{\left(1+c_1\nu_{i}^{-1}\log\frac{1}{\varepsilon}\left((\varepsilon+1-r)^{2}+\rho^{2}\right)^{3}\right)^{4}}\rho^{2}\ud\rho \ud r\\
 =&16\varepsilon^{12}\int_{0}^{\frac{\delta}{\varepsilon}}\int_{0}^{\frac{\delta}{\varepsilon}}\frac{s^{4}t^{2}\ud t\ud s}{(1+c_{1}\nu_{i}^{-1}\varepsilon^{6}\log\frac{1}{\varepsilon}((1+s)^{2}+t^{2})^{3})^{4}}\\
  \leqslant& 16\varepsilon^{12}\int_{0}^{\frac{\delta}{\varepsilon}}\int_{0}^{\frac{\delta}{\varepsilon}}s^{4}t^{2}\ud t\ud s \lesssim \varepsilon^{4}
\end{align*}

\begin{align*}
I_{3} & =\int_{1-\delta}^1 \int_0^\delta\frac{\varepsilon^{2}(1-r)^{4}(\varepsilon+1-r)^{2}}{\left(1+c_1\nu_{i}^{-1}\log\frac{1}{\varepsilon}\left((\varepsilon+1-r)^{2}+\rho^{2}\right)^{3}\right)^{4}}\rho^{2}\ud\rho \ud r\\
 & \lesssim\varepsilon^{12}\int_{0}^{\frac{\delta}{\varepsilon}}\int_{0}^{\frac{\delta}{\varepsilon}}(1+s)^{2}s^{4}t^{2}\ud t\ud s \lesssim\varepsilon^{2}
\end{align*}
and
\begin{align*}
I_{4} & =\int_{1-\delta}^1 \int_0^\delta\frac{\varepsilon^{2}\rho^{4}(1-r)^{2}}{\left(1+c_1\nu_{i}^{-1}\log\frac{1}{\varepsilon}\left((\varepsilon+1-r)^{2}+\rho^{2}\right)^{3}\right)^{4}}\rho^{2}\ud\rho \ud r\\
 & \lesssim\varepsilon^{12}\int_{0}^{\frac{\delta}{\varepsilon}}\int_{0}^{\frac{\delta}{\varepsilon}}t^{6}s^{2}\ud t\ud s\lesssim\varepsilon^{2}.
\end{align*}

Hence, putting these facts together we obtain 
\begin{align*}
 & \int_{\mathcal{A}_{\delta}(x_{i})}(B_{1}+B_{2})^{2}\ud x\\
 & \lesssim\big(\log\frac{1}{\varepsilon}\big)^{2}\bigg(I_{1}+I_{2}+I_{3}+I_{4}\bigg)\\
 & \lesssim(\log\frac{1}{\varepsilon})^{12b-2}=o(\log\frac{1}{\varepsilon}),
\end{align*}
where $b \in (0,1/4)$.

We turn to show
\begin{align*}
 & \int_{\mathcal{A}_{\delta}(x_{i})}B_{3}^{2}\ud x\\
= & \int_{\mathcal{A}_{\delta}(x_{i})}\bigg(\frac{8\varepsilon(\varepsilon+1-r)}{((\varepsilon+1-r)^{2}+\rho^{2})^{2}}\frac{1}{1+c_1\nu_i^{-1}\log\frac{1}{\varepsilon}\left((\varepsilon+1-r)^{2}+\rho^{2}\right)^{3}}\bigg)^{2}\ud x\\
\lesssim &|\Sp^2|\int_{1-\delta}^1 \int_0^\delta\bigg[\frac{8\varepsilon(\varepsilon+1-r)}{((\varepsilon+1-r)^{2}+\rho^{2})^{2}}\bigg]^{2}\rho^{2}\ud\rho \ud r\\
\lesssim & \int_{0}^{\frac{\delta}{\varepsilon}}\int_{0}^{\frac{\delta}{\varepsilon}}\frac{(1+s)^{2}t^{2}}{((1+s)^{2}+t^{2})^{4}}\ud t\ud s\\
\lesssim & \int_{0}^{\frac{\delta}{\varepsilon}}\frac{\ud s}{(1+s)^{3}}\int_{0}^{\frac{\delta}{\varepsilon(1+s)}}\frac{\tau^{2}}{(1+\tau^{2})^{4}}\ud\tau \\
\lesssim & \int_{0}^{\frac{\delta}{\varepsilon}}\frac{1}{(1+s)^{3}}ds=O(1).
\end{align*}

Consequently, it follows from Cauchy inequality that 
\begin{equation}\label{main term}
\begin{split}
 & \int_{\mathcal{A}_{\delta}(x_{i})}(\Delta_{1}u)^{2}\ud x=4\pi \left|\Sp^{2}\right|\log\frac{1}{\varepsilon}+o(\log\frac{1}{\varepsilon}).
\end{split}
\end{equation}
By \eqref{extra term-1}, \eqref{extra term-2}, \eqref{main term} and summing $i$ from $1$ to $N_m$, we obtain
\begin{align}\label{est:Delta u_L^2}
  \int_{\mathbb{B}^4}(\Delta u)^{2}\ud x=4\pi\left|\Sp^{2}\right|N_m \log\frac{1}{\varepsilon}+o(\log\frac{1}{\varepsilon}).
\end{align}

Therefore, from the assumption that
\[
b+a\left(\int_{\mathbb{B}^4}(\Delta u)^{2}\ud x+2\int_{\mathbb{S}^{3}}|\nabla u|_{\mathbb{S}^{3}}^{2}\ud \mu_{\Sp^3}\right)\geqslant\log\left(\frac{1}{2 \pi^2}\int_{\mathbb{S}^{3}}e^{3u}\ud \mu_{\Sp^3}\right)-3\bar u,
\]
putting the above estimates \eqref{area-ln}, \eqref{average-u}, \eqref{est:|Du|_L^2} and \eqref{est:Delta u_L^2} together, we conclude that
\[
a\geqslant \frac{3}{4\pi |\mathbb{S}^{2}| N_{m}}=\frac{3}{16\pi^2 N_m}.
\]
This completes our construction.
\end{proof}

\appendix

\section{An example on the sharpness of Lebedev-Milin inequality under constraints}\label{Append}

The purpose of this appendix is to demonstrate our idea in dimension four through an example for Lebedev-Milin inequality in higher order moments case, which will be presented in a concise manner, as the same idea prevails throughout the paper.

For clarity, we restate  Widom \cite{Widom} inequality, a generalization of Lebedev-Milin inequality: If $u\in H^1(\mathbb{B}^2)$ satisfies $\int_{\Sp^1}e^u p \ud\mu_{\Sp^1}=0$ for all $p \in \mathring{\mathcal{P}}_m$, then
$$\log\left(\frac{1}{2\pi}\int_{\Sp^1}e^{u-\bar u} \ud \mu_{\Sp^1}\right)\leqslant \frac{1}{4\pi N_m(\Sp^1)}\int_{\mathbb{B}^2}|\nabla u|^2 \ud x,$$
 where $\bar u=(2\pi)^{-1}\int_{\Sp^1}e^u  \ud\mu_{\Sp^1}$ and $N_m(\Sp^1)=m+1$ as shown in \cite{Chang-Hang}. For $m=1$, the above inequality was first proved by Osgood-Phillips-Sarnak \cite{OPS}. Furthermore, if we relax $m$ to define $\mathring{\mathcal{P}}_0=\emptyset$ and $N_0(\Sp^1)=1$, then the above inequality with $m=0$ reduces to Lebedev-Milin inequality. An alternative proof of the above Widom's inequality is also available in \cite[Section 6]{Chang-Hang}.
 
 The sharpness of the above inequality can be shown in the following way: If there exists $a\in \R_+$ such that for all $u\in H^1(\mathbb{B}^2)$ satisfying $\int_{\Sp^1}e^u p \ud\mu_{\Sp^1}=0$ for all $p \in \mathring{\mathcal{P}}_m$, we have
 $$a\int_{\mathbb{B}^2}|\nabla u|^{2}\ud x\geqslant\log\left(\frac{1}{2\pi}\int_{\mathbb{S}^{1}}e^{u-\bar u}\ud x\right),
$$
then 
$$a \geqslant \frac{1}{4\pi N_m}.$$

For brevity, we let $N_m=N_m(\Sp^1)$.

 In dimension two, for each $x_i:=(\cos \theta_i,\sin \theta_i)\in \Sp^1$ we define
 $$\mathcal{A}_\delta (x_i):=\left\{x=r\xi \in \mathbb{B}^n; \xi \in \Sp^{n-1}, 1-r< \delta, ~~ |\theta-\theta_i|< \delta \right\}.$$
 Near $x_i$, we use the polar coordinates
 $$|\ud x|^2=\ud r^2+r^2 \ud \theta^2$$
 for $x=(r,\theta)$ with $\xi=(\cos \theta,\sin \theta)$.
 
 Define
\[
\phi_{\varepsilon,i}(r,\theta)=-\log\left(\varepsilon^{2}+(1-r)^{2}+(\theta-\theta_i)^{2}\right).
\]

Let $\rho=\theta-\theta_i$, then
\begin{align*}
&\sum_{i=1}^{N_m}\int_{\mathbb{S}^{1}}\chi_{i}e^{\phi_{\varepsilon,i}+\log\nu_{i}}\ud\theta\\
=&\int_{-\delta}^\delta \frac{\ud \rho}{\varepsilon^2+\rho^2}+O(1)\int_\delta^{2\delta} \frac{\ud \rho}{\varepsilon^2+\rho^2}\\
=&2\varepsilon^{-1}\int_{0}^{\delta/\varepsilon}\frac{\ud t}{1+t^2}+O(1)\\
=&\pi\varepsilon^{-1}+O(1)
\end{align*}
and for all $p \in \mathring{\mathcal{P}}_m$,
\begin{align*}
 & \sum_{i=1}^{N_m}\int_{\mathbb{S}^{1}}\chi_{i}e^{\phi_{\varepsilon,i}+\log\nu_{i}}p(x)\ud \theta\\
  =&\sum_{i=1}^{N_m}\nu_i\int_{-\delta}^\delta\chi_{i}e^{\phi_{\varepsilon,i}}O(|\theta-\theta_i|^{2})\ud \theta+O(1)\\
 =&O(1).
\end{align*}


We can choose a test function as  
\begin{align*}
e^{u} =\sum_{i=1}^{N_m}\chi_{i}e^{\phi_{\varepsilon,i}+\log\nu_{i}}+c_{1}+\sum_{j=1}^{2m}\beta_{j}\psi_{j}
\end{align*}
As before, we can find $\psi_j \in C_c^\infty\left([0,2\pi]\setminus \cup_{i=1}^{m+1}(\theta_i-2\delta,\theta_i+2\delta)\right)$ and $c_1\in \R_+$ such that
$$\int_{\mathbb{S}^1}e^u p_k \ud\mu_{\mathbb{S}^1}=0 \qquad 1 \leqslant k \leqslant 2m$$
and
$$c_{1}+\sum_{j=1}^{2m}\beta_{j}\psi_{j}>1,$$
 then it in turn implies $\beta_{j}=O(1)$.

Now we check that 
\begin{align*}
\bar{u} =&\frac{1}{2\pi}\int_{\mathbb{S}^{1}}u \ud \mu_{\Sp^1}\\
  =&O(1)+\frac{1}{2\pi}\int_0^\delta \log\left(\frac{1}{\varepsilon^{2}+(\theta-\theta_i)^{2}}+O(1)\right)\ud \theta\\
 =& O(1)+\frac{1}{2\pi}\int_{-\delta}^\delta\log\left(\frac{1}{\varepsilon^{2}+\rho^{2}}\right)\ud\rho\\
 =& \frac{1}{\pi}\int_{0}^{\delta}\left(2\log\varepsilon^{-1}+\log\frac{1}{1+\left(\frac{\rho}{\varepsilon}\right)^{2}}\right)\ud\rho+O(1)\\
 =&\frac{1}{\pi}\left[2\delta\log\varepsilon^{-1}-\varepsilon\int_{0}^{\delta/\varepsilon}\log(1+t^2)\ud t\right]+O(1)\\
 =& \frac{1}{\pi}\left[2\delta\log\varepsilon^{-1}-\delta \log(1+\delta^2\ve^{-2})\right]+O(1)\\
 =&O(1).
\end{align*}

Observe that
\[
\partial_{r}\phi_{\ve,i}=\frac{2(1-r)}{\varepsilon^{2}+(1-r)^{2}+(\theta-\theta_i)^{2}}
\]
 and 
\[
\partial_{\theta}\phi_{\ve,i}=\frac{-2(\theta-\theta_i)}{\varepsilon^{2}+(1-r)^{2}+(\theta-\theta_i)^{2}},
\]
 then 
\begin{align*}
|\nabla\phi_{\ve,i}|^{2} =\frac{4\left((1-r)^{2}+(\theta-\theta_i)^{2}\right)}{\left(\varepsilon^{2}+(1-r)^{2}+(\theta-\theta_i)^{2}\right)^{2}}+\frac{4(\theta-\theta_i)^{2}}{\left(\varepsilon^{2}+(1-r)^{2}+(\theta-\theta_i)^{2}\right)^{2}}(\frac{1}{r^{2}}-1).
\end{align*}

$ $

Under the change of variables: $s=(1-r)/\varepsilon$ and $t=\rho/\varepsilon$, we obtain
\begin{align*}
 & \int_{A_{\delta}\left(x_{i}\right)}\frac{4\left((1-r)^{2}+(\theta-\theta_i)^{2}\right)}{\left(\varepsilon^{2}+(1-r)^{2}+(\theta-\theta_i)^{2}\right)^{2}}r\ud r\ud\theta\\
= & \int_0^\delta\int_{-\delta}^\delta\frac{4\left((1-r)^{2}+\rho^{2}\right)}{\left(\varepsilon^{2}+(1-r)^{2}+\rho^{2}\right)^{2}}\ud\rho\ud r+ \int_0^\delta\int_{-\delta}^\delta\frac{4\left((1-r)^{2}+\rho^{2}\right)}{\left(\varepsilon^{2}+(1-r)^{2}+\rho^{2}\right)^{2}}(r-1)\ud\rho\ud r\\
= & 2\int_{0}^{\delta/\varepsilon}\int_{0}^{\mathscr{\delta}/\varepsilon}\frac{4\left(s^{2}+t^{2}\right)}{\left(1+s^{2}+t^{2}\right)^{2}}\ud s\ud t+O(1)\\
= & 8\int_{0}^{\delta/\varepsilon}\int_{0}^{\frac{\delta}{\varepsilon\sqrt{1+s^2}}}\frac{s^{2}+\tau^{2}\left(1+s^{2}\right)}{\left(1+s^{2}\right)^{3/2}\left(1+\tau^{2}\right)^{2}}\ud\tau\ud s+O(1)\\
= &8\int_0^{\delta/\varepsilon}\frac{\ud s}{\sqrt{1+s^2}}\int_{0}^{\frac{\delta}{\varepsilon\sqrt{1+s^2}}}\frac{\ud \tau}{1+\tau^2}\\
&-8 \int_0^{ \delta/\varepsilon}\frac{\ud s}{(1+s^2)^{3/2}} \int_{0}^{\frac{\delta}{\varepsilon\sqrt{1+s^2}}}\frac{\ud \tau}{(1+\tau^2)^2}+O(1)\\
= & 4\pi\log\frac{1}{\varepsilon}+O(1)
\end{align*}
and
\begin{align*}
 & \int_{A_{\delta}\left(x_{i}\right)}\frac{4(\theta-\theta_i)^{2}}{\left(\varepsilon^{2}+(1-r)^{2}+(\theta-\theta_i)^{2}\right)^{2}}(\frac{1}{r^{2}}-1)\ud x\\
\lesssim & \varepsilon \int_{0}^{\frac{\delta}{\varepsilon}}\int_{0}^{\frac{\delta}{\varepsilon}}\frac{t^{2}s}{\left(1+s^{2}+t^{2}\right)^{2}}\ud s\ud t=O(1).
\end{align*}

Thus, we obtain
\[
\int_{\mathbb{B}^2}|\nabla u|^{2}\ud x=4\pi N_{m}\log\frac{1}{\varepsilon}+O(1).
\]

On the other hand, we have 
$$\int_{\mathbb{S}^{1}}e^{u}\ud \mu_{\Sp^1}  =\int_{-\delta}^\delta\frac{1}{\varepsilon^{2}+\rho^{2}}\ud\rho+O(1) =\pi\varepsilon^{-1}+O(1).$$
This gives
\[
\log\left(\frac{1}{2\pi}\int_{\mathbb{S}^{1}}e^{u}\ud \mu_{\Sp^1}\right)=\log\frac{1}{\varepsilon}+O(1).
\]

Hence, it follows from the assumption that
\[
a\int_{\mathbb{B}^2}|\nabla u|^{2}\ud x\geqslant\log\left(\frac{1}{2\pi}\int_{\mathbb{S}^{1}}e^{u}\ud \mu_{\Sp^1}\right)-\bar{u}.
\]
Putting these facts together, we conclude that
$$a \geqslant \frac{1}{4\pi N_m}.$$

\begin{remark}
The above strategy can  also provide an alternative of Chang-Hang's test function in \cite{Chang-Hang} as follows: we can replace the piecewise Lipschitz function $\phi_\ve(t)$ in Chang-Hang \cite[p.10]{Chang-Hang} by a global function
$$\phi_\ve(t)=-\log(\ve^2+t^2),$$
indeed, Chang-Hang's $\phi_\ve$ also occurred in Lions' example for Moser-Trudinger inequality in \cite[pp.195-199]{Lions1}.
By choosing constants $\beta_j, c_1$ and smooth cut-off functions $\chi_i,\psi_j$ as Chang-Hang's, we define a smooth test function
$$e^{2u}=\sum_{i=1}^N \chi_i e^{2\phi_\ve(\overline{xx_i})+\log \nu_i}+\sum_{j=1}^{m^2+2m}\beta_j \psi_j+c_1\log \frac{1}{\ve}.$$
We follow the same lines of Chang-Hang to give an example to show \emph{almost sharpness} of Moser-Trudinger-Onofri inequality under constraints.
\end{remark}

\begin{thebibliography}{99}

\bibitem{Ache-Chang}
A.  Ache and A. Chang, \textit{Sobolev trace inequalities of order four}, Duke Math. J. 166 (2017), no. 14, 2719-2748.
 
\bibitem{Aubin}
T. Aubin, \textit{\'{E}quations diff\'{e}rentielles non lin\'{e}aires et probl\`{e}me de Yamabe concernant la courbure scalaire}, J. Math. Pures Appl. (9) 55 (1976), no. 3, 269-296.

\bibitem{Aubin_book}
T. Aubin, \textit{Some nonlinear problems in Riemannian geometry}, Springer Monographs in Mathematics, Springer-Verlag, Berlin, 1998. 

\bibitem{Beckner}
W. Beckner, \textit{Sharp Sobolev inequalities on the sphere and the Moser-Trudinger inequality}, 
Ann. of Math. (2) 138 (1993), no. 1, 213-242. 

\bibitem{Bannai1}
E. Bannai, \textit{On some spherical t-designs}, J. Combin. Theory Ser. A 26 (1979), no. 2, 157-161.

\bibitem{Bannai2}
E. Bannai, \textit{On tight spherical designs}, J. Combin. Theory Ser. A 26 (1979), no. 1, 38-47. 

\bibitem{Bannai3}
E. Bannai and R. Damerell, \textit{Tight spherical designs. I}, J. Math. Soc. Japan 31 (1979), no. 1, 199-207.

\bibitem{Case}
J. Case, \textit{Boundary operators associated with the Paneitz operator}, Indiana Univ. Math. J. 67 (2018), no. 1, 293-327.

\bibitem{Case-Chang}
J. Case and A. Chang, \textit{On fractional GJMS operators},
Comm. Pure Appl. Math. 69 (2016), no. 6, 1017-1061. 

 
 \bibitem{Chang-Hang}
A. Chang and F. Hang, \textit{Improved Moser-Trudinger-Onofri inequality under constraints}, to appear in Comm. Pure Appl. Math., \href{https://arxiv.org/abs/1909.00431}{arXiv:1909.00431}.

\bibitem{Chang-Xu-Yang}
 A. Chang, X. Xu and P. Yang, \textit{A perturbation result for prescribing mean curvature}, Math. Ann. 310 (1998), no. 3, 473-496.
 
 \bibitem{Chang-Yang3}
 A. Chang and P. Yang, \textit{ Extremal metrics of zeta function determinants on 4-manifolds}, 
Ann. of Math. (2) 142 (1995), no. 1, 171-212. 

\bibitem{Dai-Xu}
F. Dai and Y. Xu, \textit{Approximation theory and harmonic analysis on spheres and balls}, Springer Monographs in Mathematics. Springer, New York, 2013. xviii+440 pp.

\bibitem{DGS}
P. Delsarte, J. Goethals and J. Seidel, \textit{Spherical codes and designs}, Geometriae Dedicata 6 (1977), no. 3, 363-388

\bibitem{escobar7}
 J. Escobar, \textit{Sharp constant in a Sobolev trace inequality}, Indiana Univ. Math. J. 37 (1988), no. 3, 687-698.
 
\bibitem{escobar4}
J. Escobar, \textit{The Yamabe problem on manifolds with boundary}, J. Differential Geom. 35 (1992), no. 1, 21-84.

\bibitem{Fefferman-Graham}
 C. Fefferman and C. Graham,  \textit{$Q$-curvature and Poincar\'e metrics}, Math. Res. Lett. 9 (2002), no. 2-3, 139-151.


\bibitem{GGS}
 F. Gazzola, H. Grunau and G. Sweers, \textit{Polyharmonic boundary value problems. Positivity preserving and nonlinear higher order elliptic equations in bounded domains}, Lecture Notes in Mathematics, 1991. Springer-Verlag, Berlin, 2010. xviii+423 pp. 
 
 \bibitem{Graham-Zworski}
 
C. Graham and M. Zworski, \textit{Scattering matrix in conformal geometry}, Invent. Math. 152 (2003), no. 1, 89-118.

\bibitem{Guo-Wang}
Q. Guo and X. Wang, \textit{Uniqueness results for positive harmonic functions on $\overline{\Bn}$ satisfying a nonlinear boundary condition}, Calc. Var. Partial Differential Equations 59 (2020), no. 5, Paper No. 146, 8 pp.
 
  \bibitem{Hang}
 F. Hang, \textit{A remark on the concentration compactness principle in critical dimension},  to appear in Comm. Pure Appl. Math., \href{https://arxiv.org/abs/2002.09870}{arXiv:2002.09870}.
 
 \bibitem{Hang-Wang}
 F. Hang and X. Wang, \textit{Improved Sobolev inequality under constraints}, to appear in Int. Math. Res. Not., \href{https://arxiv.org/abs/2010.10654}{arXiv:2010.10654}.
 
 \bibitem{Lebedev-Milin}
 N. Lebedev and I. Milin, \textit{ On the coefficients of certain classes of analytic functions}, (Russian) Mat. Sbornik N.S. 28(70), (1951), 359-400.
 

\bibitem{Li-Zhang}
Y. Y. Li and L. Zhang, \textit{Liouville-type theorems and Harnack-type inequalities for semilinear elliptic equations}, J. Anal. Math. 90 (2003), 27-87.

 \bibitem{li-zhu}
 Y. Y. Li and M. Zhu,  \textit{Uniqueness theorems through the method of moving spheres}, Duke Math. J. 80 2 (1995), 383-417.
 
 \bibitem{Lions1}
 P. Lions, \textit{ The concentration-compactness principle in the calculus of variations. The limit case. I}, Rev. Mat. Iberoamericana 1 (1985), no. 1, 145-201.
 
  \bibitem{Lions2}
 P. Lions, \textit{ The concentration-compactness principle in the calculus of variations. The limit case. II}, Rev. Mat. Iberoamericana 1 (1985), no. 2, 45-121.
 
\bibitem{Mysovskih}
 I. Mysovskih, \textit{A proof of minimality of the number of nodes of a cubature formula for a hypersphere}, (Russian) \v{Z}. Vy\v{c}isl. Mat i Mat. Fiz. 6 (1966), 621-630.
 
 \bibitem{OPS}
B. Osgood, R. Phillips and P. Sarnak, \textit{Extremals of determinants of Laplacians}, J. Funct. Anal. 80 (1988), no. 1, 148-211.
 
 \bibitem{Putterman}
 E. Putterman, \textit{Cubature formulas and Sobolev inequalities}, \href{https://arxiv.org/abs/2012.08109}{arXiv:2012.08109}.
 
  
  \bibitem{Stein-Weiss}
 E. Stein and G. Weiss, \textit{Introduction to Fourier analysis on Euclidean spaces}, 
Princeton Mathematical Series, No. 32. Princeton University Press, Princeton, N.J., 1971. x+297 pp. 

\bibitem{Struwe_book}
M. Struwe, \textit{Variational methods. Applications to nonlinear partial differential equations and Hamiltonian systems}, Fourth edition. Ergebnisse der Mathematik und ihrer Grenzgebiete. 3. Folge. A Series of Modern Surveys in Mathematics [Results in Mathematics and Related Areas. 3rd Series. A Series of Modern Surveys in Mathematics], 34. Springer-Verlag, Berlin, 2008. xx+302 pp.

\bibitem{Widom}
H. Widom, \textit{On an inequality of Osgood, Phillips and Sarnak}, Proc. Amer. Math. Soc. 102 (1988), no. 3, 773-774.
 
\end{thebibliography}
\end{document}